\documentclass[12pt]{amsart}

\usepackage[usenames]{color}
\input xypic

\usepackage[colorlinks=false, linkcolor=blue]{hyperref}
\makeatletter
\def\l@subsection{\@tocline{2}{0.25em}{2.3em}{0pt}{}}
\makeatother

\usepackage{verbatim}
\usepackage{url}
\usepackage{bm}
\usepackage{tikz}

\usepackage{graphicx}
\usepackage{color}
\usepackage{nicematrix}
\usepackage{dsfont}

\usepackage{amsmath}
\usepackage[centertags]{amsmath}
\usepackage{amsfonts}
\usepackage{amssymb}
\usepackage{amsthm}
\usepackage{newlfont}
\usepackage{url}
\usepackage{bm}
\usepackage{chngcntr}

\theoremstyle{plain}
\newtheorem{thm}{Theorem}
\newtheorem{cor}[thm]{Corollary}

\newtheorem{lem}[thm]{Lemma}
\newtheorem{lem*}[thm]{Lemma}
\newtheorem{prop}[thm]{Proposition}
\newtheorem{question}[thm]{Question}

\theoremstyle{definition}
\newtheorem{dfn}{Definition}
\theoremstyle{remark}
\newtheorem{rem}{Remark}
\newtheorem{rem*}{Remark}
\newtheorem{ex}[rem]{Example}
\newtheorem{example}[rem]{Example}

\numberwithin{rem}{section} 
\numberwithin{dfn}{section} 
\numberwithin{equation}{section} 
\numberwithin{thm}{section} 

\def\!{\operatorname{!}}

\def\C{\mathbb C}
\def\F{\mathbb F}
\def\G{\mathbb G}

\def\1{\bold 1}

\def\diag{\operatorname{diag}}
\def\deg{\operatorname{deg}}

\def\Hom{\operatorname{Hom}}
\def\GL{\operatorname{GL}}

\def\Ext{\operatorname{Ext}}

\def\Lie{\operatorname{Lie}}

\def\rk{\operatorname{rk}}

\def\deg{\operatorname{deg}}

\def\End{\operatorname{End}}

\def\coker{\operatorname{coker}}

\def\mod{\operatorname{mod}}

\usepackage{amssymb}
\usepackage{enumitem}

\usepackage{graphicx}
\usepackage{xy}

\makeatletter
\@namedef{subjclassname@2010}{%
	\textup{2010} Mathematics Subject Classification}
\makeatother

\usepackage[T1]{fontenc}

\newtheorem{theorem}{Theorem}[section]

\theoremstyle{definition}
\newtheorem{definition}[theorem]{Definition}
\newtheorem{conjecture}[theorem]{Conjecture}
\newtheorem{remark}[theorem]{Remark}

\numberwithin{equation}{section}

\newcommand{\Der}{\mathrm{Der} }

\newcommand{\Derin}{\Der_{in}}

\newcommand{\lra}{\longrightarrow}

\newcommand{\podwzorem}[2]{\underbrace{#1}\limits_{#2}}

\newcommand{\tm}{$\mathbf{t}$-}

\newcommand{\uplra}[1]{\stackrel{#1}{\lra}}

\newcommand{\Mat}{\mathrm{Mat}}

\newcommand{\Exp}{\mathrm{Exp}}

\newcommand{\wlozenie}{\small\left[\begin{array}{c}
		0\\1
	\end{array}\right]}
\newcommand{\rzut}{\small\left[\begin{array}{cc}
		1&0
	\end{array}\right]}

\newcommand{\block}[1]{\big[\big|#1\big|\big]}

\begin{document}

	
	\baselineskip=17pt
	
	
	\title[]{Homological Methods in the Generalization of Drinfeld Modules }

	\author{Dawid E. K\k{e}dzierski, Piotr Kraso{\'n}$^\star$}
	
	\thanks{$^\star$Corresponding author.}
	
	\date{\today }

	\address{  Institute of Mathematics, Department of Exact and Natural Sciences, University of Szczecin, ul. Wielkopolska 15, 70-451 Szczecin, Poland 
	}
	\email{dawidedmundkedzierski@gmail.com}
	
	\address{Ariel University, Ariel 40700, Israel}
	\email{ piotrkras26@gmail.com}

	\thanks{}

	\maketitle
	
	\hyphenation{mo-du-les}
	\newcommand{\functor}{\mathrm{Rep}}

	\begin{abstract}
We introduce and study a natural class of Anderson \tm mo\-dules, called
triangular \tm modules, characterized by having Drinfeld modules as their
$\tau$-composition factors.  They form a homologically meaningful
generalization of Drinfeld modules and exhibit rich arithmetic structure.\smallskip

We establish criteria for  purity, strict and almost strict, and develop a
reduction procedure that lowers the degrees of the defining biderivations.
As a consequence, every almost strictly pure triangular \tm module becomes
strictly pure after a finite base extension.\smallskip

We then investigate morphisms and isogenies between triangular \tm mo\-dules,
provide a characterization of triangular isogenies, and describe the algebra of endomorphisms, including a criterion for commutativity.  On the
analytic side, we show that all triangular \tm modules are uniformizable and
establish finiteness and purity criteria with consequences for Taelman's
conjecture.\smallskip

Finally, we develop a duality theory for triangular \tm modules and their
biderivations, proving compatibility with $\tau$-composition series and
establishing analogues of  the Weil--Barsotti
formula and the Cartier--Nishi theorem for a special class of triangular \tm modules.

.  \smallskip

	\end{abstract}

\maketitle

\tableofcontents

\section{Introduction}
Drinfeld modules and their higher-dimensional analogues, the Anderson 
\tm modules, play a central role in the arithmetic of global function fields.  
This is due to their numerous applications, for instance in class field theory, 
the theory of automorphic forms, the Langlands program, and Diophantine 
geometry (see \cite{g96}, \cite{th04}, \cite{bp20}).  
After adjoining the zero object, the category of \tm modules becomes an  
${\mathbb F}_q[t]$-linear additive category, although it is not abelian.  
Despite many analogies with the category of abelian varieties, one of the key 
differences is that the category of \tm modules is far from semisimple, making 
it natural to study the groups $\Ext^1(\Phi,\Psi)$.  The investigation of these 
extension groups was initiated in \cite{pr} and subsequently developed in 
\cite{kk04}, \cite{gkk}, and \cite{kk24}.  A central technical tool in these 
works is to interpret extensions as appropriate biderivations.  An interesting 
and subtle question is when ${\Ext}^1(\Phi,\Psi)$ itself carries a natural 
\tm module structure.

In \cite{gkk} we introduced an algorithm, called the \emph{\tm module 
reduction} algorithm, which allows one to determine \tm module structures on 
$\Ext^1(\Phi,\Psi)$ for wide classes of pairs of \tm modules.  In that paper the 
notion of a $\tau$-composition series was developed, and it was proved in 
\cite[Theorem~5.2, Corollary~5.3]{gkk} that, under suitable hypotheses, the 
\tm module reduction algorithm applies to \tm modules possessing specific 
$\tau$-composition series.  Among these, the \tm modules whose 
$\tau$-composition factors are Drinfeld modules deserve special attention.  
In a suitable basis, such \tm modules have  triangular matrix form, and hence we 
refer to them as \emph{triangular \tm modules}.  They form a new and rich class 
in the category of \tm modules and, from the homological viewpoint, they 
generalize the class of Drinfeld modules.  This leads to the following question:

\begin{question}
Which properties of Drinfeld modules are inherited by triangular \tm modules?
\end{question}

The goal of this paper is to answer this question, as well as related questions 
for the corresponding \tm motives and \tm comotives associated with triangular 
\tm modules.

\smallskip

In Section~2, we recall the basic definitions and results needed later in the 
paper.  In particular, we review both the analytic and the algebraic approaches 
to Anderson \tm modules, and discuss extensions of \tm modules via 
biderivations.  We also introduce the category of \tm motives, which is known 
to be anti-equivalent to the category of \tm modules (cf.\ 
Theorem~\ref{equiva}), as well as the dual notion of a \tm comotive.

\begin{remark}
Following \cite{gm25} we use the term \emph{\tm comotive} instead of the more 
traditional \emph{dual \tm motive}.  This terminology avoids potential 
confusion with the dual of a \tm motive in the usual categorical sense.
\end{remark}

\smallskip

In Section~3, we introduce triangular $t$-modules and develop
a reduction theory culminating in a reduced form theorem (see Theorem~\ref{prop:t-reduced_form}), which shows that any
triangular $t$-module acquires a form with biderivations of minimal degree
after a finite base extension of a controlled type. 
\smallskip

Section~4 is devoted to morphisms between triangular \tm modules.  
We introduce the notion of a triangular morphism, which behaves  
analogously to morphisms of Drinfeld modules; see 
Proposition~\ref{basictr}.  We then study triangular isogenies 
and give a characterization of when a triangular morphism is an isogeny 
(cf.\ Theorem~\ref{prop:triangular_isogenie}).  
We conclude this section by studying the algebra of endomorphisms and 
provide a sufficient condition for its commutativity.  
Proposition~\ref{prop:wyznaczanie_morfizmu_z_wyrazu_wolnego} concerns the 
effective determination of morphisms in generic characteristic.

\smallskip
In Section 5, we prove that the category of \tm modules is exact in the sense of Quillen \cite{Q}. In \cite{g24}, Q. Gazda proved that the category of \tm motives is exact. Since the category of abelian \tm modules is equivalent to the category of abelian \tm motives by Anderson's equivalence (cf. Theorem 2.1), it follows that the category of abelian \tm modules is exact. In this paper, however, we prove directly that the category of all \tm modules is exact. As a consequence, the category of triangular \tm modules is a full exact subcategory.

\smallskip

In Section~6, we develop the analytic theory of triangular \tm modules.  
In particular, we prove a general result (Theorem~\ref{thm62}) on the behavior of 
exponential functions in exact sequences.  As a consequence, every triangular 
\tm module is uniformizable (see Corollary \ref{cor:triangular_is_uniformizable}).

\smallskip

Section 7 gives necessary and sufficient conditions for triangular \tm modules to be strictly pure (see Proposition~\ref{prop:strictly_pure_characterisation}) and almost strictly pure (see Proposition~\ref{prop:almost_strictly_pure} and Remark \ref{rem:almost_stricly_pure}). We then show that triangular \tm modules are \tm finite and hence abelian. Theorem~\ref{finitness} provides a criterion for purity, and we conclude the section with corollaries concerning Taelman's conjecture for triangular \tm modules (see Corollaries~\ref{thm:taelmas_conjecture} and \ref{thm:taelmas_conjecture2}).

\smallskip

In Section~8, we recall the definition of the dual \tm module introduced in 
\cite{kk24} and study the duality of triangular 
\tm modules.  We prove that this duality is compatible with $\tau$-composition 
series and also establish the Weil--Barsotti formula 
(Theorem~\ref{thm:dual_dla_trangulowalnego_ogolnie}).  We then adapt duality to 
biderivations, defining the dual biderivation $\delta^\vee$ and its double 
dual ${\delta^\vee}^\vee$, and show that the assignments 
$\delta\mapsto\delta^\vee$ and $\delta^\vee\mapsto{\delta^\vee}^\vee$ are 
${\mathbb F}_q$-linear (Corollary~\ref{cor:dual_dla_biderwaycji}).  
Finally, we give explicit formulas for dual biderivations in the case of triangular 
\tm modules allowing for LD-biderivations.

\smallskip

Section~9 establishes an analog of the Cartier--Nishi theorem for triangular 
\tm modules allowing for LD-biderivations 
(Theorem~\ref{thm:dual_dual_case_1}).  
We also prove Corollaries~\ref{cor:Cartier_Nishi_1} and 
\ref{cor:Cartier_Nishi_2}, which describe  cases where the sequence of ranks 
of sub-quotients is strictly increasing or at least non-decreasing.  As in 
\cite{kk24}, we use the $\Ext$-definition of the dual \tm module.  
In \cite[Appendix]{kk24}  by refining the methods of \cite{ta}, it was shown,
that for strictly pure \tm modules without nilpotence the dual module of Taguchi is 
isomorphic to that given by the $\Ext$ functor.  Since triangular 
\tm modules are not strictly pure in general, our results provide a new class 
of \tm modules for which the dual module  exists and satisfies 
${\Phi^{\vee}}^{\vee} \cong \Phi$.

\smallskip

\section{Background}
Let $p$ be a prime number,  $\F_q$  a finite field with $q=p^s$ elements, and $A=\F_q[t]$ 
the ring of polynomials in one variable.
By $k$ we  denote the field of rational functions $\F_q(\theta)$ in a variable $\theta$. The  $\infty-$adic completion $k_\infty:=\F_q((1/\theta))$ can be given a normalized valuation $|\cdot|_\infty$ satisfying $|\theta|_\infty=q$. Finally, in order to ensure that the  field under consideration is algebraically closed and complete, one has to pass to $\C_\infty$ the completion with respect to $|\cdot|_\infty$ of the algebraic closure $\overline{k_\infty}.$ Then there exists a natural embedding $\iota: A\lra \C_\infty.$ 

\subsection{Drinfeld modules}
Recall that an $A-$submodule $\Lambda\subset\C_\infty$ is called an $A-$lattice if it is a finitely generated, projective $A-$module, discrete in the topology of  $\C_\infty$. For every  $A-$lattice $\Lambda$  there is  an entire function: $\exp_\Lambda:\C_{\infty}\lra\C_\infty$ defined by the formula:
$$\exp_\Lambda(z)=z\prod_{0\neq \lambda\in\Lambda}\Big(1-\dfrac z\lambda\Big).$$
The function  $\exp_\Lambda$ is $\F_q$-linear and induces the following short, exact sequence of $A-$modules: 
$$\xymatrix{
0 \ar[r] & \Lambda \ar[r] &\C_\infty \ar[r]^{\exp_\Lambda \quad} & \C_\infty \ar[r] & 0}
.$$
Let $\G_a$ be an additive algebraic group defined over $\C_\infty.$ Then $\End(\G_a)=\C_\infty\{\tau\}$, where $\C_\infty\{\tau\}$ is an $\F_q$-algebra of twisted polynomials, i.e. which fulfills the condition $\tau z=z^q \tau$ for $z\in\C_\infty$. 
V. Drinfeld showed that for every  $a\in A$ 
there is a twisted polynomial $\phi_a\in\C_\infty\{\tau\}$ such that the following functional equation:
$$\exp_\Lambda(az)=\phi_a(\exp_\Lambda(z))$$
holds true.
The map $a\mapsto \phi_a$ is an $\F_q-$algebra homomorphism called the  Drinfeld module associated with $\Lambda.$ We will denote it as $\phi_\Lambda:A\lra \C_\infty\{\tau\}.$ The rank of a  Drinfeld  module is the rank of the lattice  $\Lambda$ and will be denoted as  $\rk\phi_\Lambda$. Notice that the Drinfeld module $\phi_\Lambda:A\lra \C_\infty\{\tau\}$ is uniquely determined by its value at  $t$:
$$\phi_t=\theta+a_1\tau+a_2\tau^2+\cdots a_n\tau^n.$$
Also  $n=\deg_\tau\phi_t$ is equal to the rank of $\phi_\Lambda.$ More generally, one has  $\deg_\tau\phi_a=\deg_t a\cdot \rk\phi_\Lambda$ for every  $a\in A$.    

\subsection{Anderson \tm modules} 
Higher dimensional analogs of Drinfeld modules are Anderson \tm modules \cite{a}. They are defined as  $\F_q-$algebra  homomorphisms:  
\begin{equation}\label{t-mod-algebraicznie}
    \Phi:A\lra \Mat_d(\C_\infty\{\tau\})
\end{equation}
where
$$\Phi_t=\partial\Phi_t+A_1\tau+A_2\tau^2+\cdots A_n\tau^n,$$
and $\partial$ denotes a derivative on $\G_a^d$. Here $\partial \Phi_t=\theta I_d+N_\Phi$ where  $I_d$ is the identity matrix and  $N_\Phi$ is a nilpotent matrix. We say that  $\Phi$ has dimension $d$ and is of degree $\deg\Phi:=\deg_\tau\Phi_t.$ In this case there also exists a uniquely determined exponential function  $\Exp_\Phi:\C_\infty^d\lra\C_\infty^d$,
defined by the series:
$$\Exp_\Phi(\tau)=I_d\tau^0+\sum_{i=1}^\infty B_i\tau^i$$
 satisfying the following functional equation:
$$\Exp_\Phi\big(\partial \Phi_a\tau\big)=\Phi_a\big(\Exp_\Phi(\tau)\big).$$
 $\Exp_\Phi$ is defined on the whole of $\C_\infty^d,$ but does not need to be surjective.  If $\Exp_\Phi$ is surjective then we say that the \tm module $\Phi$ is uniformizable. For example, Drinfeld modules are one-dimensional, uniformizable \tm modules. Recall that $\Exp_\Phi$ is a  homomorphism of  $\F_q[t]-$modules and its kernel $\Lambda$ is a discrete finitely generated $\partial\Phi(A)-$ submodule of $\C_\infty^d$. We call $\Lambda$ a $\partial\Phi(A)-$lattice  for $\Phi.$ Its rank  is  defined to be the rank of $\Phi$.
 Notice that for a \tm module $\Phi$ of dimension $d\geq 2$ its degree $\deg\Phi$ is not necessarily equal to $\rk\Phi $ (cf. \cite[Example 5.9.9]{g96}).

Considering \tm modules in terms of the exponential functions and lattices yields an analytic approach to the theory. Here, one has to work with the field $\C_\infty$ or at least a sufficiently large subfield (containing the lattice of $\Phi$). If we consider a \tm module as a  homomorphism \eqref{t-mod-algebraicznie} 
then  we have an algebraic approach to the theory.
In this case we can consider \tm modules for 
an arbitrary $A-$field $K$, i.e. a field of characteristic $p$ equipped with a homomorphism  $\iota:A\lra K$ such that  $\iota(t):=\theta$. 
We say that an $A-$field $K$ is of $A-$ characterisic zero if $\iota$ is an inclusion.
If $\iota$ factors as $A\lra A/{\mathfrak p}\lra K$ where ${\mathfrak p}$ is a prime ideal of $A$ then we say that the $A-$characteristic of $K$ is ${\mathfrak p}.$ We will denote the $A-$ characterisic of $K$ as $\mathrm{char}_AK.$
 In this paper we are usually considering  \tm modules defined over an arbitrary $A-$field $K$. The only exception is section \ref{rozdzial:analityczny}, in which we investigate analytic properties of triangular \tm modules. So, we will work with the field  $\C_\infty$ therein.

Recall that  \tm modules form a category.  A morphism between \tm modules $\Phi$ and $\Psi$ of dimensions $d$ and $e$, respectively, is given by the matrix $F\in\Mat_{e\times d}\big( K\{\tau\}\big)$ 
satisfying the following condition:
$$F\Phi=\Psi F.$$
After one adjoins the zero \tm module $0:A\lra 0$ this category becomes additive but nonabelian. The category of \tm modules will be denoted as  \tm$\mod$, and the space of morphisms  as  $\Hom_\tau(\Phi,\Psi)$. 

\subsection{Extension of \tm modules}
Every \tm module $\Phi$ of dimension $d$ induces an $A-$module structure on the space  $K^d.$ This is given by the following formula: 
$$a*x:=\Phi_a(x),\quad \textnormal{for}\quad a\in A, x\in K^d.$$

This module is called the   Mordell-Weil  group of $\Phi$ and is denoted as $\Phi(K^d).$ Analogously, every morphism of \tm modules induces a morphism of  $A-$modules. In this way one obtains the inclusion functor from the category \tm$\mod$ into the category   $A-\mod$. This functor enables us to transfer the 
notion of exact sequences to the category of \tm modules in the following way. We say that the sequence of  \tm modules  
\begin{equation}\label{eq:ciag-dokladny}
    0\lra\Phi\lra \Upsilon\lra \Psi\lra 0
\end{equation}
is exact if the corresponding sequence of $A-$modules 
$$0\lra\Phi(K^d)\lra \Upsilon(K^{e+d})\lra \Psi(K^e)\lra 0$$
is exact. 
By  $\Ext^1_\tau(\Psi,\Phi)$ we will denote the set of exact sequences of the form  \eqref{eq:ciag-dokladny} modulo the isomorphic sequences. 
Similarly, to the classical situation of modules over algebras the set
 $\Ext^1_\tau(\Psi,\Phi)$ is an  $\F_q-$linear vector space. It is also an
  $A-$module where  multiplication by $a$ is given by the pullback by $\Psi_a$, or equivalently the   pushout by $\Phi_a$. 

For a short, exact sequence of  \tm modules there exists  a six-term  $\Hom-\Ext$ exact sequence , (see\cite[Theorem 10.2]{kk04}) for both the domain
\begin{align*}
    0&\lra\Hom_\tau(\zeta,\Phi)\lra \Hom_\tau(\zeta,\Upsilon)\lra \Hom_\tau(\zeta,\Psi)\lra \\
    &\lra\Ext^1_\tau(\zeta,\Phi)\lra \Ext^1_\tau(\zeta,\Upsilon)\lra \Ext^1_\tau(\zeta,\Psi)\lra 0,
\end{align*}
and also for the range of the $\Hom_{\tau}$.
\begin{align*}
    0&\lra\Hom_\tau(\Psi,\zeta)\lra \Hom_\tau(\Upsilon,\zeta)\lra \Hom_\tau(\Phi,\zeta)\lra \\
    &\lra\Ext^1_\tau(\Psi,\zeta)\lra \Ext^1_\tau(\Upsilon,\zeta)\lra \Ext^1_\tau(\Phi,\zeta)\lra 0.
\end{align*}

In this paper we will use a description of exact sequences by biderivations. This approach was suggested by M.A. Papanikolas and 
N. Ramachandran in \cite{pr}.  According to this, by an appropriate choice of coordinates in the  \tm module $\Upsilon$ in  $\eqref{eq:ciag-dokladny}$,  we can assume that the maps in this sequence are  an inclusion  on the second   coordinate followed by  projection onto the first coordinate, thereby we obtain the following sequence: 
$$0 \lra \Phi\uplra{\wlozenie} \Upsilon\uplra{\rzut} \Psi\lra 0.$$
Then $\Upsilon$ is given by the following block matrix:
		$$\left[\begin{array}{c|c}
			\Psi & 0\\\hline
			\delta & \Phi
		\end{array}\right]:\F_q[t]\lra {\mathrm{Mat}}_{e+d}(K\{\tau\}).$$
Therefore each extension is given by the $\F_q-$linear map
$$\delta:\F_q[t]\lra {\mathrm{Mat}}_{d\times e}(K\{\tau\})$$
such that
$$\delta_{ab}=\delta_a\Psi_b+\Phi_a\delta_b\quad \textnormal{for each}\quad a,b\in \F_q[t]$$
where  $\delta_a$  is an evaluation of  $\delta$ at $a$. 
The space of biderivations will be denoted as 
$\Der(\Psi,\Phi)$.
Two exact sequences given by the biderivations
 $\delta$ and $\widehat{\delta}$ are isomorphic if there is a matrix  $U\in {\mathrm{Mat}}_{d\times e}(K\{\tau\})$ such that 
$$\delta-\widehat{\delta}=U\Psi-\Phi U.$$
The space of  biderivations  of the form $\delta^{(U)}=U\Psi-\Phi U$ will be denoted as $\Derin(\Psi,\Phi)$ and its elements will be called  inner biderivations. Then we have the following isomorphism of $A-$modules
\begin{equation}\label{eq:iso_Ext_Der}
    \Ext^1_\tau(\Psi,\Phi)\cong \Der(\Psi,\Phi)/\Derin(\Psi,\Phi).
\end{equation}

In the sequel, we will denote an exact sequence 
 by
$$\delta: \quad 0\lra \Phi\lra \Upsilon\lra \Psi\lra 0,$$
if it is given by the biderivation $\delta$.
Notice that every biderivation
 $\delta$ is uniquely determined by its value at $t$. This yields an isomorphism of the  $\F_q-$linear spaces:
$\Der(\Psi,\Phi)\cong \Mat_{d\times e}(K\{\tau\})$ given by the assignment $\delta\mapsto \delta_t$. 
Recall that the
 $\F_q[t]$-module action on  $\Ext^1_\tau(\Psi,\Phi)$ on the level of biderivations is realized by means of  the following formula:
$$t*\Big(\delta+\Derin(\Psi,\Phi)\Big):=\Phi_t\delta+\Derin(\Psi,\Phi).$$
In the sequel, to simplify the notation we omit
  $+\Derin(\Psi,\Phi)$ when we consider the class $\delta+\Derin(\Psi,\Phi)$.

\subsection{Anderson \tm motives} 
For an $A$-field $K$ let $K[t,\tau]:=K\{\tau\}[t]$ be the polynomial ring over $K\{\tau\}$
with the following commutation relations:
$$tc=ct,\,\, t\tau=\tau t,\,\, \tau c=c^q\tau ,  \quad c\in K.$$
\begin{dfn}[\cite{a}]
A \tm motive $M$ is a left $K[t,\tau]$-module which is free and finitely generated as a $K\{\tau\}$-module for which there is an $l\in\mathbb N$ with
$$(t-\theta)^l(M/\tau M)=0.$$
A morphism of \tm motives is a morphism of $K[t,\tau]$-modules.
    \end{dfn}

With every \tm module $(\Phi, G_a^d)$ over $K$ one uniquely associates  a \tm motive in the following way:
$M(\Phi)=\Hom^q_K(G_a^d,G_a)\cong K\{\tau\}^d$ is the set of ${\mathbb F}_q$-linear morphisms of algebraic groups endowed with the following 
$K[t,\tau]$-action $(ct^i,m)\rightarrow c\circ m\circ \Phi(t^i).$ 
\begin{dfn}
 A  \tm motive $M$ is called abelian if $M$ is free and finitely generated over $K[t].$  Its 
 rank $\rk M$ as a $K[t]$-module is called the rank of $M.$ An abelian t-module is a \tm module for which the associated \tm motive is abelian.
\end{dfn}
The following theorem was proved by  G. Anderson \cite{a} (see also \cite[Theorem 2.9]{h}, \cite{bp20}).

\begin{thm}\label{equiva}
  The above
correspondence between abelian \tm motives and abelian \tm modules over $K$ gives an
anti-equivalence of categories.  
\end{thm}

In order to have the category $\cal C$ of motives which   is closed with respect to  tensor products and sub-quotients, one introduces the notion of pure \tm motives.
Set $M((\frac{1}{t})):=M{\otimes}_{K[t]}K((\frac{1}{t}))$ where
$K[[\frac{1}{t}]]$ is the ring of power series in the variable $1/t$ and $K((\frac{1}{t}))$ is its field of fractions. $M((\frac{1}{t}))$ is naturally a left $K[t,\tau]$-module. If $M((\frac{1}{t}))$ contains a finitely generated $K[[\frac{1}{t}]]$-submodule $H,$ which generates $M((\frac{1}{t}))$ over $K((\frac{1}{t}))$ and additionally satisfies
$$t^uH=t^vH$$
for some $u,v\in {\mathbb N},$ then $M$ is called pure of weight $w(M)=\frac{d(M)}{\rk (M)}=\frac{u}{v}.$ The category of pure \tm motives is closed with respect to tensor products and sub-quotients \cite{a},\cite[section 3]{h}. In particular $\Lambda^kM$ is pure for any pure \tm motive $M$ and $k\in {\mathbb N}.$ 
Moreover, \cite[Proposition 3.22 (i)]{h} asserts that if a pure \tm motive sits in an exact sequence
$$0\rightarrow M_1\rightarrow M \rightarrow M_2 \rightarrow 0$$ of \tm motives, then $M_1$ and $M_2$ are also pure and one has
\begin{equation}\label{weight}
w(M_1)=w(M_2)=w(M)
\end{equation}
We will use this fact in section \ref{sec:Weil-Bersotti}.

\subsection{ \tm comotives}

Let $L$ be an extension of ${\mathbb F}_q$ which is a perfect field. Let $L[t,\sigma ]$ be a polynomial ring with coefficients in $L$ subject to the following relations:
$$ct=tc,\quad {\sigma}t=t{\sigma}, \quad{\sigma}c=c^{1/q}{\sigma}, \qquad c\in L^*.$$
\begin{dfn}
A  \tm comotive $N$ is a left $L[t,\sigma ]$-module that is finitely generated and free over $L\{\sigma \}$,  such that there exists $l\in {\mathbb N}$ with $$(t-\theta)^l(N/{\sigma}N)=0.$$ A morphism of  \tm comotives is a morphism of $L[t,\sigma]$-modules.
\end{dfn}

With every \tm module $(\Phi, G_a^d)$ over $L$ one uniquely associates  a  \tm comotive in the following way:

Let $* : L\{\tau\}\rightarrow L\{\sigma\}$ be an anti-isomorphism defined as follows:
$(\sum a_i{\tau}^i)^* = \sum {\sigma}^i a_i = \sum a_i^{(-i)}{\sigma}^i.$ Extend this map  
to   $*: {\mathrm{Mat}}_{n\times m}L\{\tau\}\rightarrow {\mathrm{Mat}}_{m\times n}(L\{\sigma\})$ by defining $(B^{*})_{i,j}=B_{j,i}^*.$

$N(\Phi)=\Hom^q_K(G_a,G_a^d)\cong L\{\sigma\}^d$ is the set of ${\mathbb F}_q$-linear morphisms of algebraic groups endowed with the following 
$L[t,\sigma]$-action 
$$(ct^i,h)\rightarrow c\circ h\circ \Phi(t^i)^*,\quad \textnormal{where}\quad h\in N(\Phi).$$ 

The following theorem was proved G. Anderson  (see  \cite{bp20}).

\begin{thm}\label{equivall}
 Over a perfect field $L$, the above
correspondence between  \tm comotives and  \tm modules  gives an
equivalence of categories.  
\end{thm}

\begin{dfn}
A  \tm comotive $M$ is called  finite if $M$ is free and finitely generated over $K[t].$   A finite  \tm module is a \tm module for which the associated  \tm comotive is finite.
\end{dfn}

Recently A.Maurischat \cite{Mau} has clarified 
the situation and proved that, over a perfect field, a \tm module 
is abelian iff it is finite. We will use his result in section \ref{finite}.

\section{Triangular modules}

In this section we introduce triangular modules and prove their basic properties.

\subsection{$\tau$-composition series and $\tau$-co-composition series}
A concept of composition series for
 \tm modules was introduced in  \cite{gkk}. 
 According to  this, a \tm module is called  $\tau-$simple if it is not a middle term in any nontrivial, short, exact sequence. 
 Then every  \tm module $\Upsilon$ is either  $\tau-$simple or there exists a $\tau-$composition series of the form:  
\begin{equation}
    \label{ciag_kompozycyjny}
    0=\Upsilon_0\hookrightarrow \Upsilon_1 \hookrightarrow \Upsilon_2 \hookrightarrow \cdots 
\hookrightarrow \Upsilon_n=\Upsilon,
\end{equation}
where the quotients $\Upsilon_{l+1}/\Upsilon_l$ are $\tau-$simple for $l=0,1,\dots,n-1$, see \cite[Theorem 5.1]{gkk}.

Similarly,  one can 
consider  $\tau-$co-composition series of the 
form:  
\begin{equation}\label{coco}
0=\widehat\Upsilon_0\twoheadleftarrow\widehat\Upsilon_1 \twoheadleftarrow \widehat\Upsilon_2 \twoheadleftarrow \cdots 
\twoheadleftarrow \widehat\Upsilon_n=\Upsilon,
\end{equation}
where $\ker\Big(\widehat\Upsilon_l \twoheadrightarrow \widehat\Upsilon_{l-1}\Big),\quad l=n,n-1,\dots ,1$ are $\tau-$simple \tm modules.

\begin{definition}
We say that a \tm module  $\Upsilon$ is triangular if it is given by the  $\tau$-composition series  \eqref{ciag_kompozycyjny} where the sub-quotients  $\Upsilon_{l+1}/\Upsilon_l$ are  Drinfeld modules for $l=0,1,\dots, n-1$.
\end{definition}

From now on, for a triangular \tm module $\Upsilon$, we shall denote by the same symbol $\Upsilon$ the matrix describing the action of $t$, omitting the subscript $t$. The intended meaning will always be clear from the context. We adopt the same convention for the matrices associated with biderivations arising from exact sequences.

\begin{remark} Let $\Upsilon$ be a triangular \tm module given by the composition series  \eqref{ciag_kompozycyjny}. 
 Denote $\psi_{l+1}:=\Upsilon_{l+1}/\Upsilon_l$ for $l=0,1,2,\dots, n-1$. The $\tau-$composition series \eqref{ciag_kompozycyjny} implies 
 existence of the following exact sequences: 
$$\delta_l: \quad0\lra \Upsilon_l\lra\Upsilon_{l+1}\lra \psi_{l+1}\lra 0\quad\textnormal{for}\quad l=1,2,\dots, n-1,$$
where $\Upsilon_1=\psi_1$ is a  Drinfeld module.
Then
\begin{enumerate}
    \item the matrix for $\Upsilon$ is of the following form:
    \begin{equation}
        \label{modul_triangulowalny_macierz}
         \Upsilon=\begin{pNiceArray}{cccccc}[margin=15pt]
 \Block[draw=black]{1-1}<\large>{\psi_n}& 0&\cdots&   0&0 &0\\[2pt] 
         \Block[draw=black,fill=gray!20!]{5-1}<\large>{\delta_{n-1}}&\Block[draw=black]{1-1}<\large>{\psi_{n-1}}&\cdots&0 &0 &0\\[2pt]
         & \Block[draw=black,fill=gray!20!]{4-1}<\large>{\delta_{n-2}} &\ddots &\vdots &\vdots &\vdots \\[2pt] 
         & & & \Block[draw=black]{1-1}<\large>{\psi_3}&0 &0 \\[2pt] 
         & &\cdots & \Block[draw=black,fill=gray!20!]{2-1}<\large>{\delta_{2}}& \Block[draw=black]{1-1}<\large>{\psi_2}&0 \\[2pt] 
         & & & & \Block[draw=black,fill=gray!20!]{1-1}<\large>{\delta_{1}}& \Block[draw=black]{1-1}<\large>{\psi_1}
\end{pNiceArray},
    \end{equation}
    \item \tm modul $\Upsilon$ is determined by the pair  $\big(\underline{\psi},\underline{\delta}\big)$, where
    $$\underline{\psi}=\big(\psi_1,\psi_2,\dots, \psi_n\big)\quad\textnormal{and}\quad
\underline{\delta}= \big(\delta_1,\delta_2,\cdots, \delta_{n-1} \big).$$
Conversely, every pair 
$\big(\underline{\psi},\underline{\delta}\big)$ determines the \tm module $\Upsilon$. The fact that the triangular \tm module $\Upsilon$ is determined by the pair  $\big(\underline{\psi},\underline{\delta}\big)$ will be denoted in the following way: $\Upsilon=\Upsilon\big(\underline{\psi},\underline{\delta}\big)$.

\item 
Notice that $\delta_1$ is a skew polynomial and appears in the matrix $\Upsilon$ in the same row as $\psi_1$ and in the same column as $\psi_2$. 
So, we can write $\delta_1$ in the following way:
$$\delta_{1}(\tau)=\Big[\delta_{1\times 2}(\tau)\Big]=\Bigg[\sum\limits_{k=0}^{\deg\delta_1}d_{1\times 2,k}\tau^k \Bigg]$$
Using analogous convention, we can write 

$$\delta_2(\tau)=
\left[\begin{array}{l}
    \delta_{2\times 3}(\tau) \\[0.4cm]
    \delta_{1\times 3}(\tau)
\end{array}\right]=\left[\begin{array}{l}
    \sum\limits_{k=0}^{\deg\delta_2}d_{2\times 3,k}\tau^k \\[0.4cm]
    \sum\limits_{k=0}^{\deg\delta_2}d_{1\times 3,k}\tau^k
\end{array}\right].
 $$
 
 and in general 
 \begin{equation}\label{eq:wyrazy_dla_delta_i}
     \delta_l(\tau)=
     \left[\begin{array}{l}
    \delta_{l\times l+1}(\tau) \\[0.4cm]
    \delta_{l-1\times l+1}(\tau)\\[0.4cm]
    \vdots\\[0.4cm]
    \delta_{1\times l+1}(\tau)
\end{array}\right]=
\left[\begin{array}{l}
    \sum\limits_{k=0}^{\deg\delta_l}d_{l\times l+1,k}\tau^k \\[0.4cm]
    \sum\limits_{k=0}^{\deg\delta_l}d_{l-1\times l+1,k}\tau^k \\
    \qquad \quad \vdots \\
    \sum\limits_{k=0}^{\deg\delta_l}d_{1\times l+1,k}\tau^k
\end{array}\right]\quad 
 \end{equation}
 for $l=1,2,\dots, n-1.$ Notice that since degrees $\delta_{j\times l+1}$ are usually not equal, some of the coefficients $d_{j\times l+1,k}$ might be  zeros. 
\end{enumerate}
\end{remark}

\subsection{\tm reduced form}

We now turn to the problem of determining a reduced form for triangular
\tm modules.  In general, such a form will not exist over the original base
field $K$. It might  only exist after a suitable finite extension of $K$. 
\begin{thm}
\label{prop:t-reduced_form}
Let $\Upsilon=\Upsilon\big(\underline{\psi},\underline{\delta}\big)$ be a triangular \tm module of dimension $n$ defined over a field $K$.  There exists a finite extension $K(\Upsilon)/K$, such that $\Upsilon$ admits a reduced form 
$\Upsilon^{\mathrm{red}}=\Upsilon\big(\underline{\psi}^{\mathrm{red}},\underline{\delta}^{\mathrm{red}}\big)$
over $K(\Upsilon)$ satisfying
\[
\deg_\tau \delta_{i\times j}^{\mathrm{red}} \;<\; \max\{\rk \psi_i,\, \rk \psi_j\}
\qquad\text{for all }\quad i>j.
\]
\end{thm}
\begin{proof}
Let $\Upsilon$ be of the form  \eqref{modul_triangulowalny_macierz} where the entries ${\delta}_l(\tau)$ are given by \eqref{eq:wyrazy_dla_delta_i}. Denote by $r_k$ the rank of the module  $\psi_k$ and by  $a_{k,r_k}$ its leading term.

Notice that every lower triangular matrix
 $U\in\Mat_{n}\big( K\{\tau\}\big)$ with  diagonal terms equal to one induces an isomorphism  of the  \tm module $\Upsilon$ with a triangular   \tm module  $U\cdot \Upsilon\cdot U^{-1}$. Denote the entries in the matrix  $U$ as $u_{i\times j}$ in the same way as the biderivations in the   \tm module $\Upsilon$.

Consider the isomorphism $U_{i\times j}$ given by the matrix whose only non-zero entry below the diagonal is $u_{i\times j}.$
Then
\begin{align*}
   U_{i\times j}\cdot \Upsilon\cdot U^{-1}_{i\times j}=  \Upsilon+
   \left[\begin{array}{ccccccc}
    0 & \cdots &0 &  0 & 0&\cdots & 0\\
 \vdots & \ddots & \vdots &  \vdots & \vdots&\ddots & \vdots\\
     0 & \cdots &0 &  0 & 0&\cdots & 0\\
        u_{i\times j}\delta_{j\times n} & \cdots & u_{i\times j}\delta_{j\times j+1} &  \delta_{\psi_j,\psi_i}^{(u_{i\times j})} & 0&\cdots & 0\\ 
         0 & \cdots & 0 &  -u_{i\times j}\delta_{i-1\times i} & 0&\cdots & 0\\
          \vdots & \ddots & \vdots &  \vdots & \vdots&\ddots & \vdots\\
         0 & \cdots & 0 &  -u_{i\times j}\delta_{1\times i} & 0&\cdots & 0\\
   \end{array}\right],
\end{align*}
where $\delta_{\psi_j,\psi_i}^{(u_{i\times j})}$ is an inner biderivation in
$\Derin(\psi_j,\psi_i)$. This isomorphism allows one to reduce the biderivation $\delta_{i\times j}$ by means of the inner biderivation $\delta_{\psi_j,\psi_i}^{(u_{i\times j})}.$ As a result of this operation, the  entries to the left and also below the biderivation  $\delta_{i\times j}$ have been changed. 
Notice that putting $u_{i\times j}=\mu\cdot \tau^k$ for $\mu\in K$ and $k=0,1,2,\dots$ we obtain the following form for  $\delta_{\psi_j,\psi_i}^{(u_{i\times j})}$:
\begin{align*}
    \delta_{\psi_j,\psi_i}^{(\mu\tau^k)}&=f_k(\mu)\cdot \tau^{k+\max\{r_i,r_j\} } + \textnormal{terms of smaller degree},
\end{align*}
where $f_k(x)\in K[x]$.

Assume that $\deg_\tau \delta_{i\times j}=h>\max\{r_i,r_j\}$ and that its leading term is equal to  $b$.  
Setting $k=h-r_j$ and letting $\mu$ be a root of the equation $f_k(x)+b=0,$
we obtain that the degree of the biderivation at the place $i\times j$ in the \tm module
$U_{i\times j}\cdot\Upsilon\cdot U_{i\times j}^{-1}$
is less than $h$. If the equation
$f_k(x)+b=0$
has no solution in $K$, we pass to the extension $K(\mu)$ and continue the construction over this field.

The reduction algorithm proceeds by increasing distance from the diagonal.
At the first stage we consider all entries on the first sub-diagonal. Once these are reduced, we move to the second sub-diagonal, and so forth.
 At the $m$-th step we simultaneously reduce all entries at distance $m$ from the diagonal, proceeding inductively until all off-diagonal entries have been reduced.
\end{proof}

\begin{remark}
The field $K(\Upsilon)$ arising in the reduced form of a triangular \tm module
$\Upsilon$ is obtained by adjoining solutions of certain additive equations of the
form
\[
a\,x^{q^{k}} - b\,x + c = 0,\qquad a,b,c\in K,
\]
which naturally occur during the reduction process.  
Depending on whether
$a=0$,  $b=0$, or both are nonzero, these equations give rise, respectively,  to
 linear equations, purely inseparable Frobenius equations $x^{q^{k}}=d$, or 
additive $q^{k}$-polynomial equations which, after a suitable change of
variables\footnote{Such a change of variables may require passing to a finite
extension of $K$.}, reduce to the classical Artin--Schreier form
\[
y^{q^{k}} - y = d'.
\]
Thus each extension appearing in $K(\Upsilon)$ is a compositum of two fundamental
families of extensions in characteristic~$p$: purely inseparable Frobenius
extensions, and separable Artin--Schreier (or more generally
Artin--Schreier--Witt) extensions.  
The latter are well understood and have a complete cohomological
classification via Witt vectors (see Neukirch
\cite[Chap.~IV]{NeukirchANT}, Goss \cite[§12.1]{g96}, and
Hazewinkel \cite{HazewinkelWitt}; see also Serre
\cite[Chap.~III, Ex.~6]{SerreLF} for the basic case $k=1$).  
Consequently, although the reduction procedure may require adjoining to $K$ several
solutions of such equations, the resulting field $K(\Upsilon)$ always lies within
the classically studied framework of Frobenius and Artin--Schreier--Witt
extensions.
\end{remark}

For a triangular \tm -module $\Upsilon$ we call its reduced form
$\Upsilon^{\mathrm{red}}$ the \emph{reduced form} of $\Upsilon$.
Let us note that, from the proof of the preceding theorem, one can easily
extract a condition on the degrees of the biderivations which guarantees that
$K(\Upsilon)=K$, independently of the nature of the base field $K$. 

\begin{prop}\label{cor:warunek_K_Upsilon=K}
    Let $\Upsilon=\Upsilon\big(\underline{\psi},\underline{\delta}\big)$ be a triangular \tm module of dimension $n$ defined over a field $K$. 
 If the condition
\begin{equation}\label{eq:warunek_K_Upsilon=K}
        \forall_{i<j}\qquad 
        \deg_{\tau}\delta_{i\times j}\;\geq\; \max\{\rk\psi_i,\rk\psi_j\} 
        \quad\Longrightarrow\quad 
        \rk\psi_j \;>\; \rk\psi_i
\end{equation}
is satisfied, then $K(\Upsilon)=K$.
\end{prop}
\begin{proof}
 Condition~\eqref{eq:warunek_K_Upsilon=K} implies that, in the process of
constructing the reduced form, only inner biderivations are needed in those
steps where $\rk\psi_j > \rk\psi_i$.  In particular, the corresponding equation
$f_k(x)+b=0$ appearing in the proof of Proposition~\ref{prop:t-reduced_form} is linear over $K$.  Consequently, the
reduced form is defined over the original field $K$.
\end{proof}

The next definition will be used in Section~\ref{sec:Weil-Bersotti}.
 \begin{definition}\label{twocl}
       Let $\Upsilon=\Upsilon(\underline{\psi},\underline{\delta})$ be a triangular
\tm module of dimension $n$, and let
$\Upsilon^{\mathrm{red}}=\Upsilon(\underline{\psi}^{\mathrm{red}},
\underline{\delta}^{\mathrm{red}})$ denote its reduced form.
We say that the triangular \tm module $\Upsilon$ allows for
low-degree biderivations (abbreviated as LD-biderivations) if
\begin{equation}\label{eq:warunek_LD}
\deg_{\tau} \delta_{l}^{\mathrm{red}}
\;<\;
\rk\psi_{l+1}
\qquad\text{for } l=1,2,\dots, n-1.    
\end{equation}
 \end{definition}

From the 
Proposition \ref{cor:warunek_K_Upsilon=K}
we obtain the following:
\begin{cor}\label{cor:male_stopnie_delta}
   Let $\Upsilon=\Upsilon(\underline{\psi},\underline{\delta})$ be a triangular
\tm module of dimension $n$ defined over a field $K$, and suppose that
\[
\rk\psi_1 \;<\; \rk\psi_2 \;<\; \cdots \;<\; \rk\psi_n.
\]
Then $\Upsilon$ allows for LD-biderivations and $K(\Upsilon)=K$.
\end{cor}
\begin{proof}
Since $\Upsilon$ satisfies condition~\eqref{eq:warunek_K_Upsilon=K}, we have
$K(\Upsilon)=K$.  Moreover,
\[
\deg_{\tau}\delta_{i\times j}^{\mathrm{red}}
\;<\;\max\{\rk\psi_i,\rk\psi_j\}
=\rk\psi_j
\qquad\text{for } i<j,
\]
which proves the claim.
\end{proof}

It turns out that Corollary~\ref{cor:male_stopnie_delta} admits a natural
generalization to the case of non-strict inequalities, at the expense of a
possible extension of the base field~$K$.  In this situation the reduction no
longer requires any purely inseparable Frobenius extensions. Adjoining the
relevant Artin--Schreier--Witt roots is sufficient.  Hence,  the reduced form
arises after a purely separable base change of the type considered above, as
stated in the following corollary.
\begin{cor} \label{cor:male_stopnie_delta_z_rownosciami}
     Let $\Upsilon=\Upsilon(\underline{\psi},\underline{\delta})$ be a triangular
\tm module of dimension $n$ defined over a field~$K$, and suppose that
\[
\rk\psi_1 \,\leq\, \rk\psi_2 \,\leq\, \cdots \,\leq\, \rk\psi_n.
\]
Then $\Upsilon$ allows for LD-biderivations, and its reduced form
$\Upsilon^{\mathrm{red}}$ is defined over an Artin--Schreier--Witt extension
of~$K$.  \qed
\end{cor}

\section{Morphism between triangular \tm modules}

In this section we study  morphism between triangular \tm modules. We begin with the following observation.
Recall that there are  no nonzero morphisms between Drinfeld modules  of different ranks.
In general, determining when
$\Hom_{\tau}$ vanishes in the category of \tm modules is a delicate problem.
However, for triangular \tm modules we obtain the following generalization of
the aforementioned fact:
\begin{prop}\label{prop:hom_vanish}
    Let $\Upsilon = \Upsilon\big(\underline{\psi},\underline{\delta}\big)$ and $\widehat{\Upsilon} = \widehat{\Upsilon}\big(\underline{\widehat{\psi}},\underline{\widehat{\delta}}\big)$ be two triangular \tm modules. If $\rk \psi_l \neq \rk\widehat{\psi}_j$ for every   $l,j$ then 
    $$\Hom_\tau(\Upsilon,\widehat{\Upsilon})=\Hom_\tau(\widehat{\Upsilon}, \Upsilon)=0.$$
\end{prop}
\begin{proof}
    Applying the functor  $\Hom_\tau(\psi_l,-)$ consecutively to the sequences $\widehat{\delta}_j$ for $j=1,2,\dots, \dim\widehat{\Upsilon}-1$, we obtain that $\Hom_\tau(\psi_l,\widehat{\Upsilon})=0$ for every $l\in\{1,2,\dots, \dim\Upsilon\}$. Then we apply the functor 
    $\Hom_\tau(-,\widehat{\Upsilon})$ consecutively to the sequences $\delta_l$ for $l=1,2,\dots, \dim\Upsilon-1$ and obtain that $\Hom_\tau(\Upsilon,\widehat{\Upsilon})=0$. The second equality follows by symmetry. 
\end{proof}

\subsection{Triangular morphism}
 Now we focus on morphisms given by triangular matrices. It turns out that they share a lot of analogous properties with the morphisms of Drinfeld modules.
\begin{dfn}
     Let $U:\Upsilon\lra  \widehat{\Upsilon}$ be a morphism between two triangular \tm modules of the same dimension $n$. We say that $U$ is a triangular morphism over $K$ if $U\in\mathrm{Trian}_n\big(K\{\tau\}\big)$.
\end{dfn}
For a triangular morphism $U:\Upsilon \lra \widehat{\Upsilon}$, where $\dim \Upsilon =n$, we adopt the following enumeration of the entries of $U$.
\begin{equation}
    \label{eq:skladowe_morfizmu}
    U=\left[\begin{array}{ccccc}
         u_{n}&  0 & \cdots & 0& 0 \\
         u_{n-1\times n }& u_{n-1} & \cdots & 0&0 \\ 
         \vdots& \vdots & \ddots & \vdots&\vdots \\ 
         u_{2\times n }& u_{2\times n-1} & \cdots & u_{2}&0 \\ 
         u_{1\times n }& u_{1\times n-1} & \cdots & u_{1\times 2}&u_{1} 
    \end{array}\right].
\end{equation}
Notice that this enumeration agrees with the convention \eqref{eq:wyrazy_dla_delta_i}
for triangular \tm modules

 The notion of a triangular morphism is important since we have the following theorem:

\begin{thm}\label{prop:trojkatne_morfizmuy_to_wszystko}
     Let $V:\Upsilon = \Upsilon\big(\underline{\psi},\underline{\delta}\big)\lra \widehat{\Upsilon} = \widehat{\Upsilon}\big(\underline{\widehat{\psi}},\underline{\widehat{\delta}}\big)$ be a morphism defined over $K$. 
     If $\psi_i\not\sim \widehat{\psi}_j$, i.e. $\psi_i$ and $\widehat{\psi}_j$ are non-isogenous,  for $i\neq j$, then 
     $V$ is triangular morphism. 
\end{thm}

\begin{proof}
Let $V=(v_{i\times j})\in M_{n,n}(K\{\tau\})$ be a matrix, with the  usual indexing\footnote{That is, we adopt the conventional matrix indexing: the entry with index $1\times1$ is placed in the upper-left corner, the first index increases as we move down the rows, and the second index increases as we move to the right across the columns. This is the only place in the paper where we use this indexing convention.}, such that $$V\cdot \Upsilon=\widehat\Upsilon\cdot V.$$ Comparing the $1\times n$ terms of matrices in the above equality, we obtain $v_{1\times n}\cdot\phi_1 =\widehat\phi_n\cdot v_{1\times n}.$ This means that $v_{1\times n}$ is an isogeny between $\phi_1$ and $\widehat\phi_n$ and thus $v_{1\times n}=0.$
Then we consider the $2\times n$ terms.
We obtain the equality $v_{2\times n}\phi_1= \widehat\delta_{n-1\times n}\cdot v_{1\times n}+\widehat\phi_{n-1}\cdot v_{2\times n}.$ This in turn yields $v_{2\times n}=0.$  Continuing,  we see that $v_{i\times n}=0$ for $i=1,\dots,n-1$ and $v_{n\times n}$ is an isogeny between $\phi_1$ and $\widehat\phi_1.$ 
Then we compare the $1\times (n-1)$ terms and get $v_{1\times (n-1)}=0$ etc. So, going down (from the first  to the diagonal term ) in   consecutive columns starting from the last   we inductively obtain that $v_{i\times j}=0$ for $j>i$ and $v_{i\times i}$   is an isogeny between $\phi_{n-i+1}$ and $\widehat\phi_{n-i+1}$ for $i=n,\dots,1.$ The fact that $v_{i\times i}$ is an isogeny is true for general triangular morphisms (cf. Prop. \ref{basictr} (i) ).
\end{proof}

\begin{cor}\label{cortr}
Let $\Upsilon = \Upsilon\big(\underline{\psi},\underline{\delta})$ be a triangular \tm module  such that $\psi_i\not\sim {\psi}_j$ for $i\neq j.$ Then the endomorphism algebra $\End_{\tau}(\Upsilon)$ consists of triangular morphisms.\qed
\end{cor}
We have the following proposition:
\begin{prop}\label{basictr}
    Let $U:\Upsilon = \Upsilon\big(\underline{\psi},\underline{\delta}\big)\lra \widehat{\Upsilon} = \widehat{\Upsilon}\big(\underline{\widehat{\psi}},\underline{\widehat{\delta}}\big)$ be a triangular morphism defined over $K$ of the form \eqref{eq:skladowe_morfizmu}.  Then
    \begin{itemize}
        \item[$(i)$]  The maps $u_i:\psi_i\lra \widehat{\psi}_i$ for $i=1,2,\dots, n$ are morphisms of Drinfeld modules.
        \item[$(ii)$] For the morphism $U$, the following conditions hold true  
         $$\widehat{\delta}_{j\times i}u_{i}-u_{j}\delta_{j\times i}+\sum_{k=j+1}^{i-1}\Big(\widehat{\delta}_{j\times k}u_{k\times i}
         - u_{j\times k}\delta_{k\times i }\Big)=\delta_{\psi_i,\widehat\psi_j}^{(u_{j\times i})}\in\Derin\big(\psi_i,\widehat{\psi}_j\big)$$
         for every $j<i$. 
\end{itemize}   
\end{prop}

\begin{proof}
    Comparing the diagonal entries in equation $U\Upsilon =\widehat{\Upsilon}U$, we obtain
 $$u_i\psi_i=\widehat{\psi}_i u_i\quad\textnormal{for each}\quad i=n, \dots,  2,1.$$
 Thus $u_i$ are morphisms of Drinfeld modules, which proves $(i)$.

 For the proof of  $(ii)$ we again use the equality $U\Upsilon =\widehat{\Upsilon}U$. This time we compare the off diagonal terms  at places  $j\times i$ where $j<i.$ We count the rows starting from the bottom, and the columns going from  right to  left. 
Then for $j<i$  we obtain the following equalities: 
 \begin{align*}
     u_{j\times i}\psi _i + \sum_{k=j+1}^{i-1}u_{j\times k}\delta_{k\times i}+ u_j\delta_{j\times i} =
     \widehat{\delta}_{j\times i} u_i+\sum_{k=j+1}^{i-1}\widehat{\delta}_{j\times k}u_{k\times i}+\widehat{\psi}_j u_{j\times i}.
 \end{align*}
 from which we get  
 \begin{align}\label{eq:dowod_morfizmy_trojkatne}
     \widehat{\delta}_{j\times i}u_{i}-u_{j}\delta_{j\times i}+\sum_{k=j+1}^{i-1}\Big(\widehat{\delta}_{j\times k}u_{k\times i}
         - u_{j\times k}\delta_{k\times i }\Big)=\\
         =u_{j\times i }\psi_i-\widehat{\psi}_ju_{j\times i }\in\Derin\big(\psi_i,\widehat{\psi}_j\big).\nonumber
 \end{align}
\end{proof}

\begin{remark}
Let us note that, in general, for triangular morphisms, the family of diagonal morphisms does not uniquely determine the entire morphism. To see this, it suffices to consider morphisms between triangular $t$-modules of length two. Indeed, let
$$
U=\left(\begin{array}{cc}
u_2 & 0\\
u_{1\times2} & u_1
\end{array}\right),\,
V=\left(\begin{array}{cc}
u_2 & 0\\
v_{1\times2} & u_1
\end{array}\right):
\left(\begin{array}{cc}
\psi_2 & 0\\
\delta_{1\times2} & \psi_1
\end{array}\right)
\longrightarrow
\left(\begin{array}{cc}
\widehat\psi_2 & 0\\
\widehat\delta_{1\times2} & \widehat\psi_1
\end{array}\right)
$$
be triangular morphisms. Then
$$
U-V=
\left(\begin{array}{cc}
0 & 0\\
u_{1\times2}-v_{1\times2} & 0
\end{array}\right)
$$
is also a morphism of triangular $t$-modules. Hence, by \eqref{eq:dowod_morfizmy_trojkatne}, we obtain
$$
0=\bigl(u_{1\times2}-v_{1\times2}\bigr)\psi_2
-\widehat{\psi}_1\bigl(u_{1\times2}-v_{1\times2}\bigr),
$$
which implies that
$$
u_{1\times2}-v_{1\times2}\in
\Hom_\tau(\psi_2,\widehat{\psi}_1).
$$
\end{remark}

In the following proposition, we identify a situation in which the family of diagonal morphisms uniquely determines a triangular morphism.

\begin{prop}
  Let $U:\Upsilon = \Upsilon\big(\underline{\psi},\underline{\delta}\big)\lra \widehat{\Upsilon} = \widehat{\Upsilon}\big(\underline{\widehat{\psi}},\underline{\widehat{\delta}}\big)$ be a morphism defined over $K$.  If $\psi_i\not\sim \widehat{\psi}_j$, i.e. $\psi_i$ and $\widehat{\psi}_j$ are non-isogenous,  for $i\neq j$, then 
\begin{itemize}
    \item[$(i)$]  
      $U$ is determined uniquely by the family $\big\{u_i:\psi_i\lra \widehat{\psi}_i\big\}_{i=1}^n$ of morphisms  of Drinfeld modules. 
    \item[$(ii)$] If $\mathrm{char}_AK=0$ then  $\partial: \Hom_\tau(\Upsilon,\widehat{\Upsilon})\lra \mathrm{Trian}_{n}(K)$ is a monomorphism of $\F_q-$vector spaces.
    \end{itemize}
\end{prop}

\begin{proof}
It follows from Theorem~\ref{prop:trojkatne_morfizmuy_to_wszystko} that $U$ is a triangular morphism.
We use equation \eqref{eq:dowod_morfizmy_trojkatne} for the following pairs $(i-1,i)$ where $i=n,n-1,\dots, 2$, in order to obtain the terms in the first sub-diagonal of $U$. For example, from the condition 
$ \widehat{\delta}_{n-1\times n}u_{n}-u_{n-1}\delta_{n-1\times n}
         \in \Derin\big(\psi_n,\widehat{\psi}_{n-1}\big)$,
we see that there exists a unique $u_{n-1\times n }$
 such that
$$ \widehat{\delta}_{n-1\times n}u_{n}-u_{n-1}\delta_{n-1\times n}
         =u_{n-1\times n }\psi_n-\widehat{\psi}_{n-1}u_{n-1\times n }.$$
Thus, the proof of $(i)$ is completed by induction on terms on consecutive sub-diagonals.

 The fact  that $\partial$ is a linear map follows immediately from the definition. It remains to be shown that  $\partial$ is a monomorphism.
   Let $U\in \ker \partial$. Then
 \begin{align}\label{dlta}
      \partial U=\left[\begin{array}{ccccc}
         \partial u_{n}&  0 & \cdots & 0& 0 \\
          \partial u_{n-1\times n }& \partial u_{n-1} & \cdots & 0&0 \\ 
         \vdots& \vdots & \ddots & \vdots&\vdots \\ 
          \partial u_{2\times n }& \partial u_{2\times n-1} & \cdots & \partial u_{2}&0 \\ 
         \partial u_{1\times n }& \partial u_{1\times n-1} & \cdots & \partial u_{1\times 2}& \partial u_{1} 
    \end{array}\right]=0.
 \end{align}
 Since
 $\partial u_{i}=0$ and $u_i$ is a morphism of Drinfeld modules, then  $u_i=0$ for $i=n,\dots, 2,1$ see \cite[Proposition 3.3.4]{Papikian}. 
 By part~(i), we have $U=0$.
\end{proof}

Recall that the space $\Hom_\tau(\Upsilon,\widehat{\Upsilon})$ is an  $\F_q[t]-$module with the multiplication by 
$a\in \F_q[t]$ defined as
    $$a\star U:= U\cdot\Upsilon_a.$$
    Analogously, we can define an $\F_q[t]-$module structure on $\mathrm{Trian}_n(K).$ Then $\partial$ is a morphism of $\F_q[t]-$modules. 

\subsection{Triangular isogenies}
Recall also that a surjective morphism of 
 \tm modules with a finite kernel is called an isogeny, see \cite[Definition 5.1.]{hartl}. 
 The following theorem gives a characterization of triangular isogenies.
 
\begin{thm}\label{prop:triangular_isogenie}
      Let $U:\Upsilon = \Upsilon\big(\underline{\psi},\underline{\delta}\big)\lra \widehat{\Upsilon} = \widehat{\Upsilon}\big(\underline{\widehat{\psi}},\underline{\widehat{\delta}}\big)$ be a triangular morphism defined over $\overline{K}$, of the form \eqref{eq:skladowe_morfizmu} .  Then
         $U$ is an isogeny  iff $u_i:\psi_i\lra \widehat{\psi}_i$ are isogenies of Drinfeld modules for all $i=1,2,\dots, n$. 
\end{thm}

\begin{proof}
    Let  $U$ be an isogeny. If for some  $i$ we have $u_i=0$, then one sees that the kernel of  $U$ cannot be finite.
    Indeed, if $i=1$, then
$
(0,0,\ldots,0,x_1)^T\in\ker U
$
for every $x_1\in\overline{K}$. Hence, $\ker U$ is infinite.
Now suppose that $i\neq1$, and let $j$ be the smallest index such that $u_j=0$. Then
$
u_{j-1}\neq0,;\ldots,;u_1\neq0.
$
Consider a vector of the form
$
\bar{x}=(0,\ldots,0,x_j,x_{j-1},\ldots,x_1)^T.
$
We have $\bar{x}\in\ker U$ if and only if
\begin{align*}
u_{j-1\times j}x_j+u_{j-1}x_{j-1}&=0,\\
u_{j-2\times j}x_j+u_{j-2\times j-1}x_{j-1}+u_{j-2}x_{j-2}&=0,\\
\vdots\\
\sum_{l=j}^{n}u_{1\times l}x_l+u_1x_1&=0.
\end{align*}
Observe that $u_{j-1},\ldots,u_1$ are isogenies and therefore surjective. Consequently, for every choice of $x_j\in\overline{K}$, the above system can be solved recursively, starting from the first equation, to determine $x_{j-1},\ldots,x_1$. Thus, for every $x_j$ there exists a vector $\bar{x}\in\ker U$, which proves that $\ker U$ is infinite.

    Conversely, let  $U$ be a triangular morphism of the form \eqref{eq:skladowe_morfizmu} such that $u_i:\psi_i\lra \widehat{\psi}_i$ is an isogeny for every $i$. 
    
    First, we will show that $U$ has  finite kernel. 
    Assume $x=(x_n,x_{n-1},\dots,$ $x_2,x_1)^T\in\ker U$. The equation $Ux=0$ yields the following system of equations:
    \begin{align*}
        &u_{n}x_n=0\\
        &u_{n-1\times n}x_n+u_{n-1}x_{n-1}=0\\
        &\vdots\\
        &\sum_{i=2}^{n}u_{1\times i}x_i+u_1x_1=0.
    \end{align*}
    From the first equation, we see that $x_n$ is a root of the nonzero polynomial $u_{n}$. So, there are only a finite number of possible values for  $x_n$. For fixed $x_n$ the second equation  has finitely many solutions $x_{n-1}.$ 
   By induction it follows that $\ker U$ is finite.

  Now, we will show that $U$ is surjective. Let $y=(y_n,\dots, y_2,y_1)^T\in K^n$. The system $Ux=y$ gives the following system of $n$ equations: 
   \begin{align*}
        &u_{n}x_n=y_n\\
        &u_{n-1\times n}x_n+u_{n-1}x_{n-1}=y_{n-1}\\
        &\vdots\\
        &\sum_{i=2}^{n}u_{1\times i}x_i+u_1x_1=y_1,
    \end{align*}
    which again can be solved starting from the first equation then passing to the second etc. This shows that  $U$ is surjective with a finite kernel, i.e. an  isogeny.  
\end{proof}

A triangular morphism, which is an isogeny will be called a triangular isogeny.

For a Drinfeld module $\psi$ let $H(\psi)$ be its height (see \cite[Lemma 3.2.11]{Papikian}.  Then from  \cite[Proposition 3.3.4]{Papikian} we obtain the following corollary.
\begin{cor}
If 
    $U:\Upsilon = \Upsilon\big(\underline{\psi},\underline{\delta}\big)\lra \widehat{\Upsilon} = \widehat{\Upsilon}\big(\underline{\widehat{\psi}},\underline{\widehat{\delta}}\big)$ is a triangular isogeny then  $H(\psi_i)=H(\widehat{\psi}_i)$ for all $i$.\qed 
\end{cor}

According to  \cite[Corollary 5.15.]{hartl} for every isogeny of  \tm modules there exists a dual isogeny satisfying:
    \begin{align}\label{eq:iso_dual}
        U\circ U^\vee = \widehat{\Upsilon}_a\quad \textnormal{and}\quad  U^\vee\circ U = \Upsilon_a,\quad \textnormal{for some}\quad a\in\F_q[t].
    \end{align}
   So, we obtain the following corollary: 
\begin{cor}\label{cor:dual_do_morfizmu_trojkatnego}
     For every triangular isogeny there exists a dual isogeny which is triangular.  \qed
\end{cor}

\subsection{Endomorphism algebra of triangular morphisms}   

Recall that in characteristic zero the endomorphism algebra of a Drinfeld module is commutative.
 An analogous result holds true for triangular \tm modules with non-isogenous sub-quotients. 
\begin{prop} Let $\mathrm{char} _AK=0$ and  $\Upsilon = \Upsilon\big(\underline{\psi},\underline{\delta}\big)$ be a triangular \tm module of dimension $n$ with pairwise non-isogenous sub-quotients i.e. 
    $\psi_i\not\sim \psi_j$ for $i\neq j$.
    Then the endomorphism algebra 
    $\End_\tau(\Upsilon)$
    is commutative.    
\end{prop}

\begin{proof} By Corollary \ref{cortr}, $\End_{\tau}(\Upsilon)$  consists of triangular endomorphisms.
So, consider two triangular morphisms $W,V\in \End_\tau(\Upsilon)$. Then the commutator  $WV-VW$ is again a triangular endomorphism of the triangular \tm module  $\Upsilon$. Denote it by $U$, and its entries according to \eqref{eq:skladowe_morfizmu}. 
Since $u_i$ for $i=1,\dots, n$ is a commutator of endomorphisms  of the Drinfeld module  $\psi_i$, we get $u_i=0$ for $i=1,\dots , n$. Similarly as in the proof of Proposition  \ref{basictr}, putting in the equality \eqref{eq:dowod_morfizmy_trojkatne}  $i=l$ and $j=l-1$, we obtain that
$u_{l-1\times l}\in\Hom_\tau(\psi_l,\psi_{l-1})=0$. Therefore $u_{l-1\times l}=0$ for all $l$. Substituting  $i=l$ and $j=l-2$ into  equality \eqref{eq:dowod_morfizmy_trojkatne}, we obtain that
$u_{l\times l-2}\in\Hom_\tau(\psi_l,\psi_{l-2})=0$.  By a finite induction we obtain that $U=0$.  

\end{proof}

\begin{prop}\label{prop:wyznaczanie_morfizmu_z_wyrazu_wolnego}
      Let $\mathrm{char}_AK=0$ and $\Upsilon = \Upsilon\big(\underline{\psi},\underline{\delta}\big)$ be defined over  $K$ such that $\rk\psi_i\neq \rk\psi_j$ for $i\neq j$. Then 
      \begin{itemize}
          \item[(i)] For each $U_0\in\mathrm{Trian}_{n}\big(K^\mathrm{sep}\big)$, one can effectively decide whether there exists $U\in \End_{K^\mathrm{sep}}(\Upsilon)$ with $\partial U=U_0$, and, if so, find it.
          \item[(ii)] If $U\in \End_{K^\mathrm{sep}}(\Upsilon)$ has 
          $\partial U\in {\mathrm{Trian}}_{n}\big(K\big)$, then $U\in \End_{K}(\Upsilon)$.
      \end{itemize}
\end{prop}
For Drinfeld modules, a 
corresponding result was proved by  Kuhn and Pink in \cite[Proposition 5.3]{kuhnpink}. Before we give the proof of Proposition \ref{prop:wyznaczanie_morfizmu_z_wyrazu_wolnego}, we need the following lemma:
\begin{lem}\label{lemmat}
    Let
    \begin{equation}\label{eq:rownanie_z_lematu}
        \psi_t \cdot c(\tau) - c(\tau)\cdot \phi_t = w(\tau)
    \end{equation} be an equation,
    where $\psi$, $\phi$ are  fixed Drinfeld modules defined over $L\supset K$ of characteristic zero  and $w(\tau)$ is a fixed skew polynomial defined over $L$. Then
    \begin{itemize}
        \item[(i)] if $\partial w(\tau)\neq 0$, then the equation \eqref{eq:rownanie_z_lematu} has no solution.
        \item[(ii)] If $\partial w(\tau) = 0$, then every solution 
        $c(\tau)$ of the equation  \eqref{eq:rownanie_z_lematu} fulfills the following conditions:
        \begin{enumerate}
            \item $\deg_\tau c(\tau)\geq \deg_\tau w(\tau)-\max\{\rk\psi, \rk\phi\}$
            \item  $c(\tau)$ is uniquely determined by  $\partial c(\tau)$, i.e. if $c(\tau)=\sum\limits_{i=0}^rc_i\tau^i$, then $c_i=f_i(c_0)$ for some skew polynomials  $f_i(\tau)\in L\{\tau\}$ that depend on  $\psi$, $\phi$ and $w(\tau).$
            \item $c_0$ is a root of the system of exactly   $\max\{\rk\psi, \rk\phi\}$ polynomials that depend on $\phi$, $\psi$, $w(\tau)$ and $\deg_\tau c(\tau)$.
            \item additionally, if $\rk\psi\neq \rk\phi$ then 
            $\deg_\tau c(\tau)=\deg_\tau w(\tau)-\max\{\rk\psi, \rk\phi\}$.
        \end{enumerate}
    \end{itemize}
\end{lem}

\begin{proof}
Let $\psi=\theta+\sum\limits_{i=1}^m \alpha_i\tau^i$,
    $\phi=\theta+\sum\limits_{i=1}^n \beta_i\tau^i$ and
    $w(\tau)=\sum\limits_{i=0}^w w_i\tau^i$. Denote 
    $M=\max\{ n,m\}$  and consider  a solution of the following form:
    $c(\tau)=\sum\limits_{i=0}^r c_i\tau^i$. 
    Then the equality \eqref{eq:rownanie_z_lematu} can be rewritten as follows:
    \begin{align*}
    \sum^{r}_{k=1}\Bigg(\big(\theta-\theta^{(k)}\big)c_k&+
        \sum_{i=i_k}^{k-1}\Big(\alpha_{k-i}c_i^{(k-i)}-c_i\beta_{k-i}^{(i)} \Big)\Bigg)\tau^k+ \\
        &+ 
         \sum^{r+M}_{k=r+1}\Bigg(
        \sum_{i=i_k}^{k-1}\Big(\alpha_{k-i}c_i^{(k-i)}-c_i\beta_{k-i}^{(i)} \Big)\Bigg)\tau^k=\sum\limits_{k=0}^w w_k\tau^k,
    \end{align*}
    where $i_k=\max\{0,k-M\}$.  We put  $\alpha_j=0$ for $j>m$ and $\beta_l=0$ for $l>n$.  
   
    Comparing the coefficients at the least power we obtain $(i)$. Comparison of the degrees yields a contradiction if $r<w-M$, which proves $(1)$. Comparing the coefficients at $\tau^k$ for $k=1,2,\dots, r$ we obtain that  
    $$c_k=\dfrac{1}{\theta-\theta^{(k)}}\cdot 
    \sum_{i=i_k}^{k-1}\Big(c_i\beta_{k-i}^{(i)} -\alpha_{k-i}c_i^{(k-i)}\Big)+\dfrac{w_k}{\theta-\theta^{(k)}}\quad \textnormal{for}\quad k=1,2,\dots, r,$$
    where we put  $w_*=0$ for $*>w$. 
    This allows one to determine consecutively $c_1$, $c_2$, \dots, $c_r$ as polynomial functions
    of $c_0$. This proves $(2)$.

    Comparing the coefficients at $\tau^k$ for $k=r+1,r+2,\dots r+M$ we obtain the following equalities:
        $$\sum_{i=i_k}^{k-1}\Big(\alpha_{k-i}c_i^{(k-i)}-c_i\beta_{k-i}^{(i)} \Big)=w_k\quad \textnormal{for} \quad k=r+1,r+2,\dots, r+M.$$
    Substituting $c_i=f_i(c_0)$ we have
    $$g_k(c_0):=w_k-\sum_{i=i_k}^{k-1}\Big(\alpha_{k-i}f_i^{(k-i)}(c_0)-f_i(c_0)\beta_{k-i}^{(i)} \Big)=0\quad \textnormal{for} \quad k=r+1,\dots, r+M.$$
    Therefore $c_0$ is a root of the system of polynomials $g_{r+1}, g_{r+2},\dots, g_{r+M}$. This shows $(3)$.

    For the proof of $(4)$ assume that $m>n$. (The case $n<m$ is proved analogously.) Then the equation \eqref{eq:rownanie_z_lematu} can be rewritten in the following form:
    $$\psi_t\cdot c(\tau)-w(\tau)=c(\tau)\cdot \phi_t.$$
    If $r\geq w-M$ then comparison of the degrees in the above equation implies that  $r=w-M$.
\end{proof}

\begin{proof}[Proof of the Proposition \ref{prop:wyznaczanie_morfizmu_z_wyrazu_wolnego}]
 Assume we have a matrix $U_0\in \mathrm{Trian}_{n}\big(K^\mathrm{sep}\big).$ We would like to decide whether there exists a matrix $U$ of the form (\ref{eq:skladowe_morfizmu}) such that $\partial U=U_0.$ By \cite[Proposition 5.3]{kuhnpink}
 we can effectively decide whether there exist 
 morphisms $u_i\in \End(\psi_i),\,\,i=1,\dots ,n$ such that $\partial u_i=U_0(i,i)$. The diagonal term of $U_0$, and if so then effectively find them. If such morphisms exist then we apply consecutively Lemma 
   \ref{lemmat}. First, to the pairs $(i,j)$ of the form
$(l,l-1),\, l=n,\dots ,2 .$ Then $(l,l-2),\, l=n,\dots ,3$ etc. We finish the finite induction once we reach the pair $(n,1)$. Notice that the above quoted nontrivial result of N. Kuhn and R. Pink was needed to have control over the $\tau$-degree of the right-hand side of the equation (\ref{eq:rownanie_z_lematu}) at the initial step of induction,
and point $(4)$ of Lemma \ref{lemmat} was essential for providing  upper bounds for the degrees of the right-hand side of (\ref{eq:rownanie_z_lematu}) at every step of  induction. 
\end{proof}

\begin{remark}
Notice that in the proof of Proposition \ref{prop:wyznaczanie_morfizmu_z_wyrazu_wolnego} the assumption concerning different ranks of subquotients $\psi_i$ was necessary for bounding the degree of $c(\tau)$ in every inductive step . It would be nice to have a bound for degree of $c(\tau)$ in arbitrary situations, see Lemma \ref{lemmat} point $(3)$.  \end{remark}

\section{ Exact category structure}\label{sec.exact}
In \cite{g24}, Q.~Gazda studies, among other topics, the category ${\cal M}_R$ of $A$-motives (cf.~\cite[Definition~2.3]{g24}). In particular, he proves that the category ${\cal M}_R$ is exact \cite[Proposition~2.15]{g24}. Since, in the case $A={\mathbb F}_q[t]$, this category coincides with the category of Anderson \tm motives (see the introduction to Section~3 of \cite{g24}), Theorem~2.1 immediately implies that the category of abelian \tm modules is exact. In this section, we give a direct proof that the category of all \tm modules is exact.

Recall that an exact category is an additive category $\cal{A}$  with the fixed class $E$ of kernel-cokernel pairs 
 $(i,d)$, i.e. they are pairs of composable morphisms such that $i$ is a kernel of $d$ and $d$ is a cokernel of $i$, (c.f. \cite{Q}).
 We will call a morphism $i$, for which there exists a morphism $d$, such that $(i,d)\in E$,
an admissible monomorphism. Similarly,  a morphism $d$ for which  there exists $i$ such that $(i,d)\in E$, will be called an admissible epimorphism.
We require that $E$ is closed under isomorphisms and
satisfies the following axioms:
\begin{enumerate}
\item[(i)] For all objects $A\in {\cal{O}}{\cal A}$ the identity $1_A$ is an admissible monic,
\item[(ii)] For all objects $A\in \cal{O}{\cal A}$  the identity $1_A $ is an admissible epic,
\item[(iii)] The class of admissible monics is closed under composition,
\item[(iv)] The class of admissible epics is closed under composition,
\item[(v)] The pushout of an admissible monic $i: A \rightarrowtail B$   along an arbitrary morphism $f : A\rightarrow C$
exists and yields an admissible monic $j$:
\begin{equation}\label{diagram 2.4a}
\xymatrix{
{A\,\,} \ar [d]^f \ar @{>->} [r]^i& B \ar[d]^g\\
{C\,\,} \ar @{>->} [r]^{j} &  D\, 
}.
\end{equation}
\item[(v)] The pullback of an admissible epic $h: C \twoheadrightarrow D$   along an arbitrary morphism $g: B\rightarrow D$
exists and yields an admissible epic $k$:
\begin{equation}\label{diagram 2.4b}
\xymatrix{
A \ar  [d]^f \ar @{->>}[r]^{k} & B \ar[d]^g\\
C \ar @{->>}[r]^{h} &  D\, 
}
\end{equation}
\end{enumerate}
\begin{thm}\label{excat}
The category of \tm modules with the set $E$ of admissible pairs $(i,d)$ coming from  the $\tau-$exact sequences: 
\begin{equation}\label{eqex1}
\xymatrix{0 \ar[r] & {\Psi \,\,} \ar @{>->}[r]^i & \Upsilon \ar @{->>}[r]^d & \Psi \ar[r]  & 0}
\end{equation}
is an exact category.
\end{thm}
\begin{proof}
For any object $A\in \cal{A}$, the exact sequence $\xymatrix@C-0.5pc{0 \ar[r] & {A \,\,} \ar @{>->}[r]^{1_A} & A \ar @{->>} [r] & 0 \ar[r] & 0} $ (resp. 
$\xymatrix@C-0.5pc{0 \ar[r] & {0 \,\,} \ar @{>->}[r] & A \ar @{->>} [r]^{1_{A}} & A \ar[r] & 0} $) shows that $1_A$ is an admissible monic (resp. epic).
This proves (i) and (ii).

In order to prove that a composition of two admissible monomorphisms $i_2\circ i_1$  where $i_1: \xymatrix{A \,\,\ar @{>->} [r] & B } $ and  $i_2: \xymatrix{B \,\,\ar @{>->} [r]& C } $ is admissible,
consider the following block triangular matrices corresponding to $i_1$ (resp. $i_2$):
 $B=\left[\begin{array}{c|c}
	B/A & 0\\ \hline
	 \delta_1 & A
\end{array}\right]$
(resp.  $C=\left[\begin{array}{c|c}
	C/B & 0\\ \hline
	 \delta_2 & B
\end{array}\right]$
Notice that  any extension  can be described 
by a block lower triangular matrix (in a suitable basis) \cite[proof of Theorem 5.1]{gkk}

Then consider  the extension corresponding to the following matrix:
$$\widehat{C}=\left[\begin{array}{c|c||c}
	C/B & 0  & 0 \\ \hline
	 \delta_{2a} & B/A &0\\ \hline\hline
	  \delta_{2b}  & \delta_1 & A
\end{array}\right]\quad \textnormal{where} \quad
\delta_2=\left[\begin{array}{c}
	 \delta_{2a} \\ \hline
	  \delta_{2b} 
\end{array}\right].$$ 
From the form of the above matrix  we obtain the following  extension of \tm modules:

\begin{equation}\label{eqex1}
\xymatrix{0 \ar[r] & {A \,\,} \ar @{>->}[r]^{i _2\circ i_1}& C\ar @{->>}[r]^d & D \ar[r]  & 0}
\end{equation}
where $D$ is given by  $\left[\begin{array}{c|c}
	C/B & 0\\ \hline
	 \delta_{2a}& B/A
\end{array}\right].$
For the  composition of two admissible epics we proceed  analogously.
Let  $d_1: \xymatrix{A\ar @{->>} [r]& B} $   and  $d_2: \xymatrix{B\ar @{->>} [r]& C} $ be admissible epics.
Let
 $A=\left[\begin{array}{c|c}
	B& 0\\ \hline
	 \delta_1 & D
\end{array}\right]$
(resp.  $B=\left[\begin{array}{c|c}
	C& 0\\ \hline
	 \delta_2 & E
\end{array}\right]$) correspond to $d_1$ (resp. $d_2$), then the following matrix: 
$$\widehat{A}=\left[\begin{array}{c||c|c}
	C & 0  & 0 \\ \hline\hline
	 \delta_{2} & E&0\\  \hline
	  \delta_{1b}  & \delta_{1a} & A
\end{array}\right]\quad \textnormal{where} 
\quad \delta_1=\left[\begin{array}{c|c}
	 \delta_{1b} &
	  \delta_{1a} 
\end{array}\right]$$ 

\noindent
shows that $d_2 \circ d_1$ is an admissible epimorphism.

Parts (iv) and (v)  are the content of \cite[Theorem 10.1 ]{kk04}.
\end{proof}
Let $\mathcal{T}$ denote the full subcategory of triangular \tm modules, and let
$\mathcal{T}^{\triangle}$ denote the category of triangular \tm modules with triangular morphisms.
The following corollary is an immediate consequence of the preceding theorem.

\begin{cor}
The categories $\mathcal{T}^{\triangle}$ and $\mathcal{T}$ are exact subcategories of the category of \tm modules.
\end{cor}

\section{Exponent for an exact sequence}	\label{rozdzial:analityczny}

In this section we consider analytical properties of  triangular 
\tm modules. \ We start with the following lemma.

	\begin{lem}
		\label{lem:o_jednoznacznym_rozwiazaniu_rownania_macierzowego}
		Let $Y\in\mathrm{Mat}_{n\times m}(K)$ be an arbitrary matrix and   $M\in\mathrm{Mat}_{ m}(K)$, $N\in \mathrm{Mat}_{n}(K)$ be lower triangular nilpotent matrices. Then there exists a unique matrix $C\in \mathrm{Mat}_{n\times m}(K)$
  which fulfills the following equation
		\begin{equation}
			\label{eq:macierzowe_rownanie_z_lematu}
			C+CM+NC=Y.
		\end{equation} 
	\end{lem}
	\begin{proof}
    We determine the entries of the matrix $C$ row by row, starting from the top
row and proceeding downward.  Within each row, the entries are computed from
right to left.
	\end{proof}

    From this point until the end of this section we work over the field 
$K=\C_\infty$.  The following theorem describes the behavior of the 
exponential function associated with the middle term of a short exact sequence.  
Let us note that this fact, in the case where $\partial\delta=0$, was 
already established in \cite{pr}. 

	\begin{thm}\label{thm62} 
		Let $0\lra \Phi\lra \Upsilon \lra \Psi\lra 0$ be an exact sequence of \tm modules given by the biderivation  $\delta$. Then there exists the following commutative diagram of $\F_q[t]-$modules with exact rows:
		\begin{equation}
			\label{diag:twierdzenie_1}
			\xymatrix@+0.5pc{
			0 \ar[r]& \C_\infty^d\ar[r]^{\wlozenie} \ar[d]^{\Exp_\Phi} & \C_\infty^{d+e}  \ar[r]^{\rzut} \ar[d]^{\Exp_\Upsilon} & \C_\infty^e \ar[r] \ar[d]^{\Exp_\Psi} & 0\\ 
			0 \ar[r] & \C_\infty^d \ar[r]_{\wlozenie} & \C_\infty^{d+e}\ar[r]_{\rzut} & \C_\infty^e \ar[r] & 0,}
		\end{equation} 
	where $\Exp_\Phi,$ resp. $\Exp_\Psi$, are exponents for $\Phi,$ resp. $\Psi.$ 
	The exponent of $\Upsilon$ has the following form:
	\begin{equation}
		\label{eq:postaci_e_Upsilon}
		\Exp_\Upsilon(\tau)=\left[\begin{array}{c|c}
			\Exp_\Psi(\tau) & 0\\ \hline
			\Exp_\delta(\tau) & \Exp_\Phi(\tau),
		\end{array} \right],
	\end{equation}
where $\Exp_{\delta}(\tau)$ is uniquely determined by the following equation:
\begin{equation}
	\label{eq:rownanie_na_exp_delta}
	\Exp_\delta\big(\partial\Psi\tau^0\big)-\Phi \Exp_\delta(\tau)=\delta \Exp_\Psi(\tau)-\Exp_\Phi\big(\partial\delta\tau^0\big).
\end{equation}
	\end{thm}
\begin{proof}
Notice that if we show that
	 $\Exp_\Upsilon(\tau)$ is of the form \eqref{eq:postaci_e_Upsilon}, then commutativity of the diagram  \eqref{diag:twierdzenie_1} is obvious.

 	Since $\Exp_\Upsilon(\tau) = I\tau^0+\sum\limits_{i\geq1}C_i\tau^i$ fulfills the following equation: 
	$\Exp_\Upsilon(\partial \Upsilon\tau^0)=\Upsilon \Exp_\Upsilon(\tau)$ we obtain the following equality of the formal power series: 
	\begin{align*}
		 \sum\limits_{i\geq1}\Big( C_i\partial\Upsilon_t^{(i)} -\partial\Upsilon_tC_i \Big)\tau^i &= \sum\limits_{i=1}^{s} U_i\tau^i+ 
		 \sum\limits_{i\geq1}U_1C_i^{(1)}\tau^{i+1}+\\
		 &+\sum\limits_{i\geq1}U_2C_i^{(2)}\tau^{i+2}+ \cdots+ 
		 \sum\limits_{i\geq1}U_{s}C_i^{(s)}\tau^{i+s}
	\end{align*}
	where $\Upsilon = \partial\Upsilon\tau^0+\sum\limits_{i=1}^{s}U_i\tau^i$.
	Comparing coefficients at $\tau^k$, we obtain
	\begin{align}
	C_k\partial\Upsilon^{(k)} -\partial\Upsilon C_k&=U_k+\sum_{i=1}^{k-1} U_iC_{k-i}^{(i)} & \textnormal{for}\quad & k\leq s  \label{eq_dowod_1}\\
	C_k\partial\Upsilon^{(k)} -\partial\Upsilon C_k&=\sum_{i=1}^{s} U_iC_{k-i}^{(i)} & \textnormal{for}\quad & k> s  \label{eq_dowod_2}
	\end{align} 
	Let
	$\Phi = (\theta I + N_\Phi)\tau^0+\sum\limits_{i=1}^{\textnormal{finite}}A_i\tau^i$, 
	$\Psi = (\theta I + N_\Psi)\tau^0+\sum\limits_{i=1}^{\textnormal{finite}}B_i\tau^i$ and 
	$\delta = \sum\limits_{i=0}^{\textnormal{finite}}D_i\tau^i$. Without loss of generality we can assume that $N_{\Phi}$ and $N_\Psi$ are lower triangular matrices. 
	Then 
	$$\Upsilon=\left[\begin{array}{c|c}
		\Psi & 0\\ \hline
		\delta & \Phi
	\end{array}\right]= 
\left[\begin{array}{c|c}
	\theta I + N_\Psi & 0\\ \hline
	D_0 & \theta I + N_\Phi
\end{array}\right]\tau^0+ 
\sum\limits_{i=1}^{s}\left[\begin{array}{c|c}
	B_i & 0\\ \hline
	D_i &  A_i
\end{array}\right]\tau^i.$$
Putting $C_i=\left[\begin{array}{c|c}
	C_{1,i} & C_{2,i}\\ \hline
	C_{3,i} & C_{4,i}
\end{array}\right]$ in the equalities \eqref{eq_dowod_1} and \eqref{eq_dowod_2}, we obtain $C_{2,i}=0$ for $i\geq 1$. Indeed, comparing the blocks in the upper right corners we obtain the following relation:
$$C_{2,1}\Big(\theta^{(1)} I+N_\Phi^{(1)}\Big)-\Big(\theta I+N_\Psi\Big)C_{2,1}=0,$$
which can be rewritten in the following form:
$$\Big(\theta^{(1)}-\theta\Big)C_{2,1}+C_{2,1}N_{\Phi}^{(1)}-N_\Psi C_{2,1}=0.$$
Dividing both sides by $\theta^{(1)}-\theta$, we obtain  \eqref{eq:macierzowe_rownanie_z_lematu}, which by lemma \ref{lem:o_jednoznacznym_rozwiazaniu_rownania_macierzowego} has the unique solution $C_{2,1}=0$. 
Assume  that the subsequent matrices
 $C_{2,i}$ for $i=1,2,\cdots k-1$ are zero matrices. We will show the  $k$-th matrix is zero. Since the arguments for both $k\leq s$ and $k>s$ are analogous, 
 we will consider only the first case.
 In the formula \eqref{eq_dowod_1} all block matrices  $U_i$ and $C^{(i)}_{k-i}$ are lower triangular. Therefore the matrix  	$C_k\partial\Upsilon^{(k)} -\partial\Upsilon C_k$ is also block lower triangular.  
 Thus, we have the following equality 
$$C_{2,k}\Big(\theta^{(k)} I+N_\Phi^{(k)}\Big)-\Big(\theta I+N_\Psi\Big)C_{2,k}=0,$$
which by lemma \ref{lem:o_jednoznacznym_rozwiazaniu_rownania_macierzowego} has the only solution $C_{2,k}=0$. 

We obtained the following formula:
$$\Exp_\Upsilon=I\tau^0+\sum\limits_{i\geq 1}\left[\begin{array}{c|c}
	C_{1,i} & 0\\ \hline
	C_{3,i} & C_{4,i}
\end{array}\right]\tau^i=
\left[\begin{array}{c|c}
	I\tau^0+\sum\limits_{i\geq 1}C_{1,i}\tau^i & 0\\ \hline
	\sum\limits_{i\geq 1}C_{3,i}\tau^i & I\tau^0+\sum\limits_{i\geq 1}C_{4,i}\tau^i
\end{array}\right].$$
Substituting this into equality $\Exp_\Upsilon(\partial \Upsilon\tau^0)=\Upsilon\Exp_\Upsilon(\tau)$ and comparing the block matrices on the diagonal,
we obtain the following two conditions:
\begin{align*}
	&\Big(I\tau^0+\sum\limits_{i\geq 1}C_{1,i}\tau^i\Big)\partial \Psi\tau^0=\Psi\Big(I\tau^0+\sum\limits_{i\geq 1}C_{1,i}\tau^i\Big)\\
	&\Big(I\tau^0+\sum\limits_{i\geq 1}C_{4,i}\tau^i\Big)\partial \Phi\tau^0=\Phi\Big(I\tau^0+\sum\limits_{i\geq 1}C_{4,i}\tau^i\Big).
\end{align*}
The first implies that $I\tau^0+\sum\limits_{i\geq 1}C_{1,i}\tau^i=\Exp_{\Psi}(\tau)$, whereas the second implies that  $I\tau^0+\sum\limits_{i\geq 1}C_{4,i}\tau^i=\Exp_{\Phi}(\tau)$.
Finally, comparing the matrices in the  lower left corner and denoting
 $\Exp_\delta(\tau)=\sum\limits_{i\geq 1}C_{3,i}\tau^i$ we obtain the following equation
\begin{align*}
	\Exp_\delta\Big(\partial\Psi\tau^0\Big)+
	\Exp_\Phi\Big(\partial\delta\tau^0\Big)=
	\delta\Exp_{\Psi}(\tau)+\Phi\Exp_\delta(\tau),
\end{align*}
which is equivalent to \eqref{eq:rownanie_na_exp_delta}. 

To finish the proof one has to show that the equality  \eqref{eq:rownanie_na_exp_delta} determines $\Exp_\delta(\tau)$ uniquely. 
Recall that $\Exp_\delta(\tau)$ has zero constant term.

In order to do this, notice that the right-hand side of the aforementioned equality is of the following form:
$\sum\limits_{i\geq 1}Y_i\tau^i$ for certain matrices $Y_i$. Then putting $\Exp_\delta(\tau)=\sum\limits_{i\geq 1}E_{i}\tau^i$  and comparing the coefficients at  $\tau^k$ we obtain the following equalities
\begin{align*}
	E_k\Big(\theta I+N_{\Psi}\Big)^{(k)}-\Big(\theta I+N_{\Phi}\Big)E_k= Y_k+\sum\limits_{i=1}^{m}A_iE_{k-i}^{(i)},
\end{align*}
where $m=\min\{k-1,\deg\Phi\}$. Assuming inductively that the matrices  $E_1$, $E_2$, $\dots$, $E_{k-1}$ have already been  determined, one can find  $E_k$ by transforming the above equality to the form \eqref{eq:macierzowe_rownanie_z_lematu} and applying 
lemma \ref{lem:o_jednoznacznym_rozwiazaniu_rownania_macierzowego}.
\end{proof}	

As a consequence of this theorem, we obtain the following interesting fact 
about uniformizable \tm modules.

\begin{cor}
	The class of uniformizable \tm modules is closed with respect to extensions  and images.
\end{cor}
\begin{proof}
	Applying  the snake lemma  to the diagram \eqref{diag:twierdzenie_1} we obtain the following exact sequence:
	\begin{align}
		\label{eq:ciag_dokladny_ker_coker}
		0&	\lra \ker \Exp_\Phi \lra \ker \Exp_\Upsilon\lra \ker \Exp_\Psi \lra \\
		&\lra \coker \Exp_\Phi \lra \coker \Exp_\Upsilon\lra \coker \Exp_\Psi \lra 0,\nonumber
	\end{align} 
from which the corollary  follows immediately.
\end{proof}

Moreover, in the case of uniformizable \tm modules we obtain the following 
description of the lattice corresponding to the middle \tm module in the short 
exact sequence.

\begin{cor}\label{c64}
If in a short exact sequence
	 $0\lra\Phi\lra \Upsilon\lra \Psi\lra 0$,
	a \tm module $\Phi$ is uniformizable, then the lattice  $\Lambda_\Upsilon$ as an  $\F_q[t]-$module is isomorphic to the product of lattices $\Lambda_\Phi\times\Lambda_\Psi$.  
\end{cor}
\begin{proof}
	From \eqref{eq:ciag_dokladny_ker_coker} we obtain the following short exact sequence:
	$$0	\lra \ker \Exp_\Phi \lra \ker \Exp_\Upsilon\lra \ker \Exp_\Psi \lra 0.$$
	Since  $\ker \Exp_\Psi=\Lambda_\Psi$ is a projective $\F_q[t]-$module this sequence splits. Thus,  
	$\Lambda_\Upsilon=\ker\Upsilon\cong \ker \Exp_\Phi\times \ker \Exp_\Psi=\Lambda_\Phi\times\Lambda_\Psi,$ as $\F_q[t]-$modules. 
\end{proof}

Next, we may extend the last two corollaries by induction to 
$\tau$-composition series of \tm modules.

\begin{cor}
Let 
$$\Upsilon_1 \hookrightarrow \Upsilon_2 \hookrightarrow \cdots 
\hookrightarrow \Upsilon_n=\Upsilon,$$
be a $\tau-$composition series where the quotients $\Psi_i=\Upsilon_{i+1}/\Upsilon_i$ are uniformizable \tm modules for $i=0,1,\dots, n-1$. Then 
	\begin{itemize}
		\item[$(i)$] $\Upsilon$ is a uniformizable \tm module.
		\item[$(ii)$] $\Lambda_\Upsilon\cong\prod_{i=0}^{n-1}\Lambda_{\Psi_i}$ as $\F_q[t]-$modules.
	\end{itemize} 
\end{cor}
\begin{proof}
	Consider the following exact sequences:
	$$\eta_i:\quad0\lra	\Upsilon_{i}\lra \Upsilon_{i+1}\lra \Psi_i\lra 0$$
	for $i=1,2,\dots, n-1$. Applying the last two corollaries successively  to  sequences
	$\eta_1$, $\eta_2$, $\dots$, $\eta_{n-1}$ 
 we obtain the assertion. 
\end{proof}

\begin{cor}\label{cor:triangular_is_uniformizable}
    Every triangular \tm module is uniformizable. \qed
\end{cor}
\begin{cor}
     Triangular \tm module $\Upsilon = \Upsilon\big(\underline{\psi},\underline{\delta}\big)$ has rank $r=\sum\limits_i\rk\psi_i$. 
\end{cor}
\begin{proof}
Let $r_i=\rk\psi_i$ and
${\lambda}_{i,k_i},\, k_i=1,\dots, \rk\psi_i$ be the $({\partial}\psi_i)(A)$ generators of $\Lambda_{{\psi}_i}.$ Applying an easy induction using Corollary
\ref{c64} one sees that the $(\partial\Upsilon)(A)$-basis of $\Lambda_{\Upsilon}$ is given by 
\begin{align*}
    &\big[{\lambda}_{n,1},{\alpha}_{n,1,2},\dots, {\alpha}_{n,1,n}\big]^T,\dots ,\big[{\lambda}_{n,r_n},{\alpha}_{n,r_n,2},\dots, {\alpha}_{n,r_n,n}\big]^T,\\
    &\big[0,{\lambda}_{n-1,1},{\alpha}_{n-1,1,3}\dots ,{\alpha}_{n-1,1,n}\big]^T ,\dots ,\big[0,\dots, 0,{\lambda}_{1,k_1}\big]^T
\end{align*}
for some $\alpha_{s,u,v}\in {\mathbb C}_{\infty}$. The corollary follows.
\end{proof}    

\section{Finiteness and purity}\label{finite}

In this section, we prove that triangular modules are finite. We also establish a necessary and sufficient condition for purity. As a consequence, we derive several conclusions concerning Taelman’s conjecture.

\subsection{Strict purity and almost strict purity}
Recall that a $t$-module $\Upsilon$ with
\[
\Upsilon \;=\; \partial\Upsilon\,\tau^0 \;+\; \sum_{i=1}^{s} A_i \tau^i
\]
is called \emph{strictly pure} if the leading matrix $A_s$ is invertible
(cf.\ \cite[Section~1.1]{np}), and \emph{almost strictly pure} if there exists
an integer $n$ such that the leading matrix of the polynomial
$\Upsilon^n$ is invertible.  
Both strictly pure and almost strictly pure $t$-modules play a significant
role in the arithmetic and analytic theory of $t$-modules.  
Almost strictly pure $t$-modules arise naturally in the study of periods,
quasi-periods, and logarithms
\cite{LucasPurityAS,np},  
whereas strictly pure $t$-modules are precisely the class for which one can 
construct natural $t$-module structures on extension groups 
$\Ext^1_{\tau}(\Phi,\Psi)$ \cite{kk04,gkk}, and for which duality theorems and 
Weil--Barsotti type formulas have been established 
\cite{ta,pr,kk24}.  

From the matrix expression \eqref{modul_triangulowalny_macierz}, we therefore
obtain the following characterization.
\begin{prop}\label{prop:strictly_pure_characterisation}
    Let $\Upsilon = \Upsilon\big(\underline{\psi},\underline{\delta}\big)$ be a triangular module of dimension  $n$. Then 
    $\Upsilon$ is strictly pure iff $$\deg_\tau\Upsilon=\rk\psi_1=\rk\psi_2=\cdots=\rk\psi_n.$$ \qed
\end{prop}
The situation for almost strictly pure \tm modules is more subtle. The following proposition holds.

\begin{prop}\label{prop:almost_strictly_pure}
Each triangular almost strictly pure \tm module $\Upsilon=\Upsilon\big(\underline{\psi},\underline{\delta}\big)$ satisfies
$$\rk\psi_1=\rk\psi_2=\cdots=\rk\psi_n.$$
\end{prop}

\begin{proof}
Suppose that $\operatorname{rk}\psi_i\neq \operatorname{rk}\psi_j$
for some $i\neq j$. Let $r$ be the maximal rank among the Drinfeld modules $\psi_i$. Write the \tm module $\Upsilon$ in the form
$$
\Upsilon
=
\sum_{i=0}^{r}A_i\tau^i
+
\sum_{j=r+1}^{R}N_j\tau^j,
$$
where the matrices $A_i$ are lower triangular and the matrices $N_j$ are strictly lower triangular, and hence nilpotent.
Consider the product
$$
\Upsilon^k=\underbrace{\Upsilon\cdots\Upsilon}_{k\text{ factors}}.
$$
In its expansion, the matrices $A_i$ and $N_j$ may occur in arbitrary order. Thus, we obtain all possible words of length $k$ formed from the matrices $A_i,\qquad 0\leq i\leq r,$
and
$N_j,\qquad r+1\leq j\leq R,$
with the appropriate Frobenius twists. More explicitly,
$$
\begin{aligned}
\Upsilon^k
={}&
\sum_{i_1=0}^{r}\cdots\sum_{i_k=0}^{r}
A_{i_1}A_{i_2}^{(i_1)}
\cdots
A_{i_k}^{(i_1+\cdots+i_{k-1})}
\tau^{i_1+\cdots+i_k}
\\
&+
\sum_{i_1=0}^{r}\cdots
\sum_{i_{k-1}=0}^{r}
\sum_{j_k=r+1}^{R}
A_{i_1}A_{i_2}^{(i_1)}
\cdots
A_{i_{k-1}}^{(i_1+\cdots+i_{k-2})}
N_{j_k}^{(i_1+\cdots+i_{k-1})}
\\
&\hspace{60mm}\cdot
\tau^{i_1+\cdots+i_{k-1}+j_k}
\\
&+
\sum_{i_1=0}^{r}\cdots
\sum_{i_{k-2}=0}^{r}
\sum_{j_{k-1}=r+1}^{R}
\sum_{j_k=r+1}^{R}
A_{i_1}A_{i_2}^{(i_1)}
\cdots
N_{j_{k-1}}^{(i_1+\cdots+i_{k-2})}
\\
&\hspace{60mm}\cdot
N_{j_k}^{(i_1+\cdots+i_{k-2}+j_{k-1})}
\tau^{i_1+\cdots+i_{k-2}+j_{k-1}+j_k}
\\
&+\cdots
\\
&+
\sum
A_{i_1}
N_{j_2}^{(i_1)}
A_{i_3}^{(i_1+j_2)}
\cdots
N_{j_k}^{(i_1+j_2+\cdots)}
\tau^{i_1+j_2+\cdots}
\\
&+\cdots
\\
&+
\sum_{j_1=r+1}^{R}\cdots\sum_{j_k=r+1}^{R}
N_{j_1}
N_{j_2}^{(j_1)}
\cdots
N_{j_k}^{(j_1+\cdots+j_{k-1})}
\tau^{j_1+\cdots+j_k}.
\end{aligned}
$$
The omitted terms correspond to all other possible interlacings of the matrices $A_i$ and $N_j$.

Since the matrices $A_i$ are lower triangular and the matrices $N_j$ are strictly lower triangular, every summand containing at least one $N_j$ is strictly lower triangular, and hence nilpotent. It may of course happen that some, or even all, of these summands vanish.
Therefore, it suffices to consider the coefficient of $\tau^{kr}$, which is the highest power of $\tau$ appearing in the first sum. Indeed, since
$0\leq i_\ell\leq r,$ the equality
$i_1+\cdots+i_k=kr$
is possible only when
$i_1=\cdots=i_k=r.$
Consequently, the coefficient of $\tau^{kr}$ is
$$
(\star)\qquad
A_rA_r^{(r)}
A_r^{(2r)}
\cdots
A_r^{(kr-r)}.
$$
Let $i_0$ be the index of a Drinfeld module of maximal rank, and let $a$ denote its leading coefficient. Since $A_r$ is lower triangular, the $(i_0,i_0)$-entry of the matrix $(\star)$ is
$$
aa^{(r)}a^{(2r)}\cdots a^{(kr-r)}\neq0.
$$
On the other hand, since the ranks of the Drinfeld modules $\psi_i$ are not all equal, the matrix $A_r$ is non-invertible. Because $A_r$ is lower triangular, at least one of its diagonal entries is zero. Let the corresponding index be $i_1$. Since Frobenius preserves zero, the $(i_1,i_1)$-entry of every Frobenius twist
$A_r^{(mr)},$ $ m=0,\ldots,k-1,$

is also zero. Hence the $(i_1,i_1)$-entry of the product $(\star)$ is zero. Thus $(\star)$ is non-invertible.
We have therefore shown that the coefficient of $\tau^{kr}$ in $\Upsilon^k$ is non-invertible, whereas all terms of $\Upsilon^k$ containing at least one $N_j$ are strictly lower triangular. Consequently, $\Upsilon^k$ cannot be invertible for any $k\geq1$. Hence $\Upsilon$ cannot be almost strictly pure.
\end{proof}

\begin{remark}\label{rem:almost_stricly_pure}
The argument in the proof of the preceding proposition shows that every triangular almost strictly pure \tm module is of the form
$$
\Upsilon=\sum\limits_{i=0}^{r}A_i\tau^i+
\sum\limits_{j=r+1}^R N_j\tau^j,
$$
where $\sum\limits_{i=0}^{r}A_i\tau^i$ is a triangular strictly pure \tm module, while the matrices $N_j$ are lower triangular nilpotent matrices satisfying a certain system of matrix equations.

Fix a strictly pure \tm module
$$
\sum\limits_{i=0}^{r}A_i\tau^i.
$$
Starting with $k=2,3,\dots$, one can determine the polynomial relations that the matrices $N_j$ must satisfy in order for the corresponding \tm module to be almost strictly pure. Since the coefficient of $\tau^{kr}$ in the product $\Upsilon^k$ is invertible, all coefficients of $\tau^j$ with $j>kr$ must vanish.

At first glance, one might expect that, for fixed $R$, the number of equations grows linearly with $k$. This is, however, not the case. Indeed, the matrices $N_j$ are lower triangular with zero diagonal entries, and therefore any product of $n=\dim\Upsilon$ such matrices is equal to the zero matrix. Moreover, the matrices $A_i$ are lower triangular, so both products $A_iN_j^{(i)}$ and $N_jA_i^{(j)}$ also have zero diagonal entries. Consequently, if the product
$$
A_{i_1}N_{j_2}^{(j_2)}\cdots A_{i_{k-1}}^{(i_1+\cdots+i_{k-2})}N_{j_k}^{(i_1+i_2+\cdots+i_{k-1})}
$$
contains at least $n$ matrices of the form $N_j$ (not necessarily consecutive), then this product is the zero matrix.

As a consequence, if $k\geq n$, then
$$
\deg\Upsilon^k=r\cdot k+(n-1)\cdot (R-r).
$$
Hence, for every $k\geq n$, we always obtain the same number, namely $(n-1)(R-r)$, of matrix equations that the matrices $N_j$ must satisfy.
 Note that the matrix equations satisfied by the matrices $N_j$ can be viewed as systems of polynomial equations in the entries of the matrices $N_j$, where the defining polynomials are additive. We illustrate this phenomenon with the following example.

\end{remark}

\begin{ex}
Consider the triangular \tm module
$${\Upsilon^\vee}^\vee=
\left(\begin{array}{ccc}
        \theta+\tau^3 & 0 & 0\\
        a\tau^2 & \theta+\tau^3 & 0\\
        b\tau^3 & c\tau & \theta+\tau^3
    \end{array}\right)+
    \left(\begin{array}{ccc}
        0 & 0 & 0\\
        \alpha\tau^4 & 0 & 0\\
        \beta\tau^4 & \gamma\tau^4 & 0
    \end{array}\right).
$$
A direct computation shows that the leading coefficient matrix of $\Upsilon^2$ is invertible if and only if the coefficient matrices of $\tau^7$ and $\tau^8$ vanish, namely
\begin{align*}
\left(\begin{array}{ccc}
         0 & 0 & 0\\
        \alpha+\alpha^{(3)} & 0 & 0\\
        \beta+\beta^{(3)} & \gamma+\gamma^{(3)} & 0
    \end{array}\right)&=0,\\
\left(\begin{array}{ccc}
          0 & 0 & 0\\
        0 & 0 & 0\\
        \alpha^{(4)}\gamma & 0 & 0
    \end{array}\right)&=0.
\end{align*}

Similarly, the leading coefficient matrix of $\Upsilon^3$ is invertible if and only if
\begin{align*}
\left(\begin{array}{ccc}
         0 & 0 & 0\\
        \alpha+\alpha^{(3)}+\alpha^{(6)} & 0 & 0\\
        \beta+\beta^{(3)}+\beta^{(6)} & \gamma+\gamma^{(3)}+\gamma^{(6)} & 0
    \end{array}\right)&=0,\\
\left(\begin{array}{ccc}
          0 & 0 & 0\\
        0 & 0 & 0\\
        \alpha^{(4)}\gamma+\alpha^{(7)}\Big(\gamma+\gamma^{(3)}\Big) & 0 & 0
    \end{array}\right)&=0.
\end{align*}
    \end{ex}

Moreover, Propositions \ref{prop:strictly_pure_characterisation}, \ref{prop:almost_strictly_pure} and \ref{cor:male_stopnie_delta_z_rownosciami} imply the following corollary.
\begin{cor}
Let $\Upsilon=\Upsilon(\underline{\psi},\underline{\delta})$ be a triangular \tm module, defined over a field~$K$,  such that
$\rk\psi_1=\rk\psi_2=\cdots=\rk\psi_n$. Then there exists
an Artin--Schreier--Witt extension $K(\Upsilon)/K$ such that $\Upsilon$ becomes
isomorphic over $K(\Upsilon)$ to a strictly pure \tm module. In particular, this applies to all almost strictly pure \tm modules.
\end{cor}
\begin{proof}
  Since $\Upsilon$ allows for LD-biderivations, its reduced form
$\Upsilon^{\mathrm{red}}$ is a strictly pure \tm module.
\end{proof}

\subsection{Proofs of finiteness and purity}
We now prove that triangular \(t\)-modules are finite. To this end, we use the results of A. Maurischat \cite{Mau}. The finiteness criterion developed there is based on the Smith normal form of a matrix over a suitable skew ring.
We recall that, by \cite[p. 41]{jac43}, any matrix $M\in\operatorname{Mat}_{n\times m}(\mathcal O),$
where \(\mathcal O\) is both a left and a right Euclidean domain, can be reduced to diagonal form by elementary operations. More precisely, there exist
$P\in\operatorname{GL}_n(\mathcal O),$
$Q\in\operatorname{GL}_m(\mathcal O),$
such that $M_{\mathrm{red}}=PMQ$
 has the form 
 such that
$(M_{\mathrm{red}})_{i,i}=d_i,
$ $
(M_{\mathrm{red}})_{i,j}=0
\quad\text{otherwise},
$
for
$i=1,\ldots,\min(m,n).$
The invariant factors $d_i$ satisfy the left and right divisibility relations
$d_1\mid d_2\mid\cdots\mid d_{\min(m,n)}.$

Now, consider the twisted polynomial ring $K\{\tau\},$ i.e., a polynomial ring in the variable $\tau$ with the commutativity relation
$\tau\cdot\alpha=\alpha^q\cdot\tau.$
Assume that $K$ is perfect. Let, as in \cite{Mau},
$K\{\{\sigma\}\}=
\sum_{i=0}^{\infty}\alpha_i\sigma^i,\quad \alpha_i\in K
,$
where $\sigma=\tau^{-1},$ be the ring of power series in $\sigma$ with coefficients in $K.$ Then
$\sigma\cdot\alpha=\alpha^{1/q}\cdot\sigma, \alpha\in K,$
and the ring of Laurent series
$$K(\{\sigma\})=
\sum_{i=i_0}^{\infty}\alpha_i\sigma^i,\quad \alpha_i\in K
$$
is a skew field.
The ring $K\{\tau\}$ can be naturally embedded in $K(\{\sigma\})$ by the formula
\begin{equation}\label{embd}
\sum_{i=0}^{n}\alpha_i\tau^i
\longrightarrow
\sum_{i=0}^{n}\alpha_i\sigma^{-i}.
\end{equation}

Notice that in the polynomial ring $K(\{\sigma\})[t]$ there exist both a left and a right division algorithm. Therefore, any matrix
$M\in {\mathrm{Mat}}_{n\times m}(K(\{\sigma\})[t])$
can be reduced via elementary operations to Smith normal form.

Let $\Upsilon$ be a \tm module of dimension $d$, and let
$\Upsilon\in {\mathrm{Mat}}_{d\times d}(K\{\tau\})$
be the matrix corresponding to multiplication by $t.$ Then we can consider $\Upsilon$ as a matrix in
${\mathrm{Mat}}_{d\times d}(K(\{\sigma\}))$
via the embedding \eqref{embd}.
Let $\lambda_i,$ where $i=1,\dots,d,$ be the invariant factors of $\Upsilon$ obtained by diagonalizing the matrix
$tI_d-\Upsilon.$
We also have a definition analogous to the commutative case.

\begin{dfn}
 A Newton polygon of a polynomial $p(t)=\sum_{i=0}^na_it^i\in K(\{ \sigma \})[t]$ is the lower convex hull of the set of points $P_i=(i,{\mathrm{ord}}_{\sigma}(a_i)), i=0,\dots , n.$  \end{dfn}
 In the non-commutative case, the invariant factors $\lambda_i$ are not unique but the orders in $\sigma$ of their coefficients are unique once we assume that $\lambda_i$'s are monic \cite[p. 3]{Mau}. 

In this setting, A. Maurischat proved the following remarkable results:
\begin{thm}\cite[Theorem 6.6]{Mau}\label{crit}
For the \tm module $(E, \Upsilon)$ with $\Upsilon$ being represented by the matrix $D$ the following are equivalent
\begin{enumerate}
    \item[(1)] E is abelian
    \item[(2)] E is finite
    \item[(3)] the Newton polygon of the last invariant factor $\lambda_d$ of the matrix $D$ has positive slopes only.
\end{enumerate}
\end{thm}
\begin{thm}\cite[Theorem 7.2]{Mau}\label{crit1}
Let $(E,\Upsilon)$ be an abelian \tm module of dimension $d$ and $D\in {\mathrm{Mat}}_{d\times d}(K\{\tau\})$ the matrix representing $\Upsilon$ with respect to a fixed coordinate system. Let $M$ be the \tm motive of $E$. The \tm motive is pure if and only if the Newton polygon of the last invariant factor $\lambda_d$ of $D$ has exactly one edge.

In this case, the weight of $M$ equals the reciprocal of the slope of the edge.
\end{thm}

The idea of the proof of next theorem relies on reducing  the matrix $C=tI_n-\Upsilon$ to the Smith normal form. We consider the matrix $C$ in $K(\{\sigma\})[t]$- the polynomial ring over the skew field $K(\{\sigma\}).$  

 We begin with the following lemma
 \begin{lem}\label{lt}
Let $R$ be a ring and let $A=(a_{ij})\in M_n(R)$ be a lower triangular matrix. Regard $R^n$ as a right $R$-module of column vectors and consider the right $R$-linear map $A:R^n\to R^n$, $x\mapsto Ax$. We define $\operatorname{coker}(A)=AR^n\backslash R^n$.

Then $\operatorname{coker}(A)$ admits a filtration $0=M_0\subset M_1\subset\cdots\subset M_n=\operatorname{coker}(A)$ such that $M_i/M_{i-1}\cong a_{n-i+1,n-i+1}R\backslash R$ as right $R$-modules.

In particular, the composition factors of $\operatorname{coker}(A)$ depend only on the diagonal entries $a_{11},\ldots,a_{nn}$.
\end{lem}
\begin{proof}
Let $\overline e_1,\ldots,\overline e_n$ denote the images of the standard basis vectors in $\operatorname{coker}(A)$. For $i=1,\ldots,n$, let $M_i$ be the right $R$-submodule generated by $\overline e_{n-i+1},\ldots,\overline e_n$.

Put $j=n-i+1$. Since $A$ is lower triangular, its $j$-th column is
$Ae_j=e_ja_{jj}+e_{j+1}a_{j+1,j}+\cdots+e_na_{nj}$.
Since $Ae_j$ belongs to $AR^n$, its class in $\operatorname{coker}(A)$ is zero. Modulo $M_{i-1}$, this gives the relation $\overline e_j a_{jj}=0$.

Thus $M_i/M_{i-1}$ is generated by the class of $\overline e_j$ with the relation determined by $a_{jj}$. Hence
$M_i/M_{i-1}\cong a_{jj}R\backslash R=a_{n-i+1,n-i+1}R\backslash R$ as right $R$-modules.

Finally, $M_n$ is generated by $\overline e_1,\ldots,\overline e_n$, and therefore $M_n=\operatorname{coker}(A)$.
\end{proof}
\begin{lem}\label{equiv}
Let $R$ be a ring and let $A,B\in M_n(R)$. Suppose that $A$ and $B$ are equivalent, i.e. there exist $P,Q\in GL_n(R)$ such that $B=PAQ$. Then
$\operatorname{coker}(A)\cong\operatorname{coker}(B)$
as right $R$-modules.
\end{lem}

\begin{proof}
We regard $R^n$ as a right $R$-module and define
$\operatorname{coker}(A)=AR^n\backslash R^n$.
Consider the right $R$-linear automorphism
$P:R^n\to R^n$, $x\mapsto Px$.
Since $B=PAQ$, we have
$P(AR^n)=PAR^n=PAQR^n=BR^n$.
Hence $P$ induces an isomorphism
$\operatorname{coker}(A)\to\operatorname{coker}(B)$,
given by
$x+AR^n\mapsto Px+BR^n$.
\end{proof}
\begin{lem}\label{isos}
Let $\mathcal{L}$ be a skew field and let $t$ be central in $\mathcal{L}[t]$. For $\alpha,\beta\in \mathcal{L}$, consider the right $\mathcal{L}[t]$-modules
$M_\alpha=(t-\alpha)\mathcal{L}[t]\backslash \mathcal{L}[t]$
and
$M_\beta=(t-\beta)\mathcal{L}[t]\backslash \mathcal{L}[t]$.

Then $M_\alpha\cong_{\mathcal{L}[t]}M_\beta$ if and only if $\alpha$ and $\beta$ are similar in $\mathcal{L}$, i.e. if and only if there exists $c\in \mathcal{L}^\times$ such that
$\alpha=c^{-1}\beta c$.
\end{lem}
\begin{proof}
Let $\overline{1}_\alpha$ and $\overline{1}_\beta$ denote the classes of $1$ in $M_\alpha$ and $M_\beta$, respectively. Since $t$ is central and $\overline{1}_\alpha t=\overline{1}_\alpha\alpha$, every element of $M_\alpha$ is uniquely of the form $\overline{1}_\alpha c$, $c\in \mathcal{L}$. Thus $M_\alpha$ is a one-dimensional right $\mathcal{L}$-vector space. Similarly, $M_\beta$ is a one-dimensional right $\mathcal{L}$-vector space.

Suppose that $f:M_\alpha\to M_\beta$ is an isomorphism of right $\mathcal{L}[t]$-modules. Since $f$ is in particular right $\mathcal{L}$-linear, there exists $c\in\mathcal{L}K^\times$ such that
$f(\overline{1}_\alpha)=\overline{1}_\beta c$.

Since $\overline{1}_\alpha t=\overline{1}_\alpha\alpha$, we have
$f(\overline{1}_\alpha t)=f(\overline{1}_\alpha)t$.
Hence
$\overline{1}_\beta c\alpha=(\overline{1}_\beta c)t$.
Since $t$ is central,
$(\overline{1}_\beta c)t=\overline{1}_\beta\beta c$.
Therefore $c\alpha=\beta c$, and hence
$\alpha=c^{-1}\beta c$.

Conversely, suppose that $\alpha=c^{-1}\beta c$ for some $c\in \mathcal{L}^\times$. Then $c\alpha=\beta c$. Since $M_\alpha$ and $M_\beta$ are one-dimensional right $\mathcal{L}$-vector spaces, there is a unique right $\mathcal{L}$-linear map
$f:M_\alpha\to M_\beta$
such that
$f(\overline{1}_\alpha)=\overline{1}_\beta c$.

We claim that $f$ is $\mathcal{L}[t]$-linear. Indeed,
$f(\overline{1}_\alpha t)
=f(\overline{1}_\alpha\alpha)
=\overline{1}_\beta c\alpha
=\overline{1}_\beta\beta c
=(\overline{1}_\beta c)t
=f(\overline{1}_\alpha)t$.
Since $f$ is right $\mathcal{L}$-linear and commutes with the action of $t$, it is a right $\mathcal{L}[t]$-module homomorphism.

Similarly, the unique right $\mathcal{L}$-linear map
$g:M_\beta\to M_\alpha$
satisfying
$g(\overline{1}_\beta)=\overline{1}_\alpha c^{-1}$
is $\mathcal{L}[t]$-linear. Clearly $g=f^{-1}$. Hence $M_\alpha\cong_{\mathcal{L}[t]}M_\beta$.
\end{proof}

Now we are ready to prove the following
\begin{thm}\label{finitness}
Let $\Upsilon = \Upsilon\big(\underline{\psi},\underline{\delta}\big)$ be a triangular \tm module of dimension  $n$ over a perfect field $K$.
Then 
\begin{enumerate}
    \item[(i)] $\Upsilon$ is finite and, therefore, also abelian,
\item[(ii)] $\Upsilon$ is pure if and only if 
$\rk\psi_1=\dots =\rk\psi_n$
\end{enumerate}
\end{thm}
\begin{proof}
Let
\begin{equation}
C_t=tI_n-\Upsilon_t=
\label{eq_skladowe_morfizmu}
\left[
\begin{array}{ccccc}
(t-\psi_{n,\sigma})&0&\cdots&0&0\\
\xi_{n-1\times n}&(t-\psi_{n-1,\sigma})&\cdots&0&0\\
\vdots&\vdots&\ddots&\vdots&\vdots\\
\xi_{2\times n}&\xi_{2\times n-1}&\cdots&(t-\psi_{2,\sigma})&0\\
\xi_{1\times n}&\xi_{1\times n-1}&\cdots&\xi_{1\times2}&(t-\psi_{1,\sigma})
\end{array}  
\right] 
\end{equation}
We view  $C_t$ as a matrix with entries in $K(\{\sigma\})[t]$ i.e.  $C_t\in {\mathrm{Mat}}_{n\times n}K(\{\sigma\})[t].$ Since the ring $R=K(\{\sigma\})[t]$ is both left and right Euclidean there exist $P,Q\in \GL_n(R)$ such that $S:=PC_tQ$ is the  Smith normal form of $C_t,$ namely $S=\operatorname{diag}(d_1,\ldots,d_n)$, where $d_i$ divides $d_{i+1}$ from both the left and the right \cite[p.41]{jac43}. Matrices $P$ and $Q$ correspond to elementary operations.

We consider the right $R$-modules $\operatorname{coker}(C_t)$ and $\operatorname{coker}(S)$. By Lemma~\ref{equiv}, $\operatorname{coker}(C_t)\cong\operatorname{coker}(S)$ as right $R$-modules.

Since $C_t$ is lower triangular, Lemma~\ref{lt} gives a filtration of
$\operatorname{coker}(C_t)$ whose successive subquotients are
$(t-\psi_i)R\backslash R$, $i=1,\ldots,n$. Hence the composition factors of
$\operatorname{coker}(S)$ are, up to isomorphism, precisely the modules
$(t-\psi_i)R\backslash R$.

On the other hand,
$\operatorname{coker}(S)\cong\bigoplus_{i=1}^n  d_iR\backslash R$.
Factor each $d_i$ into irreducible factors. Since the modules
$(t-\psi_i)R\backslash R$ are simple, every linear factor
$(t-\alpha_i)R\backslash R$ occurring in these factorizations is isomorphic
to some $(t-\psi_j)R\backslash R$.

By Lemma~\ref{isos}, the modules
$(t-\alpha)R\backslash R$ and $(t-\beta)R\backslash R$ are isomorphic if and
only if $\alpha$ and $\beta$ are similar in $K(\{\sigma\})$. We call the
linear monic polynomials $(t-\alpha)$ and $(t-\beta)$ similar if $\alpha$ and
$\beta$ are similar in $K(\{\sigma\})$.

Consequently, the multiset of similarity classes of the linear factors occurring in the invariant factors $d_1,\ldots,d_n$ is ${[t-\psi_1],\ldots,[t-\psi_n]}$. Since $d_i$ divides $d_n$ for every $i$, all these linear factors occur, up to similarity, in the last invariant factor $\lambda_n=d_n$.

Now suppose that $\psi_i$ is a Drinfeld module. The Newton polygon of $t-\psi_{i,\sigma}$ consists of a single edge of positive slope $\lambda_i=\operatorname{rk}\psi_i$. The same holds for every linear factor in the similarity class of $t-\psi_{i,\sigma}$, since conjugation by a unit of $R$ preserves the $\sigma$-valuation of the coefficients. Hence every linear factor occurring in $\lambda_n$ has positive slope.

By Proposition~2.3 of \cite{Mau}, the Newton polygon of a product is obtained by concatenating the edges of the Newton polygons of its factors in non-decreasing order of their slopes. Consequently, all slopes of the Newton polygon of $\lambda_n$ are positive. This proves (i).

For (ii), the Newton polygon of $\lambda_n$ consists of a single edge if and only if all its linear factors have the same slope. The slope of a linear factor in the similarity class of $t-\psi_{i,\sigma}$ is $\operatorname{rk}\psi_i$. Since every invariant factor divides the last invariant factor, this is equivalent to $\operatorname{rk}\psi_1=\cdots=\operatorname{rk}\psi_n$.
\end{proof}

\begin{remark}
Recently, A. Maurischat developed an  algorithm \cite{Mau1} for finding $t$-bases for both \tm motives and \tm comotives corresponding to finite \tm modules. Theorem \ref{finitness} (i) shows that his algorithm is applicable for triangular \tm modules. 
Determining the $t$-bases of  motives and comotives is necessary for effectively finding  a rigid analytic trivialization once we know the period lattice of a \tm module
\cite[section 4]{np} (see also \cite[section 2]{hj}).
\end{remark}
\subsection{Taelman's conjecture} Let $K={\mathbb F}_q(\theta)$ and let $F$ be a finite extension of $K.$ Let ${\cal O}_F$ be an integral closure of $A$ in $F$ and $F_{\infty}=F\otimes_K K_{\infty}.$ Consider 
a \tm module  $\Upsilon$ over ${\cal O}_F$, i.e.
an ${\mathbb F}_q$-algebra homomorphism ${\Upsilon} : A\rightarrow {\mathrm{Mat}}_{d\times d}({\cal O}_F)\{\tau \} $ such that for every $\Upsilon_a=\partial_{\Upsilon}(a) +\sum_{i=1}^n\Upsilon_{a,i}{\tau}^i$ one has $(\partial_{\Upsilon}(a)-a I_d)^d=0.$

Set  $W_{\Upsilon}(F_{\infty}):= {\mathrm{Lie}}({\Upsilon})(F_{\infty})/(\partial_{\Upsilon}(\theta)-\theta I_d){\mathrm{Lie}}({\Upsilon})(F_{\infty})$ and let $$w: {\mathrm{Lie}}({\Upsilon})(F_{\infty})\rightarrow W_{\Upsilon}(F_{\infty})$$ be the projection. In \cite{an}  the authors consider the following form of Taelman's conjecture (cf. \cite{tae}).

\begin{conjecture}\label{tc}(Taelman's conjecture \cite{an}) With the above notation, suppose that ${\Upsilon}$ is  finite and uniformizable. Then there exist an element $a\in A\backslash \{0\}$ and a sub-A-module $Z\subset {\mathrm{Lie}}({\Upsilon})(F_{\infty})$ of rank 
$r:=\dim_{ K_{\infty}}W_{\Upsilon}(F_{\infty})$ such that
\begin{enumerate}
    \item[1)] $\exp_{\Upsilon}(Z)\subset {\mathrm{Lie}}({\Upsilon})({\cal O}_F)$
     \item[2)] ${\bigwedge}_A^r w(Z)=a\cdot L({\Upsilon}/{\cal O}_F)\cdot {\bigwedge}_A^r W_{\Upsilon}({\cal O}_F)$
\end{enumerate}
    \end{conjecture}
\begin{remark}
 In the above conjecture $$L({\Upsilon}/{\cal O}_F):={\prod}_{\mathfrak{p}}\frac{[{\mathrm{Lie}}({\Upsilon})({\cal O}_F/{\mathfrak{p}})]_A}{[{\Upsilon}({\cal O}_F)/{\mathfrak{p}}]_A}$$  where ${\mathfrak p}$ runs through the set of maximal ideals of ${\cal O}_F$ and the right-hand side converges in $K_{\infty}$ (cf. \cite{fa}).
\end{remark}
Since triangular \tm modules   are finite and uniformizable, it is natural to ask whether the Taelman's conjecture holds  for them.
J. Fang \cite{fa} proved that Conjecture \ref{tc} holds  for \tm modules with no nilpotence. Thus, we have the following:
\begin{cor}
    \label{thm:taelmas_conjecture}
Let ${\Upsilon}=\Upsilon\big(\underline{\psi},\underline{\delta}\big)$ be a triangular \tm module defined over ${\cal O}_F$ with no nilpotence, i.e. none of $\delta_{i,j}$ appearing in $\underline{\delta}$  has a constant term.  Then  Conjecture \ref{tc}, holds true. 
\end{cor} 

By \cite[Theorem C]{an} we also have the following:
\begin{cor}\label{thm:taelmas_conjecture2}
 Let  ${\Upsilon}=\Upsilon\big(\underline{\psi},\underline{\delta}\big)$ be a triangular \tm module defined over ${\cal O}_F$, and let $C$ be the Carlitz module. Then Conjecture \ref{tc} holds for ${\Upsilon}\otimes C/{\cal O}_F $.  
\end{cor}
\section{Weil-Barsotti formula} \label{sec:Weil-Bersotti}
In this section, we study the notion of duality for triangular $t$-modules.  
We introduce the concept of a dual biderivation and show its 
${\mathbb F}_q$-linearity.  We also establish an analogue of the classical 
Weil--Barsotti formula for triangular $t$-modules.

\subsection{Dual modules for triangular \tm modules}
Let $\Phi$ and $\Psi$ be \tm modules. Consider a biderivation $\delta\in \Der(\Phi,\Psi)$. We say that $\delta$ has no constant term if $\partial\delta=0$. The space of all biderivations without constant term will be denoted as $\Der_0(\Phi,\Psi).$ By $\Der_{in,0}(\Phi,\Psi)$ we denote the space of inner biderivations with no constant term. An important example of an $\F_q[t]-$module is 
\begin{equation}\label{def:ext_0}
    \Ext_0(\Phi,\Psi):= \Der_0(\Phi,\Psi)/\Der_{in,0}(\Phi,\Psi)
\end{equation}
for which there exists the following exact sequence of 
 $F_q[t]-$modules
$$0\lra \Ext_0(\Phi,\Psi) \lra \Ext_\tau^1(\Phi,\Psi)\lra \Ext^1\big(\Lie(E),\Lie(F)\big),$$
see \cite[Corollary 2.3]{pr}. 

Let $C$ be the Carlitz  module, i.e. $C_t=\theta+\tau$. For a  \tm module $\Phi$, we define its dual module $\Phi^\vee$, in the following manner:
\begin{equation}\label{Phi_dual}
    \Phi^\vee:=\Ext_0(\Phi,C).
\end{equation}
A dual module to a  \tm module is an  $\F_q[t]-$module, but it does not necessarily have  a structure of a \tm module, see \cite[Example 4.1]{kk24}. 
We are interested in the cases where 
 $\Phi^\vee$ is a  \tm module. This is a case where $\Phi$ is a Drinfeld module of  at least rank two \cite[Theorem 1.1]{pr}, or is a product of such Drinfeld modules \cite[Corollary 6.2]{kk04}, or is a strictly pure \tm module of rank at least two \cite[Theorem 8.4]{kk04}. 
 We will show that an analogous assertion for triangular \tm modules holds true  .  

We start with the following:
\begin{lem}\label{lem:dual_Drinfelda_jest_prosty}
    Let $\psi$ be a  Drinfeld module. Then $\psi^\vee$ is a  $\tau-$simple module. 
\end{lem}
\begin{proof}
Let $\psi=\theta+\sum\limits_{i=1}^r a_i\tau^i$. 
Recall from \cite[Section 5]{ta} that $\psi^\vee$ has the following form:
\begin{align}\label{eq:postac_psi_dual}
   \psi^\vee =  \left[ 
   \begin{array}{cccccc}
        \theta& 0 & 0 & \cdots & 0 & -\dfrac{a_1}{a_r}\tau+ \dfrac 1{a_r^{(1)}} \tau^2  \\[1pc] 
        \tau & \theta & 0 & \cdots & 0&   -\dfrac{a_2}{a_r}\tau\\[1pc]
        0& \tau & \theta &  \cdots & 0 &  -\dfrac{a_3}{a_r}\tau\\
        \vdots & \vdots & \vdots & \ddots &\vdots &\vdots  \\
        0 & 0 & 0& \cdots & \theta&   -\dfrac{a_{r-2} }{a_r}\tau\\[1pc]
        0 & 0 & 0& \cdots & \tau& \theta  -\dfrac{a_{r-1} }{a_r}\tau
   \end{array}\right].
\end{align}
Notice that $\psi^\vee\cong {\bigwedge}^{r-1}\psi. $ 
Assume that $\psi^{\vee}$ is not $\tau$-simple.
Thus there exists a nontrivial exact sequence of \tm motives:
$$0\rightarrow M_1 \rightarrow {\bigwedge}^{r-1}M_\psi\rightarrow M_2\rightarrow 0$$
where $\rk M_1, \rk M_2 < \rk{\bigwedge}^{r-1} M_\psi$ and $d(M_1), d(M_2)< d({\bigwedge}^{r-1} M_\psi).$ But by \cite[Theorem 5.1]{ta}
$w(\psi^\vee)=(r-1)/r.$ Since $\gcd (r,r-1)=1$ for 
$r=\rk\psi \geq 2,$ (\ref{weight}) yields a contradiction. Thus $\psi^\vee$ is $\tau$-simple.
\end{proof}

Notice that the concept of duality is consistent with the $\tau-$composition series.
\begin{thm}
\label{thm:dual_dla_trangulowalnego_ogolnie}
  Let $\Upsilon=\Upsilon\big(\underline{\psi},\underline{\delta}\big)$ be a triangular \tm module 
  with no nilpotence such that  $\rk\psi_i>1$ for $i=1,\dots, n$.
  \\ Then 
  \begin{itemize}
      \item[(i)] dual $\Upsilon^\vee$ is a  \tm module given by the $\tau-$co-composition series of the following form: 
      $$\Upsilon_1^\vee \twoheadrightarrow \Upsilon_2^\vee \twoheadrightarrow \cdots 
\twoheadrightarrow \Upsilon_n^\vee=\Upsilon^\vee,$$
    \item[(ii)] there exists the following exact sequence of  \tm modules  
        $$0\lra \Upsilon^\vee\lra \Ext^1_\tau(\Upsilon, C)\lra \G_a^n\lra 0$$
      \item[(iii)] Every map $f:\Upsilon\lra \widehat{\Upsilon}$  of \tm modules that
      satisfy the assumptions of this theorem induces a map
       $f^\vee:\widehat{\Upsilon}^\vee\lra \Upsilon^\vee$.
  \end{itemize}
\end{thm}
\begin{proof}

Consider the following exact sequences:
 $$\delta_l:\quad 0\lra \Upsilon_l\lra \Upsilon_{l+1}\lra \psi_{l+1}\lra 0$$
for $l=1,2,\dots, n-1$, where $\Upsilon_1=\psi_1$ is a Drinfeld module. 
 Notice that from Proposition \ref{prop:hom_vanish} it follows that  $\Hom_\tau(\Upsilon_l,C)$ vanishes for every  $l=1,2,\dots n-1$. 

Applying the functor $\Hom_\tau(-,C)$ to the sequence $\delta_1$ we obtain the following exact sequence of $\F_q[t]$-modules:
\begin{align*}
    &0\lra \Ext^1_\tau(\psi_2,C)\lra \Ext^1_\tau(\Upsilon_2,C)\lra \Ext^1_\tau(\psi_1,C)\lra 0.  
\end{align*} 
From \cite[Lemma 3.2]{gkk}, it follows that this sequence is also a sequence of \tm modules. The existence of the exact sequences of \tm modules (see \cite[Theorem 1.1]{pr}) 
\begin{align*}
    0\lra \psi_l^\vee\lra  \Ext^1_\tau(\psi_l,C) \lra \G_a\lra 0 \quad \textnormal{for} \quad l=1,2
\end{align*}
implies that there exists the following commutative diagram of \tm modules with exact rows: 
$$\xymatrix{
& 0 \ar[d] & 0 \ar[d] & 0 \ar[d] & \\ 
0\ar[r] & \psi_2^\vee\ar[r] \ar[d] & \Upsilon_2^\vee \ar[r]\ar[d] & \psi_1^\vee\ar[r] \ar[d] & 0 \\
0\ar[r]& \Ext^1_\tau(\psi_2,C) \ar[r] \ar[d]& \Ext^1_\tau(\Upsilon_2,C) \ar[r] \ar[d]& \Ext^1_\tau(\psi_1,C) \ar[r]\ar[d]& 0\\
0\ar[r] & \G_a \ar[r]\ar[d] & \G_a\times \G_a\ar[r]\ar[d] & \G_a \ar[r] \ar[d] & 0\\
& 0 & 0 & 0 & \\},$$
where the surjectivity of the middle vertical map follows from the Snake Lemma.
Therefore, we obtain the following exact sequence:

$$0\lra \Upsilon_2^\vee\lra \Ext_\tau^1(\Upsilon_2, C)\lra \G_a^2\lra 0.$$
Repeating this reasoning for the exact sequences  $\delta_l$ for $l=2,3,\dots, n-1$ we obtain the fact that  there are the following exact sequences of \tm modules:
$$\delta_l^\vee: \quad 0\lra \psi_{l+1}^\vee\lra \Upsilon_{l+1}^\vee \lra \Upsilon_{l}^\vee\lra 0$$ 
for $l=1,2,\dots, n-1$.  
 From these exact sequences, in the same way as before, we obtain the following exact sequences of \tm modules 
$$0\lra \Upsilon_{l}^\vee\lra \Ext_\tau^1(\Upsilon_{l}, C)\lra \G_a^{l}\lra 0$$
for $l=2,\dots, n$. The last sequence for  $l=n$ shows $(ii)$.  

From lemma  \ref{lem:dual_Drinfelda_jest_prosty} the modules $\psi^\vee$ are $\tau-$simple. Therefore, the sequences  $\delta_l^\vee$ for $l=1,2,\dots, n-1$ show $(i)$.

Part $(iii)$. From $(ii)$ there exist exact sequences of the following form:
\begin{align*}
    0\lra &\Upsilon^\vee\lra \Ext^1_\tau(\Upsilon, C)\lra \G_a^n\lra 0\\
0\lra &\widehat{\Upsilon}^\vee\lra \Ext^1_\tau(\widehat{\Upsilon}, C)\lra \G_a^m\lra 0.
\end{align*}

Arguing similarly as in the proof of 
 \cite[Proposition 3.2]{kk24}, we obtain that there exists  the following commutative diagram with the exact rows:
$$\xymatrix{0\ar[r] & \widehat\Upsilon^\vee \ar[d]^{ f^\vee} \ar[r]  &\Ext^1_\tau( \widehat\Upsilon, C) \ar[d]^{\overline{f} }  \ar[r] &\G_a^m \ar[d]^{\partial f}  \ar[r] & 0\\
0\ar[r] &{\Upsilon}^\vee\ar[r] &  \ar[r] \Ext^1_\tau({\Upsilon}, C) & \G_a^n \ar[r] &  0 }.$$
\end{proof}

\begin{remark}
If $\Upsilon=\Upsilon\big(\underline{\psi},\underline{\delta}\big)$ is a triangular \tm module fulfilling the assumptions of the theorem \ref{thm:dual_dla_trangulowalnego_ogolnie}, then the dual module 
     $\Upsilon^\vee$ has the following matrix form: 
         \begin{equation}
        \label{eq:postac_dla_duala}
    \Upsilon^\vee=\begin{pNiceArray}{cccccc}[margin=15pt]
 \Block[draw=black]{1-1}<\large>{\psi_n^\vee}&\Block[draw=black,fill=gray!20!]{1-5}<\large>{\delta_{n-1}^\vee} &&   & &\\[2pt] \hline 
         0&\Block[draw=black]{1-1}<\large>{\psi_{n-1}^\vee}&\Block[draw=black, fill=gray!20!]{1-4}<\large>{\delta_{n-2}^\vee}  & & &\\[2pt]
         \vdots& &\ddots & &\vdots & \\[2pt] 
         0&0 &\cdots & \Block[draw=black]{1-1}<\large>{\psi_3^\vee}&\Block[draw=black, fill=gray!20!]{1-2}<\large>{\delta_{2}^\vee} & \\[2pt] 
         0&0 &\cdots & 0& \Block[draw=black]{1-1}<\large>{\psi_2^\vee}&\Block[draw=black, fill=gray!20!]{1-1}<\large>{\delta_1^\vee} \\[2pt] 
         0&0 &\cdots & 0& 0& \Block[draw=black]{1-1}<\large>{\psi_1^\vee}
\end{pNiceArray},
    \end{equation}
\end{remark}

In what follows, the following proposition will be useful:
\begin{prop} \label{prop:hom_dual_Carlitz}
\phantom{a}
\begin{itemize}
    \item[$(i)$] If $\psi$ is a Drinfeld module such that $\rk\psi>1$ then 
        $$\Hom_\tau \big(\psi^\vee, C\big)=0.$$
    \item[$(ii)$] If $\Upsilon=\Upsilon\big(\underline{\psi},\underline{\delta}\big)$ is a triangular \tm module 
  with no nilpotence such that  $\rk\psi_i>1$ for $i=1,\dots, n$, then
  $$\Hom_\tau(\Upsilon^\vee, C)=0.$$
\end{itemize}   
\end{prop}

\begin{proof}
   Let $u\in\Hom_\tau(\psi^\vee,C)$ be a nonzero morphism, i.e.  $u\psi^\vee=C u$. Let $\psi^\vee$ be of the form \eqref{eq:postac_psi_dual}. Then  $u=[u_1,u_2,\dots, u_{r-1}]\in \Mat_{1\times r-1}\big(K\{\tau\}\big)$. The equality $u\psi^\vee=C u$ yields the following relations: 
$$ u_i\cdot \theta+ u_{i+1}\cdot \tau=(\theta+\tau)\cdot u_i\quad \textnormal{for}\quad i=1,2,\dots, r-2. $$
and 
$$ u_1\cdot \Big(-\dfrac{a_1}{a_r}\tau+ \dfrac 1{a_r^{(1)}} \tau^2 \Big) - u_2\cdot \dfrac{a_2}{a_r}\tau -\cdots- u_{r-1}\cdot \dfrac{a_{r-1}}{a_r}\tau = (\theta+\tau)\cdot u_{r-1}.  $$
Comparing the highest degree terms in the first
 $r-2$ equalities we obtain that
$$\deg_\tau u_1=\deg_\tau u_2=\cdots =\deg_\tau u_{r-1}.$$
But then by comparing the degrees in the last equality, we obtain a contradiction. Thus $\Hom_\tau(\psi^\vee,C)=0$.  This shows $(i)$.

For the proof of  $(ii)$, recall that according to  Theorem  \ref{thm:dual_dla_trangulowalnego_ogolnie} there exist the following exact sequences:
   $$\delta_l^\vee: \quad 0\lra \psi_{l+1}^\vee\lra \Upsilon_{l+1}^\vee \lra \Upsilon_{l}^\vee\lra 0\quad \textnormal{for}\quad l=1,2,\dots, n-1.$$ 
   
Acting on $\delta_1^\vee$ by the functor $\Hom_\tau(-,C)$, we obtain the six-term exact sequence \cite[Theorem 10.2]{kk04}:
\begin{align*}
    0&\lra \Hom_\tau(\psi_1^\vee , C)\lra \Hom_\tau (\Upsilon_2^\vee , C)\lra \Hom_\tau (\psi_2^\vee , C)\lra \cdots 
\end{align*}
 This implies   $\Hom_\tau (\Upsilon_2^\vee , C)=0$. Reasoning analogously for $\delta^\vee_l,$ and $l=2,3,\dots, n-2$ we obtain $\Hom_\tau(\Upsilon^\vee, C)=0.$

\end{proof}

\subsection{\tm reduction algorithm}
For the convenience of the reader, we briefly recall the \tm reduction algorithm introduced in \cite{gkk}.
Let  $\Phi$ and $\Psi$ be \tm modules of dimensions $d$ and $e$, resp. By  $E_{i\times j}$ we denote the matrix   of type $e\times d$, where the only nonzero entry is  $1$ at the place $i\times j$. 
For a matrix $N=\big[n_{i,j}\big]\in\mathrm{Mat}_{e\times d}(\mathbb{Z}_{\geq 0})$, whose entries are non-negative integers we denote:
$$\mathrm{Mat}_{e\times d}(K\{\tau\})_{< N}=\Big\{	
\big[w_{i,j}(\tau)\big]\in \mathrm{Mat}_{e\times d}(K\{\tau\})\mid \deg_{\tau}w_{i, j}(\tau)<n_{i, j}\quad \forall i,j
\Big\}$$

The \tm reduction algorithm is performed in the following steps \medskip \\
\textbf{Step 1.} \smallskip\\
We use an isomorphism of  $\F_q[t]-$modules  
$$\Ext^1_{\tau}(\Phi, \Psi)\cong \Der(\Phi, \Psi)/\Derin(\Phi, \Psi),$$
and identify every element of  $\Ext^1_{\tau}(\Phi,\Psi)$ with the corresponding (class) of biderivation $\delta:\F_q[t]\lra \mathrm{Mat}_{e\times d}(K\{\tau\})$. Since every biderivation is uniquely determined by  $\delta_t:=\delta(t)$,
we can view every element of 
  $\Ext^1_{\tau}(\Phi,\Psi)$ as the class of the matrix $\delta_t\in\mathrm{Mat}_{e\times d}(K\{\tau\})$. \medskip \\
\textbf{Step 2.} \smallskip\\
 We find a basis $\big(U_{i\times j}\big)$ of the space $\mathrm{Mat}_{e\times d}(K)$ such that the inner biderivations 
  $\delta^{(c\tau^kU_{i\times j})}$ for $c\in K$ i $k=0,1,2,\dots$  fulfill the following conditions:
 \begin{itemize}
 	\item[A.] the inner biderivations $\delta^{(c\tau^k U_{i\times j})}$ form a basis of the  
  $\F_q-$linear space $\Derin(\Phi, \Psi)$;
 	\item[B.] the matrix $\delta^{(c\tau^kU_{i\times j})}_t$ has at the $i\times j$-entry  a polynomial whose coefficient 
	at the biggest power of  $\tau$ is equal to $c$. 
 	\item[C.] for every biderivation $\delta\in \Der(\Phi,\Psi)$ there exists a sequence of inner biderivations  $\delta^{(c_l\tau^{k_l} U_{i_l\times j_l})}$ wbere $c_l\in K$ and $k_l\in\mathbb{Z}_{\geq 0}$   for $l=1,2,\dots, r$ 
	such that the reduced biderivation:
 		$$\delta^{\mathrm{reduce}}=\delta -\sum_{l=1}^s \delta^{(c_l\tau^{k_l}\cdot U_{i_l\times j_l})}\in \mathrm{Mat}_{e\times d}(K\{\tau\})_{< N}, $$
 		where
 		$N=\Big[ 
 		\deg_{\tau}\delta^{(c\tau^0U_{i\times j})}_t\Big]_{i,j}.$
 		
 	 \item[D.]  different reduced biderivations in 
 	  $\mathrm{Mat}_{e\times d}(K\{\tau\})_{< N}$ represent different elements in $\Ext^1_{\tau}(\Phi,\Psi)$. 
 			 
 \end{itemize}
Points  A., B., C. and  D. imply that there exists an isomorphism of  $\F_q-$linear spaces:
\begin{equation}
	\label{eq:algorytm_reduckji_step2}
	\Ext^1_{\tau}(\Phi, \Psi)\cong \mathrm{Mat}_{e\times d}(K\{\tau\})_{< N}.
\end{equation}
\noindent
\textbf{Step 3.} \smallskip\\
In the space $\mathrm{Mat}_{e\times d}(K\{\tau\})_{< N}$ we choose  the  set of generators of the following form:
$$E_{i\times j}c\tau^k\quad \textnormal{for}\quad c\in K,\quad\textnormal{and} \quad k=0,1,2,\dots, n_{i\times j}-1$$
\noindent
\textbf{Step 4.} \smallskip\\
 In order to  determine the structure of a \tm module on $\Ext^1_{\tau}(\Phi,\Psi)$
 we transfer the  $\F_q[t]-$module  structure from $\Ext^1_{\tau}(\Phi,\Psi)=\Der(\Phi,\Psi)/\Derin(\Phi,\Psi)$ to the space  $\mathrm{Mat}_{e\times d}(K\{\tau\})_{< N}$ via the isomorphism \eqref{eq:algorytm_reduckji_step2}.
 To achieve this it is enough to find the value of  multiplication by 
$t$  on the generators from Step 3. 
Recall that multiplication of $\delta\in \Der(\Phi,\Psi)/\Derin(\Phi,\Psi)$ by  $t$ is given by the following formula cf. \cite{pr}: 
$t*\delta = \Psi_t\cdot \delta.$
Thus in this step we determine  the values of multiplication by  $t$ on the generators chosen in Step 3:
$$t*E_{i\times j}c\tau^k=\Psi_t\cdot E_{i\times j}c\tau^k.$$
\textbf{Step 5.} \smallskip\\
If  $t*E_{i\times j}c\tau^k\notin \mathrm{Mat}_{e\times d}(K\{\tau\})_{< N}$ then we reduce the terms of too big degrees 
by means of the inner biderivations from Step 2. We  finish the process when the reduced biderivation is in $\mathrm{Mat}_{e\times d}(K\{\tau\})_{< N}$.\medskip\\
\textbf{Step 6.} \smallskip\\
We write the reduced values of $t*E_{i\times j}c\tau^k$ in the coordinate system $E_{i\times j}\tau^k$. 
\smallskip\\
\textbf{Step 7.} \smallskip\\
Since every expression of the form
 $c^{q^m}$ can be written as the evaluation of the monomial  $\tau^m$ at $\tau=c$ the coordinate vector for every value of  $t*E_{i\times j}c\tau^k$ can be written as the vector of polynomials in  $K\{\tau\}$ evaluated in  $\tau=c$. 
\medskip\\
\textbf{Final Step} \smallskip\\    
We form a matrix  $\Pi_t$, whose columns are the  vectors of polynomials computed in  Step 7. Matrix $\Pi_t$ 
determines the \tm module structure on $\Ext^1_{\tau}(\Phi,\Psi)$. \\

Note that the \tm reduction algorithm is not always executable; see \cite[Remark~3.1]{gkk}. Nevertheless, the following proposition holds.

\begin{prop}[ Proposition 3.1,\cite{gkk} ] \label{prop:step2towystarczy}
	Let $\Phi$ and $\Psi$ be \tm modules such that  $\rk\Phi>\rk\Psi$ and for which there exists a basis 
	 $\big(U_{i\times j}\big)_{i,j}$ satisfying the conditions of  \textnormal{\textbf{Step 2.}} of \tm reduction algorithm. 
	 Then the \tm reduction algorithm is executable.
\end{prop}
If the \tm module structure  on $\Ext^1_{\tau}(\Phi,\Psi)$ can be obtained from the described above  \tm reduction algorithm we will say that such a  \tm module structure  comes from  \tm reduction.

\subsection{Dual biderivation}

Our next goal is to determine those  triangular \tm modules for which the  Cartier-Nishi theorem holds true, i.e. ${\Upsilon^\vee}^\vee\cong \Upsilon$. 
In order to do this, we have to determine the exact form of the  \tm modules $\Upsilon^\vee$. 
This amounts to finding the biderivations
 $\delta_i^\vee$ for $i=1,\dots n-1$ in \eqref{eq:postac_dla_duala}.
  First, we will investigate the map $\delta_i\longmapsto\delta_i^\vee$. However, our reasoning will be applied in much more general situation.  We start with the following lemma.
 \begin{lem}\label{lem:rozczepialne_przechodza_na_rozczepialne}
      Let  
      $\delta:\, 0\lra\Psi\lra \Upsilon\lra \Phi\lra 0$ be a short exact sequence of \tm modules, given by a  $\delta$ that splits. 
     Then
     \begin{itemize}
         \item[$(i)$] for each \tm module $\zeta$ the exact sequence of $\F_q[t]-$modules 
         \begin{equation}\label{eq:ciaqg:extow:r}
             0\lra \Ext^1_\tau(\Phi,\zeta)\lra \Ext^1_\tau(\Upsilon,\zeta)\lra \Ext^1_\tau(\Psi,\zeta)\lra 0
         \end{equation}
          also splits.
         \item[$(ii)$]  the exact sequence of  $\F_q[t]-$modules
         $$0\lra \Phi^\vee\lra \Upsilon^\vee\lra \Psi^\vee\lra 0 $$  splits as well. 
     \end{itemize}
 \end{lem}
 \begin{proof}
  We start with the proof of $(i).$
  Since the sequence $\delta$ splits, there exists a matrix of skew polynomials  $U$ such that  $\delta=U\Psi-\Psi U$. 
 Applying the functor $\Hom_\tau(-,\zeta)$ to the sequence  $\delta$ by \cite[Theorem 10.2]{kk04} we obtain the exact sequence of the form 
 $$ \cdots \lra\Hom_\tau(\Psi,\zeta)\uplra{-\circ \delta} \Ext^1_\tau(\Phi,\zeta)\lra \Ext^1_\tau(\Upsilon,\zeta)\lra \Ext^1_\tau(\Psi,\zeta)\lra 0.$$
 
 Since $\delta$ is an inner biderivation, the map $-\circ\delta$ is identically zero. Therefore, we obtain the short exact sequence of the form \eqref{eq:ciaqg:extow:r}.
In the language of biderivations, it can be written as:
     $$0\lra \dfrac{\Der(\Phi,\zeta)}{\Derin(\Phi,\zeta)}\uplra{i} \dfrac{\Der(\Upsilon,\zeta)}{\Derin(\Upsilon,\zeta)}\uplra{\pi} \dfrac{\Der(\Psi,\zeta)}{\Derin(\Psi,\zeta)}\lra 0, $$
     where $i(\delta_1)=\big[-\delta_1,0\big]$ and $\pi\Big(\big[\delta_1,\delta_2\big]\Big)=-\delta_2$. We will show that  $\pi$ is a retraction.

     Let 
     $g:\Der(\Psi,\zeta)\lra \Der(\Upsilon,\zeta)$ be the $\F_q$-linear map given by the following formula $g(\delta_2)=\big[\delta_2 U, -\delta_2\big]$.
     Additionally 
      $g$ maps inner biderivations into inner biderivations. Indeed, if $\delta_2\in\Derin(\Psi,\zeta),$ i.e. $\delta_2=V\Psi-\zeta V$ for some matrix  $V$, then $g(\delta_2)=\delta_{\Upsilon,\zeta}^{\big([ VU, - V]\big)}\in \Derin(\Upsilon, \zeta)$.
     So $g$ induces homomorphism of  $\F_q[t]$-modules $\overline{g}:\dfrac{\Der(\Psi,\zeta)}{\Derin(\Psi,\zeta)}\lra \dfrac{\Der(\Upsilon,\zeta)}{\Derin(\Upsilon,\zeta)}$. 
     Moreover $\pi\circ \overline {g}=id$, so $\pi$ is a retraction and the sequence  \eqref{eq:ciaqg:extow:r} splits.

Since $g(\delta_2)\in \Der_0(\Upsilon,\zeta)$ for every
$\delta_2\in \Der_0(\Psi,\zeta)$, part~(ii) is obtained by applying part~(i) with $\zeta=C$.
  
 \end{proof}

Notice that a description of an exact sequence
 $0\lra \Psi\uplra i\Upsilon\uplra \pi \Phi\lra 0$ depends on the choice of  coordinate system in the module $\Upsilon$. In \cite{pr} Papanikolas and Ramachandran fixed the choice of coordinates in such a way that  $i=[0,1]^T$ is an embedding into the second coordinate and $\pi=[1,0]$ is a projection onto the first coordinate. 
 In our work we followed this convention.
 However, one can assume that $i=[1,0]^T$is an embedding into the first coordinate and $\pi=[0,1]$ is a projection onto the second coordinate. \footnote{Such a convention for description of extensions is common in the theory of representations of algebras, see \cite{r76,r98,mel07,km20} and many others.} Then the  \tm-module $\Upsilon$ has the following bloc form: 
\begin{equation*}
    \Upsilon=\left[\begin{array}{c|c}
        \Psi & \delta^{op} \\ \hline
         0 & \Psi
    \end{array}\right]
\end{equation*}
 In order to distinguish, this space of biderivations will be denoted as $\Der^{op}(\Phi,\Psi)$. Then the biderivations $\delta^{(U)}=U\Phi-\Psi U,$  for some matrix $U$ of skew polynomials correspond to splitting sequences. Analogously, the inner biderivations will be denoted  as  $\Derin^{op}(\Phi,\Psi)$. Then there is an  $F_q[t]-$module
 isomorphism $$\Ext^1_\tau(\Phi,\Psi)\cong \dfrac{\Der^{op}(\Phi,\Psi)}{\Derin^{op}(\Phi,\Psi)}$$
This notation will be useful further on.

 Consider an exact sequence of \tm modules of the form
$\delta:\quad 0\lra \Psi\lra \Upsilon\lra \Phi\lra 0$ and a  \tm module $\zeta$.
According to  \cite[Lemma 3.5]{gkk} if $\Hom_\tau(\Psi,\zeta)=0,$  $\Ext^1_\tau(\Psi,\zeta)$ and $\Ext^1_\tau(\Phi,\zeta)$ are \tm modules coming from the algorithm of \tm reduction then $\Ext^1_\tau(\Upsilon,\zeta)$ has the  structure of a \tm module. Moreover, from the proof of \cite[Proposition 3.2]{gkk} it follows that the  \tm module structure of  $\Upsilon$ is given by the following matrix:
$$\Pi_{\Upsilon,\zeta}= \left[ \begin{array}{c|c}
   \Pi_{\Psi,\zeta}  & \Pi_{\delta} \\ \hline
   0  & \Pi_{\Phi,\zeta}
\end{array}\right],$$
where $\Pi_{\delta}\in \Der^{op}\Big(\Ext^1_\tau(\Psi,\zeta),\Ext^1_\tau(\Phi,\zeta) \Big)$ 
and $\Pi_{*,\zeta}$ is the matrix defining the  \tm module structure of $\Ext^1_\tau(*,\zeta)$.

Define the following maps:
\begin{itemize}
    \item[$(i)$] $\Pi_{(-)}:\Der(\Phi,\Psi)\lra \Der^{op}\Big(\Ext^1_\tau(\Psi,\zeta),\Ext^1_\tau(\Phi,\zeta) \Big)$ defined by the assignment  $\delta\longmapsto \Pi_\delta,$
    \item[$(i)$] $\Pi^{op}_{(-)}:\Der^{op}(\Phi,\Psi)\lra \Der\Big(\Ext^1_\tau(\Psi,\zeta),\Ext^1_\tau(\Phi,\zeta) \Big)$ defined by the assignment $\delta^{op}\longmapsto \Pi^{op}_{\delta^{op}}.$
\end{itemize}

\begin{lem}\label{lem:odwzorwanie_Pi_jest_liniowe}

    Let $\Phi,\Psi$ and $\zeta$ be \tm modules such that  $\Hom_\tau(\Psi,\zeta)=0$ and both 
    $\Ext^1_\tau(\Psi,\zeta)$ and $\Ext^1_\tau(\Phi,\zeta)$ are \tm modules 
coming from the algorithm of \tm reduction. Then the maps $\Pi_{(-)}$ i $\Pi^{op}_{(-)}$ are $\F_q$-linear.
    \end{lem}
\begin{proof}
     We will give a proof for the map
     $\Pi_{(-)}$. The proof for $\Pi^{op}_{(-)}$ is analogous. Let $\dim \Phi=d$, $\dim\Psi=e$ and $\dim\zeta=f$.

    Let $\delta_1,\delta_2\in \Der(\Phi,\Psi)$.  First, we will show that $\Pi_{(\delta_1+\delta_2)}=\Pi_{\delta_1}+\Pi_{\delta_2}$. Denote by $\Upsilon_i$ for $i=1,2$ the middle term of the exact sequence corresponding to the biderivation  $\delta_i$ and denote by  $\Upsilon_+$ the middle term of the exact sequence corresponding to  $\delta_1+\delta_2$.
    Since both \tm module structures for  $\Ext^1_\tau(\Psi,\zeta)$ and $\Ext^1_\tau(\Phi,\zeta)$ come from the algorithm of \tm reduction there exist the following bases:
    $$U_{i\times j}^{\Phi,\zeta}\in \Mat_{f\times d}\big(K\{\tau\}\big) \quad \textnormal{and}\quad U_{i\times l}^{\Psi,\zeta}\in \Mat_{f\times e}\big(K\{\tau\}\big) $$
    which fulfill   \textbf{step 2.} of the  \tm reduction algorithm. From the proof of  \cite[Lemma 3.5]{gkk}, it follows that the basis 
    $$\Big[ U_{i\times j}^{\Phi,\zeta} \mid 0\Big],\quad \Big[0\mid  U_{i\times j}^{\Psi,\zeta} \Big]$$ fulfills   \textbf{step 2.} of the  \tm reduction algorithm. Therefore, the spaces
    $\Ext^1_\tau(\Upsilon_*, \zeta)$ are \tm modules for $*=1,2,+$. 
    Now, the argument is similar to that of 
    \cite[Proposition 3.2]{gkk}. Let  
    $\delta_*^{(U)}$ denote the inner biderivation belonging to $\Derin(\Upsilon_*,\zeta)$ for $*=1,2,+$.
    \begin{align*}
        \delta_*^{\big(  \mu\tau^k\big[ U_{i\times j}^{\Phi,\zeta} \mid 0\big]\big)} &=\Big[\delta^{\big(  \mu\tau^k  U_{i\times j}^{\Phi,\zeta}\big)}_{\Phi,\zeta} \mid 0\Big] \quad \textnormal{for}\quad *=1,2,+\\
        \delta_*^{\big(  \mu\tau^k \big[0 \mid U_{i\times l}^{\Psi,\zeta} \big]\big)} &=\Big[ \mu\tau^k 
        U_{i\times l}^{\Psi,\zeta}\cdot \delta_* \mid \delta^{\big(  \mu\tau^k  U_{i\times l}^{\Psi,\zeta}\big)}_{\Psi,\zeta}  \Big] \quad \textnormal{for}\quad *=1,2\\
        \delta_+^{\big(  \mu\tau^k \big[0 \mid U_{i\times l}^{\Psi,\zeta} \big]\big)} &=\Big[ \mu\tau^k 
        U_{i\times l}^{\Psi,\zeta}\cdot (\delta_1+\delta_2) \mid \delta^{\big(  \mu\tau^k  U_{i\times l}^{\Psi,\zeta}\big)}_{\Psi,\zeta}  \Big],
    \end{align*}
    where $\delta^{(  U )}_{\star,\zeta} \in \Derin(\star,\zeta)$. 
    Since we are interested only in the terms
    $\Pi_{\delta_1}, \Pi_{\delta_2}$ and $\Pi_{\delta_1+\delta_2}$ we compute the action
    \begin{align*}
        t*\big[ 0 \mid E_{i\times l }c\tau^k\big]= \zeta_t\cdot \big[ 0 \mid E_{i\times l }c\tau^k\big] = \big[ 0 \mid \zeta_t\cdot E_{i\times l }c\tau^k\big].
    \end{align*}
    In order to determine the 
     \tm module structure on  $\Ext^1_\tau(\Upsilon_*, \zeta),$ reducing by means of inner biderivations, we have to ensure that  
    $$\big[ 0 \mid \zeta_t\cdot E_{i\times l }c\tau^k\big]\in \Mat_{f\times d}\big(K\{\tau\}\big)_{<N} \times \Mat_{f\times e}\big(K\{\tau\}\big)_{<M},$$
    where $N=\Big[\deg_{\tau} \delta^{\big(  \mu\tau^{0}  U_{i\times j}^{\Phi,\zeta}\big)} \Big]_{i\times j}$ and $M=\Big[\deg_{\tau} \delta^{\big(  \mu\tau^{0}  U_{i\times l}^{\Psi,\zeta}\big)} \Big]_{i\times l}$. This means that the degrees of the terms in the first bloc at the  $i\times j$-th position are smaller than  $\deg_\tau \delta^{\big(  \mu\tau^{0}  U_{i\times j}^{\Phi,\zeta}\big)} $, and in the second bloc at the  $i\times l$-th position are smaller than
    $\deg_\tau \delta^{\big(  \mu\tau^{0}  U_{i\times l}^{\Psi,\zeta}\big)}$.
    Since the basis $U_{i\times l}^{\Psi,\zeta}$ fulfills \textbf{step 2.} of the  \tm reduction algorithm, there exists a sequence of inner biderivations of the form
    $$ \delta^{\big(  \mu_s\tau^{k_s}  U_{i_s\times l_s}^{\Psi,\zeta}\big)}_{\Psi,\zeta}\quad \textnormal{for}\quad  s\in J$$
    such that 
    $$\zeta_t\cdot E_{i\times l} c\tau^k -\sum_{s\in J} \delta^{\big(  \mu_s\tau^{k_s}  U_{i_s\times l_s}^{\Psi,\zeta}\big)}_{\Psi,\zeta}\in  \Mat_{f\times e}\big(K\{\tau\}\big)_{<M}.$$ 
    Therefore, in all three cases we will reduce by means of a sequence of  biderivations:
    $$\delta_*^{\big(  \mu_s\tau^{k_s} \big[0 \mid U_{i_s\times l_s}^{\Psi,\zeta} \big]\big)}\quad\textnormal{for}\quad s\in J\quad\textnormal{and}\quad *=1,2,+.$$
    For $\Pi_{\delta_i},\,\, i=1,2$, we obtain 
    $$t*\big[ 0 \mid E_{i\times l} c\tau^k\big]= \Big[ -\sum_{s\in J}\mu_s\tau^{k_s}U_{i_s\times l_s}\cdot \delta_i \Big|\quad \textnormal{reduced part}  \Big].$$
    For  $i=1,2$ denote $w_i(\tau):=-\sum_{s\in J}\mu_s\tau^{k_s}U_{i_s\times l_s}^{\Psi,\zeta}\cdot \delta_i$. Since the basis $U_{i\times j}^{\Phi,\zeta}$ fulfills \textbf{step 2.} of the \tm reduction algorithm, there exists a sequence of  inner biderivations of the form 
    $$\delta_{\Phi,\zeta}^{\big(\mu_r\tau^{k_r}U_{i_r\times j_r}^{\Phi,\zeta} \big) }\quad r\in I_i\quad \textnormal{for}\quad i=1,2$$  
    such that 
    $$w_i(\tau)-\sum_{r\in I_i}\delta_{\Phi,\zeta}^{\big(\mu_r\tau^{k_r}U_{i_r\times j_r}^{\Phi,\zeta} \big) }\in \Mat_{f\times d}\big(K\{\tau\}\big)_{<N}.$$
    Therefore, in the reduction process, we use  inner biderivations of the following form: 
    $$\delta_{i}^{\big(\mu_r\tau^{k_r} \big[U_{i_r\times j_r}^{\Phi,\zeta} \mid 0 \big]\big) }\quad r\in I_i\quad \textnormal{for}\quad i=1,2.$$  
   In the case of   $\Pi_{\delta_1+\delta_2}$ we have
    \begin{align*}
        t*\big[ 0 \mid E_{i\times l} c\tau^k\big] &= \Big[ -\sum_{s\in J}\mu_s\tau^{k_s}U_{i_s\times l_s}\cdot (\delta_1+\delta_2) \Big|\quad \textnormal{reduced part}  \Big] \\
        &=   \Big[ w_1(\tau)+w_2(\tau) \Big|\quad \textnormal{reduced part}  \Big].
    \end{align*}
    Then for the reduction we can use the former sequences
     
    $$\delta_{1}^{\big(\mu_r\tau^{k_r} \big[U_{i_r\times j_r}^{\Phi,\zeta} \mid 0 \big]\big) }\quad r\in I_1\quad \textnormal{and}\quad \delta_{2}^{\big(\mu_r\tau^{k_r} \big[U_{i_r\times j_r}^{\Phi,\zeta} \mid 0 \big]\big) }\quad r\in I_2.$$
    Since for $\Pi_{\delta_1+\delta_2}$ we used the same inner biderivations as for   $\Pi_{\delta_1}$ and $\Pi_{\delta_2}$, we obtain the same reduced form. This shows that   $\Pi_{\delta_1+\delta_2}=\Pi_{\delta_1}+\Pi_{\delta_2}$.

In order to show that  $x\cdot \Pi_\delta=\Pi_{x\cdot \delta}$ for $x\in \F_q$, we can proceed similarly. First, we determine the inner biderivations for $\Pi_\delta$ and 
$\Pi_{x\cdot \delta}$. Then we show that for computation of the action 
$ t*\big[ 0 \mid E_{i\times l} c\tau^k\big]$, 
we can use the same inner biderivations. We leave  the details to the reader.

\end{proof}

From Lemmas \ref{lem:rozczepialne_przechodza_na_rozczepialne} and \ref{lem:odwzorwanie_Pi_jest_liniowe} we obtain the following theorem.
\begin{thm}\label{thm:homomorfizm_biderywacji}
    Let $\Phi,\Psi$ and $\zeta$ be \tm modules such that $\Hom_\tau(\Psi,\zeta)=0$ and both 
    $\Ext^1_\tau(\Psi,\zeta)$ and $\Ext^1_\tau(\Phi,\zeta)$ are \tm modules coming from the  \tm reduction algorithm.
    \begin{itemize}
        \item[(i)] Then the following maps: $$\Pi_{(-)},\Pi^{op}_{(-)}: \Ext^1_\tau(\Phi,\Psi)\lra \Ext_\tau^1\Big(\Ext^1_\tau(\Psi,\zeta),\Ext^1_\tau(\Phi,\zeta) \Big)$$
    are $\F_q$-linear.
        \item[(ii)] If the \tm modules $\Phi,\Psi$ and $\zeta$ have no nilpotency then the above maps induce the following  $\F_q-$linear maps:
        $$\Pi_{(-)},\Pi^{op}_{(-)}: \Ext_{0}(\Phi,\Psi)\lra \Ext_0\Big(\Ext_0(\Psi,\zeta),\Ext_0(\Phi,\zeta) \Big)$$
    \end{itemize}
\end{thm}

\begin{proof}
From additivity of  $\Pi_{(-)}$ and $\Pi^{op}_{(-)}$ and the fact that they map inner biderivations into inner biderivations, it follows that they induce the maps 
    $\Ext^1_\tau(\Phi,\Psi) \rightarrow \Ext_\tau^1\Big(\Ext^1_\tau(\Psi,\zeta),\Ext^1_\tau(\Phi,\zeta) \Big)$. Linearity of these maps is the assertion of Lemma  \ref{lem:odwzorwanie_Pi_jest_liniowe} which shows $(i)$. 
   We will prove   $(ii)$ for the map  $\Pi_{(-)}$. The argument for  $\Pi^{op}_{(-)}$ is analogous. It is enough to show that if $\delta\in \Der(\Phi,\Psi)$ fulfills $\partial\delta=0$, then also $\partial\Pi_{(\delta)}=0.$ But this follows directly from the   \tm reduction algorithm, since if both $\Upsilon$ and  $\zeta$ are \tm modules with no nilpotence, then  the \tm module  $\Ext^1_\tau(\Upsilon, \zeta)$ also has no nilpotence.   
\end{proof}

\begin{cor}\label{cor:odwzorowania_pi_i_pi_op}
Let $\Phi$ and $\Psi$ be \tm modules with no nilpotence such that $\Hom_\tau(\Psi,C)=0$ and 
 both   $\Ext^1_\tau(\Psi,C)$ and $\Ext^1_\tau(\Phi,C)$ are \tm modules coming from the \tm reduction algorithm.
    Then the maps $$\Pi_{(-)},\Pi^{op}_{(-)}: \Ext_0(\Phi,\Psi)\lra \Ext_0\Big(\Psi^\vee, \Phi^\vee \Big)$$
    are $\F_q-$linear.
    \end{cor}
For Corollary \ref{cor:odwzorowania_pi_i_pi_op}
we will write  $\delta^\vee$ instead of  $\Pi_\delta$ and call it a biderivation dual to  $\delta$. If  $\Psi^\vee$ and $\Phi^\vee$  again fulfill the assumptions of  Corollary \ref{cor:odwzorowania_pi_i_pi_op} then we can consider the following composition  of maps
$$\Ext_0(\Phi,\Psi)\uplra{\Pi_{(-)}} \Ext_0\Big(\Psi^\vee, \Phi^\vee \Big) \uplra{\Pi^{op}_{(-)}} \Ext_0\Big( {\Phi^{\vee}}^\vee, {\Psi^{\vee}}^\vee \Big),$$
which will be denoted as $$\delta\longmapsto {\delta^\vee}^\vee.$$

For finding an exact formula for a dual module to a triangular \tm module we need the following corollary.
Applying this corollary considerably simplifies the procedure of finding the aforementioned dual module.
   
.

\begin{cor}\label{cor:dual_dla_biderwaycji}
 Let $\Upsilon = \Upsilon\big(\underline{\psi},\underline{\delta}\big)$ be a triangular \tm module with no nilpotence of the form \eqref{modul_triangulowalny_macierz} where $\rk\psi_i>1$.  Then 
    \begin{itemize}
        \item[(i)] the correspondence   $\delta_i\longmapsto \delta_i^\vee$ is $\F_q-$linear.\smallskip
        \item[(ii)] the correspondence $\delta_i\longmapsto {\delta_i^\vee}^\vee$ is $\F_q-$linear. 
    \end{itemize}
\end{cor}

\begin{proof} Assume that  $\Upsilon$ has a $\tau-$composition series of the form \eqref{ciag_kompozycyjny}. 
Then every pair of \tm modules $\psi_{j+1}$ and $\Upsilon_j$ for $j=1,\dots n-1$ meets the assumptions of Corollary \ref{cor:odwzorowania_pi_i_pi_op}. From Proposition \ref{prop:hom_vanish} we obtain  that $\Hom_\tau(\Upsilon,C)$ vanishes since $\rk\psi_i>1$ for all $i=1,2,\dots, n$. Moreover, by \cite[Corollary 5.3.]{gkk} we see that both $\Ext^1_\tau(\psi_{j+1},C)$ and $\Ext^1_\tau(\Upsilon_j,C)$ are \tm modules coming from the  \tm reduction algorithm.
 Then part $(i)$ follows by induction applied to Corollary \ref{cor:odwzorowania_pi_i_pi_op}.

By Theorem \ref{thm:dual_dla_trangulowalnego_ogolnie} the \tm module $\Upsilon^\vee$ has a $\tau-$ cocomposition series of the form  \eqref{coco}. Now we will show that the pairs
$\Upsilon_j^\vee,$  $\psi_{j+1}^\vee$  fulfill the assumptions of Corollary \ref{cor:odwzorowania_pi_i_pi_op}.

By Proposition \ref{prop:hom_dual_Carlitz} the space $\Hom_\tau(\psi_{j+1}^\vee,C)$ vanishes. So the fact that  $\Ext^1_\tau(\psi_{j+1}^\vee,C)$ and $\Ext^1_\tau(\Upsilon_j^\vee,C)$ are \tm modules coming from Proposition \ref{prop:doubledualjesttmodulem} which will be proven in the next section.  
The proof of
$(ii)$ again then follows    by induction  applied to Corollary \ref{cor:odwzorowania_pi_i_pi_op}.
\end{proof}

\begin{remark}\label{rem:jak_liczyc_duala}Notice an important fact that 
     $\F_q-$linearity is valid only for a specific biderivation $\delta_i$, not for different biderivations.  
   This means that for a determination of a dual module to that of the following form: 
    \begin{align*} 
        \Upsilon=
        \left[\begin{array}{ccc}
            \theta+\tau^3 & 0 &0\\
             \tau^4+\tau^5 & \theta+\tau^4 &0\\
             \tau+\tau^2& \tau^5+\tau^6  &\theta+\tau^2
        \end{array}\right],
    \end{align*}
    we have to proceed as follows. First, we compute a dual module to  the following \tm module:
    \begin{align*}
        \left[\begin{array}{cc}
              \theta+\tau^4 &0\\
              \tau^5+\tau^6  &\theta+\tau^2
        \end{array}\right]
    \end{align*}
by applying  additivity $\big(\tau^5+\tau^6\big)^\vee=\big(\tau^{5}\big)^\vee+\big(\tau^6\big)^\vee$.  
Then we compute the dual module to  $\Upsilon$ by determining the dual modules for the following \tm modules:
$$
        \left[\begin{array}{ccc}
            \theta+\tau^3 & 0 &0\\
             \tau^4 & \theta+\tau^4 &0\\
             \tau& \tau^5+\tau^6  &\theta+\tau^2
        \end{array}\right]\quad \textnormal{and} \quad 
        \left[\begin{array}{ccc}
            \theta+\tau^3 & 0 &0\\
             \tau^5 & \theta+\tau^4 &0\\
             \tau^2& \tau^5+\tau^6  &\theta+\tau^2
        \end{array}\right]
$$
where we use the following  $\F_q-$linearity:
$$ \left(\begin{array}{c}
           \tau^4+\tau^5 \\
             \tau+\tau^2
        \end{array}\right)^\vee=
        \left(\begin{array}{c}
             \tau^4 \\
             \tau
        \end{array}\right)^\vee+\left(\begin{array}{c}
             \tau^5 \\
             \tau^2
        \end{array}\right)^\vee$$
\end{remark}

 \subsection{Explicit formulas for dual modules and dual biderivations}
 In the next section, we address the problem relating
 to which triangular \tm modules the Cartier-Nishi  theorem holds true for, i.e. whether a double dual  of a triangular \tm module  is isomorphic to  an initial \tm module.
 It turns out that for  \tm modules with nilpotence the answer  is negative. 

We establish the Cartier-Nishi Theorem  for \tm modules allowing for LD-biderivations. In order to prove this theorem,
we need to find exact forms of dual biderivations for this case. This will again be achieved  by applying the \tm reduction algorithm from \cite{gkk}. We start
with the following lemma.

\begin{lem}\label{lem:iso_dla_duala}
     Let $\Upsilon=\Upsilon\big(\underline{\psi},\underline{\delta}\big)$ be a triangular \tm module of dimension $n$ with no nilpotence and $r_j:=\rk\psi_{j}>1$ for $j=1,2,\dots,n$.
     Then there is an $F_q$-linear isomophism
     \begin{equation}
        \label{eq:iso:dowod}
        \Upsilon^\vee\cong \Big[ K\{\tau \}_{\langle 1, r_n )}, \cdots, 
    K\{\tau \}_{\langle 1, r_1 )} \Big],
    \end{equation}
     where $K\{\tau \}_{\langle 1, r_i )}$ denotes the space of polynomials with monomials of degrees belonging to the interval   $\langle 1, r_i )$. 
\end{lem}

\begin{proof}
    Recall that $\Upsilon^\vee=\Ext_0(\Upsilon, C)$
    is isomorphic to  $\Der_0(\Upsilon, C)/\Der_{in,0}(\Upsilon, C)$ as an $\F_q[t]-$module. Because of an isomorphism  
    $\Der(\Upsilon, C)\cong \Mat_{1\times n}(K\{\tau \})$ every biderivation $\delta\in \Der(\Upsilon, C)$ can be regarded as a $1\times n$ matrix of twisted polynomials.
      
      In what follows, to simplify  the notation, we will enumerate columns from  right to  left. According to this convention, $E_i$ will denote the matrix of dimensions  $1\times n$, where  entry at position   $1\times n-i+1$  equals $1$ and all  other entries are $0$.  
    
    For the proof, we use the following inner biderivations from $\Derin(\Upsilon, C)$:
     \begin{align}\label{eq:biderywacje_wewnetrze_trojkatny_modul}
       \delta^{(\mu \tau^k E_n)}_{\Upsilon,C} &= \Big[ \delta^{(\mu\tau^k)}_{\psi_n,C}\  \Big|\  0\  \Big|\  \cdots\  \Big|\  0\   \Big] \\ \nonumber
       \delta^{(\mu \tau^k E_{n-1})}_{\Upsilon,C} &= \Bigg[\ \sum_{l=1}^{\deg \delta_{n-1}} \mu\cdot d_{n-1\times n,l}^{(k)}\tau^{k+l}\  \Big| \ \delta^{(\mu\tau^k)}_{\psi_{n-1},C} \ \Big|\  0\  \Big|\  \cdots\  \Big|\  0\   \Bigg] \\ \nonumber
       \delta^{(\mu \tau^k E_{n-2})}_{\Upsilon,C} &= \Bigg[ \ \sum_{l=1}^{\deg \delta_{n-1}} \mu\cdot d_{n-2\times n,l}^{(k)}\tau^{k+l}\  \Big|\  \sum_{l=1}^{\deg \delta_{n-2}} \mu\cdot d_{n-2\times n-1,l}^{(k)}\tau^{k+l} \  \Big|\\\nonumber
       &\qquad \qquad \qquad \qquad \qquad \qquad  \Big|\   \delta^{(\mu\tau^k)}_{\psi_{n-2},C}\  \Big|\ 0\  \Big|\  \cdots\  \Big|\  0\   \Bigg] \\ \nonumber
        &\vdots  \\ \nonumber
        \end{align}
        \begin{align*}
        \delta^{(\mu \tau^k E_1)}_{\Upsilon,C} &=  \Bigg[\ \sum_{l=1}^{\deg \delta_{n-1}} \mu\cdot d_{1\times n,l}^{(k)}\tau^{k+l}\  \Big|\ \sum_{l=1}^{\deg \delta_{n-2}} \mu\cdot d_{1\times n-1,l}^{(k)}\tau^{k+l} \ \Big|\  \cdots\ \\ \nonumber
        &\qquad \qquad \qquad \qquad  \qquad \cdots\Big|\ 
       \sum_{l=1}^{\deg \delta_{1}}\mu\cdot  d_{1\times 2,l}^ {(k)}\tau^{k+l}\  \Big|\  \delta^{(\mu\tau^k)}_{\psi_1,C}\   \Bigg],  
    \end{align*}
     where  $\delta^{(\mu\tau^k)}_{\psi_l,C}=\mu\tau^k\cdot\psi_l-C\cdot \mu\tau^k\in\Derin(\psi_l,C)$ for $l=1,2,\dots, n$. Since $\Upsilon$ has no 
    nilpotence,  then all these biderivations belong to $\Der_0(\Upsilon, C)$. Notice that the inner biderivation
    $\delta^{(\mu\tau^kE_s)}$ will be used for reduction of    terms at the $s-$th position counted from  right to left. Similarly, every biderivation   
    $\delta\in\Der_0(\Upsilon, C)$ will be reduced from the right to the left.  
    Therefore, starting with an arbitrary  biderivation $\delta\in\Der_0(\Upsilon, C)$ one can, by means 
    of the inner biderivations $\delta^{(\mu\tau^kE_s)}$, reduce $\delta$ in such a way that the polynomial at the $s$-th  place (counted from  right to left) in the biderivation $\delta$ after reduction has degree less than   $\rk\psi_s$.
\end{proof}

From Lemma \ref{lem:iso_dla_duala} it follows that  the basis  $E_n$,\dots,  $E_2$, $E_1$ satisfies \textbf{step 2.} of the \tm reduction algorithm. Thus $\Ext^1_\tau\big(  \Upsilon ,C \big)$ has a \tm module structure coming from the   \tm reduction algorithm.
We also have the following proposition.

\begin{prop}\label{cor:dual_bez_nilpotentnosci}
     Let $\Upsilon=\big(\underline{\psi},\underline{\delta}\big)$ be a triangular \tm module with no nilpotence, such that $\rk\psi_l>1$ for 
    $l=1,2,\dots, n$. Then \tm module $\Upsilon^\vee$ has no nilpotence. 
\end{prop}
\begin{proof}
Notice that if  $\Upsilon$ has no nilpotence then the inner biderivations \eqref{eq:biderywacje_wewnetrze_trojkatny_modul} satisfy the conditions $\partial  \delta^{(\mu \tau^k E_j)}=0$ for $j=1,2,\dots n$. Therefore the \tm module $\Ext^1_\tau\big(  \Upsilon,C \big)$ has no nilpotence as well as the \tm module  $\Upsilon^\vee\subset \Ext^1_\tau\big(  \Upsilon ,C \big)$. 
\end{proof}

\begin{remark}
In the process of determining  a dual module as well as a double dual module, we will use Corollary \ref{cor:dual_dla_biderwaycji} as in  Remark \ref{rem:jak_liczyc_duala}. So we need to ensure that if a \tm module $\Upsilon$  allows for LD-biderivations (c.f. Definition \ref{twocl}), then auxiliary  \tm modules coming from monomials of the biderivations for  $\Upsilon$ also have  this property. And indeed this is true. The converse is also true once we limit ourselves to the summation of monomials in one column, i.e. those within one biderivation.  
\end{remark}

In the proposition below we determine an explicit formula for the dual module in the 
case of a \tm module allowing for LD-biderivations.

\begin{prop}\label{prop:dual_postac}
    Let $\Upsilon=\Upsilon\big(\underline{\psi},\underline{\delta}\big)$ be a triangular \tm module of dimension $n$, with no nilpotence and $r_j:=\rk\psi_{j}>1$ for $j=1,2,\dots,n$. If $\Upsilon$ allows for LD-biderivations, then
    $\delta^\vee_{i}$ has the following bloc matrix form over $K(\Upsilon)$ 
    \begin{equation}
        \label{eq:delta_vee_block_matrix_form}
         \delta^\vee_{i}=\Big[ D_{i+1,i} \mid\dots \mid D_{i+1,2} \mid  D_{i+1,1}\Big],
    \end{equation}
    where $D_{i+1, j}$ is a matrix of  size $r_{i}-1\times r_{j}-1$
     for $j=1,2\dots, i$.\\
     Moreover, if the biderivations $\delta_{i}$ are of the form \eqref{eq:wyrazy_dla_delta_i}, then  
        \begin{equation}
            \label{eq:postac_delta_vee_latwiejsza}
            D_{i+1, j}=-\dfrac{1}{a_{j,r_j}}
            \left[\begin{array}{cccc}
             0&\cdots& 0 & d_{j\times i+1,1}  \\
             0&\cdots& 0 & d_{j\times i+1,2} \\
             \vdots & \ddots & \vdots & \vdots \\
             0& \cdots & 0&  d_{j\times i+1,\deg\delta_{i+1}} \\
             0&\cdots  & 0 & 0 \\
             \vdots & \ddots & \vdots & \vdots \\
             0&\cdots & 0 & 0 \\
        \end{array}\right]\tau,  
        \end{equation} 
        where $a_{j,r_j}$ is the leading term of $\psi_j$.
\end{prop}
\begin{proof}
Since $\Upsilon$ allows for LD-biderivations, after changing to the field
$K(\Upsilon)$, we may assume that $\Upsilon$ satisfies
condition~\eqref{eq:warunek_LD}.

Determining  the dual biderivations 
     $\delta_i^\vee$ for $i=1,2,\cdots n-1$ will be done inductively. In the proof we will use the inner biderivations of the form  \eqref{eq:biderywacje_wewnetrze_trojkatny_modul}.
    
    The starting point of induction is determination  of $\delta_1^\vee$. Let
    $\Upsilon_2=\left[\begin{array}{c|c}
         \psi_2&0  \\ 
         \hline
         \delta_1&\psi_1 
    \end{array}\right],$ where $\rk\psi_i=r_i>1$ for $i=1,2$. From Corollary \ref{cor:dual_dla_biderwaycji} it follows that the correspondence  $\delta_1\longmapsto\delta_1^\vee$ is $\F_q-$linear. Therefore it is enough to find the dual biderivation for generators of the form $d\tau^h$ for $1\leq h< r_2$. 
     Lemma \ref{lem:iso_dla_duala} yields an isomorphism of $\F_q-$linear spaces:
    $$\Upsilon^\vee\cong  V_2=\Big[ K\big\{ \tau \big\}_{\langle 1, r_2 )}, K\big\{ \tau \big\}_{\langle 1, r_1 )} \Big].$$
    Recall that, following the convention from Lemma \ref{lem:iso_dla_duala}, we denote 
    $E_2=[1,0]$ and $E_1=[0,1]$ since we are counting terms from  left to  right.

    After  choosing  the following coordinate system for $V_2$:
         $$\Big(E_2 \tau^k \Big)_{k=1}^{r_2-1}, \Big(E_1 \tau^k \Big)_{k=1}^{r_1-1},$$ 
    we can determine the  \tm module structure on  $\Upsilon^\vee$. To do this, one has to compute the results of the action of $t$ on the generators  $c\tau^k E_i$ for $i=2,1$, and express the  obtained results in the chosen coordinate system.  
    
    Since we are only interested in determination of 
     $\delta^\vee_i$ we will compute the action of  $t$ solely on these generators, which have influence on the  form of $\delta^\vee_i$.

    If  $k\in\{1,\dots, r_1-2\}$, then
    $$t*c\tau^k E_{1}= (\theta+\tau)c\tau^k E_{1} = \Big(\theta c\tau^k +c^{(1)} \tau^{k+1} \Big)E_{1}\in V_2,$$
   then  
    $t*c\tau^k E_{1}$ in the fixed basis has the following coordinates:
    $$\Big[\  \boldsymbol{0} \ \Big| \ 
    \podwzorem{0,\dots, 0}{k-1\  \textnormal{terms} }, \theta , \tau , 0, \dots, 0\ \Big],$$
    where $\boldsymbol{0}$ denotes the zero row of appropriate size.
    This shows  that all  terms except those in the last column of the matrix $\delta_1^\vee$ are zero. 

    Now consider  the last column of
    $\delta_1^\vee$. The terms there result from the action of  $t$ on the generator $c\tau^{r_1-1}E_{1}$. Then
    \begin{align*}
        t*c\tau^{r_1-1}E_{1}&= (\theta+\tau)c\tau^{r_1-1}E_{1}=\Big(\theta c\tau^{r_1-1} + c^{(1)}\tau^{r_1}\Big) E_{1}\notin V_2. 
    \end{align*}
    Since the degree of $\theta c\tau^{r_1-1} + c^{(1)}\tau^{r_1}$ equals 
    $r_1$ we reduce it by means of the following inner biderivation
    $$\delta^{( \mu\tau^0 E_{1})}_{\Upsilon,C}=\Big[\, \mu d\tau^h\, \big|\, \delta^{(\mu\tau^0)}_{\psi_1,C}     \,\Big ] \quad\textnormal{for}\quad 
    \mu =\dfrac{c^{(1)}}{a_{1, r_{1} }}.$$ 
    As a result one gets  the following equality:
    \begin{align*}
        t*c\tau^{r_1-1}E_{1}&= \Bigg[\ -\dfrac{c^{(1)}}{a_{1, r_{1} }}\cdot  d \tau^{h}\ \Big|\  \textnormal{terms for } \psi_{1}^\vee\  \Bigg]
    \end{align*}

    Since $h< r_1$ the value $t*c\tau^{r_1-1}E_{1}$ can be expressed in the aforementioned fixed basis. We obtain
    \begin{align*}
         t*c\tau^{r_1-1}E_{1}&=\Bigg[\ 0,\dots, 0,\podwzorem{-\dfrac{d }{a_{1, r_{1} }} }{h-th \textnormal{ place}}\tau, 0, \dots, 
       0 \ \Big| \textnormal{ terms for } \psi_{1}^\vee\ \Bigg]_{\big|_{\tau=c} }.
    \end{align*}
    Thus 
    $$\Big(d\tau^h\Big)^\vee = -\dfrac{d }{a_{1, r_{1} }}E_{h\times r_1-1} \tau,$$
    where the matrix  $E_{a\times b}$ has the only nonzero entry equal to one at the position  $a\times b$. 
     
Then for
$\delta_{1}=\sum\limits_{k=0}^{\deg\delta_1}d_{1\times 2,k}\tau^k$ we obtain 
    \begin{align*}
        \delta_{1}^\vee=\Big(\sum\limits_{k=0}^{\deg\delta_1}d_{1\times 2,k}\tau^k\Big)^\vee=
        \sum\limits_{k=0}^{\deg\delta_1} \Big( d_{1\times 2,k}\tau^k\Big)^\vee=
        \sum\limits_{k=0}^{\deg\delta_1} -\dfrac{d_{1\times 2,k}}{a_{1, r_{1} }}E_{k\times r_1-1} \tau.
    \end{align*}
   This finishes the determination  of $\delta^\vee_1$ and the first step of induction.

    Assume that  the consecutive dual biderivations
      $\delta_l^\vee$ for
       \tm modules $\Upsilon_{l+1}$ where  $l=1,2,\dots, n-2$ have been determined and that they satisfy the following conditions: 
       
      \begin{itemize}
          \item[\textbf{a.}] every $\delta_l^\vee$ is of the form \eqref{eq:postac_delta_vee_latwiejsza},
          \item[\textbf{b.}] in the process of determining   $\delta_l^\vee$, all the values  $t*c\tau^kE_j$ gave the  columns with entries equal to zero for  $1\leq k<r_j-1$ and $j=1,2,\dots l-1$,
          \item[\textbf{c.}] in the process of determining   $\delta_l^\vee$ the values $t*c\tau^{r_j-1}E_j$ were obtained after a single reduction by the inner  biderivation 
     $\delta^{( \mu\tau^0 E_{j})}_{\Upsilon_l,C}$ for $\mu =\dfrac{c^{(1)}}{a_{j, r_{j} }}$
     and $j=1,2,\dots, l-1$.
      \end{itemize}

    Now we will find the form of the dual biderivation
      $\delta^\vee_{n-1}$ for the  \tm module $\Upsilon=\Upsilon_n$. By the linearity of the correspondence $\delta_{n-1}\longmapsto \delta_{n-1}^\vee$, it is enough to find the dual biderivations  for the generators of the following form:
     $$\delta_{n-1}=\Big[ d_{n-1}\tau^{h_{n-1}},\dots ,d_2\tau^{h_2},d_1\tau^{h_1} \Big]^T.$$ 
     Denote by  $V_n$ the space $\Big[ K\{\tau \}_{\langle 1, \rk\psi_n )}, \cdots, 
    K\{\tau \}_{\langle 1, \rk\psi_1 )} \Big]$. By \eqref{eq:iso:dowod} we have $\Upsilon^\vee\cong V_{n}$. Then we have the following basis for $V_n$:
     \begin{align}
    \label{eq:baza:dowod}
        \Big(E_n \tau^k \Big)_{k=1}^{\rk{\psi_n}-1},\cdots, 
        \Big(E_2 \tau^k \Big)_{k=1}^{\rk{\psi_2}-1}, \Big(E_1 \tau^k \Big)_{k=1}^{\rk{\psi_1}-1}
    \end{align}  
   Now we will find the form of the matrix  $D_{n,j}$ for fixed $j\in\{1,2,\dots, n-1\}$. 
   Again, $t*c\tau^k E_j$ for $1\leq k<r_j-1$ gives us  zero-value columns.

    Consider now the last column of
    $D_{n, j}$. The terms there come from the action of  $t$ on the generator $c\tau^{k}E_{j}$ for $k=r_j-1$. Then
    \begin{align*}
        t*c\tau^{r_j-1}E_{j}&= (\theta+\tau)c\tau^{r_j-1}E_{j}=\Big(\theta c\tau^{r_j-1} + c^{(1)}\tau^{r_j}\Big) E_{j}\notin V_n. 
    \end{align*}
    Since the degree of $\theta c\tau^{r_j-1} + c^{(1)}\tau^{r_j}$ equals 
    $\rk\psi_{j}$ we reduce it by means of the following inner biderivation:
    $$\delta^{( \mu\tau^0 E_{j})}_{\Upsilon, C}=\Big[\  \mu d_j\tau^{h_j}  \ \Big|\ \delta^{(\mu\tau^0E_j)}_{\Upsilon_{n-1},C}  \ \Big] \quad\textnormal{for}\quad 
    \mu =\dfrac{c^{(1)}}{a_{j, r_{j} }}.
    \footnote{Notice that in the above formula we used the symbol
    $E_j$ for two different matrices. In the left-hand side the matrix $E_j$ is of type $1\times n$ whereas in  the right-hand side it is of type $1\times (n-1)$. In both cases the only nonzero term is 1 at the  $j$-th position counted from the left. This should not lead to any ambiguity.  
    } $$
    
    As a result, one gets  the following equality:
    \begin{align}\label{eq:wartosc_dowod}
        t*c\tau^{r_j-1}E_{j}&= \Big[\  -\mu d_j\tau^{h_j}  \ \Big|\ \textnormal{terms for } \Upsilon_{n-1}^\vee  \ \Big]\in V_n,
    \end{align}
where the right-hand side has, by inductive hypotheses, the summands for
    $\Upsilon_{n-1}^\vee$.  Since $h_j< \rk\psi_{j}$, the value $t*c\tau^{r_j-1}E_{j}$ can be expressed in the basis \eqref{eq:baza:dowod}. We obtain
    
    \begin{align*}
     t*c\tau^{r_j-1}E_{j}&=\Bigg[\ 0,\dots, 0,\podwzorem{-\dfrac{d_j }{a_{j, r_{j} }} }{h_j-th \textnormal{ place}}\tau, 0, \dots, 
       0 \ \Big| \textnormal{ terms for } \Upsilon^\vee\ \Bigg]_{\big|_{\tau=c} }.
    \end{align*}
      Then the bloc $D_{n,j}=-\dfrac{d_j }{a_{j, r_{j} }}E_{h_j\times r_j-1}\tau$, 
    which shows \eqref{eq:postac_delta_vee_latwiejsza}.
\end{proof}

From the reasoning presented in the proof of Proposition \ref{prop:dual_postac}, we obtain the following corollary for any triangular \tm module whose quotient ranks are greater than $1$.

\begin{cor} \label{cor:postac_dual_w_ogolnosci}
    Let $\Upsilon=\Upsilon\big(\underline{\psi},\underline{\delta}\big)$ be a triangular \tm module of dimension $n$, with no nilpotence and $\rk\psi_{j}>1$ for $j=1,2,\dots,n$.
    Then  $\delta^\vee_{i}$ has the following bloc matrix form over $K(\Upsilon)$
    \begin{equation}
        \label{eq:delta_vee_block_matrix_form}
         \delta^\vee_{i}=\Big[ D_{i+1,i} \mid\dots \mid D_{i+1,2} \mid  D_{i+1,1}\Big],
    \end{equation}
    where $D_{i+1, j}$ for $j=1,2\dots, i$  are matrices of size $r_{i}-1\times r_{j}-1$, 
    in which only the last column has nonzero entries. 
     
\end{cor}

\section{Double dual \tm module and Cartier-Nishi theorem}

 Classical  Cartier–Nishi theorem asserts that one may recover
an abelian variety from its dual via the $\Ext$ functor.  
In the category of \tm modules an analogous result is known to hold for the 
duality of Drinfeld modules \cite{pr}, as well as for strictly pure \tm 
modules without nilpotence \cite{kk24}.  
We will show that the same phenomenon occurs for triangular \tm modules 
allowing for LD-biderivations.

\subsection{Double dual \tm module}
In the
 first step, we show the double dual analogue of Theorem \ref{thm:dual_dla_trangulowalnego_ogolnie}. 
We note that the double dual module
$\Big(\Upsilon^\vee\Big)^\vee$ will be denoted as ${\Upsilon^\vee}^\vee$, as this notation allows a more concise presentation.

In order to prove that the Cartier-Nishi theorem holds for triangular \tm 
modules allowing for LD-biderivations, we first need to determine the precise 
form of the double dual module for this class of \tm modules.  
We begin with the following lemma:

\begin{lem}\label{lem:ext_upsilon_dual_carlitz}
    Let $\Upsilon=\Upsilon\big(\underline{\psi},\underline{\delta}\big)$ be a triangular \tm module with no nilpotence, such that $r_l:=\rk\psi_l>1$ for 
    $l=1,2,\dots, n$.
    Then 
    \begin{align*}
    \quad(i)&\quad    \Ext_\tau^1(\Upsilon^\vee, C)\cong
    \Big[\podwzorem{K\{\tau\}_{<1},\dots,K\{\tau\}_{<1}, K\{\tau\}_{<2}}{\rk\psi_n-1 - \textnormal{terms} }\ \big|\ \cdots \\
    &\cdots\ \big|\ \podwzorem{K\{\tau\}_{<1},\dots,K\{\tau\}_{<1}, K\{\tau\}_{<2}}{\rk\psi_2-1 - \textnormal{terms} }\ 
    \big|\ \podwzorem{K\{\tau\}_{<1},\dots,K\{\tau\}_{<1}, K\{\tau\}_{<2}}{\rk\psi_1-1 - \textnormal{terms} }\Big] \\
    \quad(ii)&\quad    {\Upsilon^\vee}^\vee\cong
    \Big[\podwzorem{0,0,\dots, 0, K\tau}{\rk\psi_n-1 - \textnormal{terms}}\ \big|\ \cdots
    \ \big|\ \podwzorem{0,0,\dots, 0, K\tau}{\rk\psi_2-1 - \textnormal{terms}} 
    \ \big|\ \podwzorem{0,0,\dots, 0, K\tau}{\rk\psi_1-1 - \textnormal{terms}}  \Big]
        \end{align*}
\end{lem}

\begin{proof}
    Recall that $\Ext^1_\tau(\Upsilon^\vee, C)\cong \Der(\Upsilon^\vee, C)/\Derin(\Upsilon^\vee, C),$ where \linebreak $ \Der(\Upsilon^\vee, C)\cong \prod\limits_{l=1}^n\Mat_{1\times \rk\psi_l-1}\big(K\{\tau \}\big).$ 
    Every element of $\Der(\Upsilon^\vee, C)$ will  be identified  with the matrix consisting of $n$ blocs. 
    We note that in this proof, and in the remainder of the paper, blocs will always be numbered from right to left, whereas the elements within each block will be indexed in the usual manner from left to right.
    Denote by $E_{i\mid l}$ the matrix belonging to  $\Mat_{1\times \rk\psi_l-1}(K),$ which as the $i-$th entry (counted from left to right) has  $1$ and has zeros at all other places.
    Additionally,  by $\block{E_{i\mid l}}$ we denote the bloc matrix from $\prod\limits_{l=1}^n\Mat_{1\times \rk\psi_l-1}\big(K\{\tau \}\big),$
    where the $l-$th bloc from the right is equal to  $E_{i\mid l}$ and all the other blocs are equal to zero.
    Consider the following generators of the space
     $\Derin(\Upsilon^\vee, C)$ for $\mu\in K$ and $k=0,1,2,\dots$
    \begin{align}\label{eq:wewnetrzna_dla_duala}
        \delta^{\big(\block{E_{i\mid l}} \mu \tau ^k\big)}&= \mu\tau^k \block{E_{i\mid l}} \Upsilon^\vee - C \mu\tau^k \block{E_{i\mid l}} =\\
        &= \Big[\ \podwzorem{\boldsymbol{0}\ \big|\  \cdots\ \big|\  \boldsymbol{0}}{ (n-l)-\textnormal{blocks}}  \ \big|\ 
        \delta_{\psi_l^\vee, C}^{(E_{i\mid l}\mu\tau^k)}\ \big| 
        E_{i\mid l}\cdot \mu \tau ^k D_{l,l-1}\ \big|\ \cdots\ \big| E_{i\mid l}\cdot \mu \tau ^k D_{l,1}\Big],  \nonumber
        \end{align}
   where  $\delta_{\psi_l^\vee, C}^{(E_{i\mid l}\mu\tau^k)}\in\Derin(\psi_l^\vee, C)$ and 
   $\delta^\vee_{l-1}=\Big[ D_{l,l-1}\ \big|\ \dots\ \big|\ D_{l,2}\ \big|\   D_{l,1}\Big]$.

    From Proposition \ref{cor:dual_bez_nilpotentnosci} we get $\delta^{\big(\block{E_{i\mid l}} \mu \tau ^k\big)}\in\Der_0(\Upsilon^\vee, C)$ for every  $k=0,1,2,\dots$. Thus $(ii)$ follows from $(i)$.

    Notice that from  Corollary  \ref{cor:postac_dual_w_ogolnosci} it follows that in every matrix
    $$E_{i\mid l}\cdot  \mu \tau ^k D_{l,j}\quad \textnormal{for} \quad j=l-1,l-2, \dots, 2,1$$ only the last column is nonzero.
      
    Let
$\psi_l=\theta+\sum\limits_{s=1}^{r_l}a_{l,s}\tau^s$. Then from the form of \eqref{eq:postac_psi_dual} it follows that the inner biderivation $\delta_{\psi_l^\vee, C}^{(E_{i\mid l} \mu \tau^k)}$
    is one of  two types.
    For $i=1$ we obtain 
    \begin{align*}
        \delta_{\psi_l^\vee, C}^{(E_{1\mid l}\cdot\mu\tau^k)}&=
        \Big[ \mu\cdot \big(\theta^{(k)}-\theta\big)\tau^k - \mu^{(1)}\tau^{k+1}, 0,\dots, 0, -\dfrac{\mu\cdot a_{l,1}^{(k)}}{a_{l,r_l}^{(k)}}\tau^{k+1} + \dfrac \mu{a_{l,r_l}^{(k+1)}} \tau^{k+2}\Big]
    \end{align*}
    whereas for $i=2,\cdots, r_l-1$ one gets 
    \begin{align*}
        \delta_{\psi_l^\vee, C}^{(E_{i\mid l}\cdot \mu\tau^k)}&=
        \Big[\podwzorem{0,\dots, 0}{(i-2)-\textnormal{terms} } , \mu\cdot \tau^{k+1}, 
        \mu\cdot \big(\theta^{(k)}-\theta\big)\tau^k -\mu^{(1)}\tau^{k+1} , 0,\dots, 0 , -\dfrac{\mu\cdot a_{l,i}^{(k)}}{a_{l,r_l}^{(k)}}\tau^{k+1}
        \Big],
    \end{align*}
    where for $i=r_l-1$ the last term is of the following form:  
    $$\mu\cdot \big(\theta^{(k)}-\theta\big)\tau^k - \mu^{(1)}\tau^{k+1}  -\dfrac{\mu\cdot a_{l,i}^{(k)}}{a_{r_l}^{(k)}}\tau^{k+1}$$
   Thus, we have two types of biderivations $\delta^{\big(\block{E_{i\mid l}} \mu \tau ^k\big)}$ for the chosen $l-$th bloc. 
   The biderivation of the first type for  $i=1$, by means of which the last column (counted from the left) of the $l$-th bloc will be reduced. Then the biderivations of the second type for $i=2,3,\dots, \rk\psi_l-1$, 
   will be used for reduction of the
    $i-1$-th column in the  $l-$th bloc. As a result of such a reduction at the   $i$-th position,  a term of the same degree as that just reduced will appear. Therefore the biderivations  of the second type "shift" the highest term by one place to the right. Additionally, the new terms will appear only in the last columns of the blocks located to the right of the $l-$th bloc. 

Now we can explain the reduction of any biderivation from
    $\Der(\Upsilon^\vee, C)$. During the reduction process  the consecutive blocs from  left to  right will be changed. Assume that we have reduced all the blocs to the left of the $l-$th bloc. The reduction in the  $l$-th bloc will be performed according to the following procedure:
    \begin{enumerate}
        \item let $M$ be the  maximal degree 
        of terms in the $l-$th bloc
        \item denote by $m$ the index of the first from the left column with the term of degree  $M$. 
        \item if $M=2$ and $m$ is the index of the last column, then  we have finished the reduction of the bloc.
        \item if  $m$ is the index of the last column, then we reduce this term by means of the biderivation of the first type  $\delta^{\big(\block{E_{1\mid l}} \mu \tau ^k\big)}$ for appropriate  $\mu\in K$ and $k\in\mathbb{Z}_{\geq 0}.$ Then we go to step (1).
        \item if $m$ is not an index of the last column, then we reduce it by means of the biderivation of  the second type of the form
        $\delta^{\big(\block{E_{m+1\mid l}} \mu \tau ^k\big)}$ for appropriately chosen  $\mu\in K$ and $k\in\mathbb{Z}_{\geq 0}$. As  a result of this operation, a term of maximal degree is now in the $m+1-$st column. We return to step (3). 
    \end{enumerate}
    As a result of applying this procedure, we see  that all terms in the $l-$th bloc have expected degrees. Since
     $$\deg_\tau \delta_{\psi_l^\vee, C}^{(E_{1\mid l}\mu\tau^k)}= k+2,\quad 
    \deg_\tau \delta_{\psi_l^\vee, C}^{(E_{i\mid l}\mu\tau^k)}= k+1,\quad\textnormal{for}\quad i=2,3,\cdots \rk\psi_l-1,$$
    we see that two different bloc matrices in
     $$\Big[\podwzorem{K\{\tau\}_{<1},\dots,K\{\tau\}_{<1}, K\{\tau\}_{<2}}{\rk\psi_n-1 - \textnormal{terms} }\ \big|\ \cdots 
    \big|\ \podwzorem{K\{\tau\}_{<1},\dots,K\{\tau\}_{<1}, K\{\tau\}_{<2}}{\rk\psi_1-1 - \textnormal{terms} }\Big]$$
    represent different elements in  $\Ext_\tau^1(\Upsilon^\vee, C).$ 
    \end{proof}
As an immediate consequence of the preceding lemma, we obtain the following proposition.

\begin{prop}\label{prop:doubledualjesttmodulem}
Let $\Upsilon=\Upsilon_{\psi,\delta}$ be a triangular $t$-module with no
nilpotence such that $r_i:=\operatorname{rk}(\psi_i)>1$ for
$i=1,\ldots,n$. Then
$$\operatorname{Ext}^1_{\tau}(\Upsilon^\vee,\mathbb{C})
\quad
\textnormal{and}
\quad
{\Upsilon^\vee}^{\vee}
=
\operatorname{Ext}^0_{\tau}(\Upsilon^\vee,\mathbb{C})
$$
carry natural structures of $t$-modules.
\end{prop}

\begin{proof}
By the preceding lemma, the basis constructed in its proof satisfies
\textbf{Step 2} of the \tm reduction algorithm. Therefore, the conclusion follows
immediately from Proposition \ref{prop:step2towystarczy}.
\end{proof}

We are now ready to prove the dual version of Theorem \ref{thm:dual_dla_trangulowalnego_ogolnie}.

\begin{thm}\label{thm:doubledual_dla_trangulowalnego_ogolnie}
  Let $\Upsilon=\Upsilon\big(\underline{\psi},\underline{\delta}\big)$ be a triangular \tm module 
  with no nilpotence such that  $r_i:=\rk\psi_i>1$ for $i=1,\dots, n$.
  \\ Then 
  \begin{itemize}
      \item[(i)] double dual ${\Upsilon^\vee}^\vee$ is a  \tm module given by the $\tau-$composition series of the following form: 
      $${\Upsilon^\vee}^\vee_1 \hookrightarrow {\Upsilon^\vee}^\vee_2 \hookrightarrow \cdots 
\hookrightarrow {\Upsilon^\vee}^\vee_n={\Upsilon^\vee}^\vee,$$
    \item[(ii)] there exists the following exact sequence of  \tm modules  
        $$0\lra {\Upsilon^\vee}^\vee\lra \Ext^1_\tau(\Upsilon^\vee, C)\lra \G_a^s\lra 0$$
        where $s=\sum\limits_{l=1}^{n}\rk\psi_l -n$
      \item[(iii)] Every map $f:\Upsilon\lra \underline{\Upsilon}$,  of \tm modules that
      satisfy the assumptions of this theorem, induces a map
       ${f^\vee}^\vee:{ {\Upsilon}^\vee} ^\vee\lra {\underline{\Upsilon}^\vee}^\vee$.
  \end{itemize}
\end{thm}
\begin{proof}
By Theorem~\ref{thm:dual_dla_trangulowalnego_ogolnie}, there exist short exact sequences
$$\delta_l^\vee: \quad 0\lra \psi_{l+1}^\vee\lra \Upsilon_{l+1}^\vee \lra \Upsilon_{l}^\vee\lra 0,$$ 
for $l=1,2,\dots,n$. Since Proposition~\ref{prop:hom_dual_Carlitz} implies that
$\Hom_\tau(\psi_{l+1}^\vee,C)=0.$
Applying the functor $\Hom_\tau(-,C)$ to the exact sequence $\delta_l^\vee$ yields the following short exact sequence of $\F_q[t]$-modules:
\begin{align*}
0\lra \Ext^1_\tau(\psi_{l}^\vee\,C)\lra \Ext^1_\tau(\Upsilon_{l+1}^\vee,C) \lra \Ext^1_\tau(\psi_{l+1}^\vee,C)\lra 0.   
\end{align*}
By Proposition~\ref{prop:doubledualjesttmodulem}, this is a short exact sequence of \tm modules.
 The existence of the exact sequences of \tm modules (see \cite[Theorem 1.1]{pr}) 
\begin{align*}
  0\lra {\psi_l^\vee}^\vee\lra \Ext^1_\tau(\psi^\vee_l,C)\lra \G_a^{\rk\psi_l-1}\lra 0, \quad\textnormal{for}\quad l=1,2 
\end{align*}
implies that there exists the following commutative diagram of \tm modules with exact rows: 
$$\xymatrix{
& 0 \ar[d] & 0 \ar[d] & 0 \ar[d] & \\ 
0\ar[r] & {\psi_1^\vee}^\vee\ar[r] \ar[d] & {\Upsilon_2^\vee}^\vee \ar[r]\ar[d] & {\psi_2^\vee}^\vee \ar[r] \ar[d] & 0 \\
0\ar[r]& \Ext^1_\tau(\psi_1^\vee,C) \ar[r] \ar[d]& \Ext^1_\tau(\Upsilon_2^\vee,C) \ar[r] \ar[d]& \Ext^1_\tau(\psi_2^\vee,C) \ar[r]\ar[d]& 0\\
0\ar[r] & \G_a^{\rk\psi_1-1} \ar[r]\ar[d] & \G_a^{\rk\psi_2-1}\times \G_a^{\rk\psi_1-1}\ar[r]\ar[d] & \G_a^{\rk\psi_2-1} \ar[r] \ar[d] & 0\\
& 0 & 0 & 0 & \\},$$
where the surjectivity of the middle vertical map follows from the Snake Lemma.
Therefore, we obtain the following exact sequence:
$$0\lra {\Upsilon_2^\vee}^\vee\lra \Ext_\tau^1(\Upsilon_2^\vee, C)\lra \G_a^{\rk\psi_1+\rk\psi_2-2}\lra 0.$$
Repeating this reasoning for the exact sequences  $\delta_l^\vee$ for $l=2,3,\dots, n-1$ we obtain the fact that  there are the following exact sequences of \tm modules:
$${\delta_l^\vee}^\vee: \quad 0\lra {\Upsilon_{l}^\vee}^\vee\lra {(\Upsilon_{l+1} )^\vee}^\vee\lra {(\psi_{l+1})^\vee}^\vee\lra 0$$ 
for $l=1,2,\dots, n-1$.  
 From these exact sequences, in the same way as before, we obtain the following exact sequences of \tm modules 
$$0\lra {\Upsilon_{l}^\vee}^\vee\lra \Ext_\tau^1(\Upsilon_{l}, C)\lra \G_a^{s}\lra 0,$$
where $s=\sum\limits_{i=1}^l\rk\psi_i-l$ for $l=2,\dots, n$. The last sequence for  $l=n$ shows $(ii)$.  

From \cite[Theorem 1.1]{pr} ${\psi_l^\vee}^\vee\cong \psi$  are $\tau-$simple. Therefore, the sequences  ${\delta_l^\vee}^\vee$ for $l=1,2,\dots, n-1$ show $(i)$.

Part $(iii)$. From $(ii)$ there exist exact sequences of the following form:
\begin{align*}
    0\lra &{\Upsilon^\vee}^\vee\lra \Ext^1_\tau(\Upsilon^\vee, C)\lra \G_a^{s_1}\lra 0\\
0\lra &\Big(\widehat{\Upsilon}^\vee\Big)^\vee\lra \Ext^1_\tau(\widehat{\Upsilon}^\vee, C)\lra \G_a^{s_2}\lra 0.
\end{align*}

Arguing similarly as in the proof of 
 \cite[Proposition 3.2]{kk24}, we obtain that there exists  the following commutative diagram with the exact rows:
$$\xymatrix{0\ar[r] & {\Upsilon^\vee}^\vee \ar[d]^{ f^\vee} \ar[r]  &\Ext^1_\tau( \Upsilon^\vee, C) \ar[d]^{\overline{f} }  \ar[r] &\G_a^{s_1} \ar[d]^{\partial f}  \ar[r] & 0\\
0\ar[r] &\Big(\widehat{\Upsilon}^\vee\Big)^\vee\ar[r] &  \ar[r] \Ext^1_\tau(\widehat{\Upsilon}^\vee, C) & \G_a^{s_2} \ar[r] &  0 }.$$
\end{proof}

\begin{remark}
Our goal is to obtain the Cartier–Nishi theorem in a form analogous to that obtained
for Drinfeld modules \cite[Theorem~1.1]{pr}, and  also for strictly pure 
\tm modules \cite[Theorem~1.2]{kk24}.  
By Theorem~\ref{thm:doubledual_dla_trangulowalnego_ogolnie}, it suffices to 
show that ${\Upsilon^\vee}^\vee \cong \Upsilon$.
\end{remark}

\begin{remark}\label{rem:postac_macierzowa_double_duala}
From the previous theorem it follows that the \tm module 
${\Upsilon^\vee}^\vee$ has the following matrix form:
\begin{equation*}
         {\Upsilon^\vee}^\vee_t=\begin{pNiceArray}{cccccc}[margin=15pt]
 \Block[draw=black]{1-1}<\large>{{\psi_n^\vee}^\vee}& 0&\cdots&   0&0 &0\\[2pt] 
         \Block[draw=black,fill=gray!20!]{5-1}<>{{(\delta_{n-1})^\vee}^\vee}&\Block[draw=black]{1-1}<>{ { {(\psi_{n-1})}^\vee}^\vee  }&\cdots&0 &0 &0\\[2pt]
         & \Block[draw=black,fill=gray!20!]{4-1}<>{{(\delta_{n-2})^\vee}^\vee} &\ddots &\vdots &\vdots &\vdots \\[2pt] 
         & & & \Block[draw=black]{1-1}<>{{\psi_3^\vee}^\vee }&0 &0 \\[2pt] 
         & &\cdots & \Block[draw=black,fill=gray!20!]{2-1}<>{{\delta_{2}^\vee}^\vee }& \Block[draw=black]{1-1}<>{{\psi_2^\vee}^\vee}&0 \\[2pt] 
         & & & & \Block[draw=black,fill=gray!20!]{1-1}<>{ {\delta_{1}^\vee}^\vee }& \Block[draw=black]{1-1}<>{{\psi_1^\vee}^\vee}
\end{pNiceArray},
    \end{equation*}
Recall that from the proof of \cite[Theorem 3.4]{kk24}, it follows that
$${\psi_l^\vee}^\vee = \dfrac 1{a_{l,r_l}}\cdot \psi_l\cdot  a_{l,r_l}\quad \textnormal{for} 
\quad l=1,2,\dots, n.$$ 
 Thus the module
${\Upsilon^\vee}^\vee$ has  sub-quotients of the form  
$\dfrac 1{a_{l,r_l}}\cdot \psi_l\cdot  a_{l,r_l}$. In order to return to the initial sub-quotients  $\psi_l$, we will consider  the conjugate of ${\Upsilon^\vee}^\vee$ by the diagonal matrix $D =  \diag\big(a_{n,r_n},\cdots, a_{1,r_{1}}\big)$, where $a_{j,r_j}$ is the leading term of $\psi_j$:
$$ {}^D {\Upsilon^\vee}^\vee :=\diag\big(a_{n,r_n},\cdots, a_{1,r_{1}}\big) \cdot {\Upsilon^\vee}^\vee\cdot  \diag\big(a_{n,r_n},\cdots, a_{1,r_{1}}\big)^{-1}.$$
Then the biderivations for 
 ${}^D {\Upsilon^\vee}^\vee$ are given by the following formulas:
$${}^D {\delta_l^\vee}^\vee=\diag\big(a_{l,r_l},\cdots, a_{1,r_{1}}\big) \cdot {\delta_l^\vee}^\vee \cdot \dfrac 1{a_{l+1,r_{l+1}}}\quad \textnormal{for} \quad l=1,2,\dots, n-1.$$
In the case of $t$-modules allowing for  LD-biderivations we will show that
\[
\delta_l = {}^D {\delta_l^\vee}^\vee
\]
from which it follows immediately that 
$\Upsilon \cong {\Upsilon^\vee}^\vee.$
\end{remark}

\subsection{The case of LD-biderivations.}
We consider triangular \tm modules that allow for LD-biderivations.

Recall the procedure for determining  
 ${\psi^\vee}^\vee$ for a Drinfeld module of the form
$\psi=\theta+\sum\limits_{s=1}^{r}a_{s}\tau^s$. Since ${\psi^\vee}^\vee\cong \big[0,0,\cdots, 0, K\tau\big]$,  in order to determine the  \tm module structure on ${\psi^\vee}^\vee$, we have to find the result of the action  
$t*\big[0,0,\cdots, 0, c\tau\big]=\big[0,0,\cdots, 0,\theta\tau+ c^{(1)}\tau^2\big]$. 
For this we reduce by means of the following inner biderivation:
\begin{align}\label{eq:biderywacje_wewnetrzne_cykl}
    \delta_{\psi^\vee,C}^{\big( V_{\psi}(\mu)\big)}&=
    \Bigg[ 0, \dots, 0, \dfrac{\mu}{a_{r}^{(1)} }\tau^2 - \sum\limits_{s=1}^{r}\dfrac {a_{s}}{a_{r} } \mu^{(s-1)}\cdot\tau\Bigg]\quad \textnormal{where}\\
    &\qquad \mu=c^{(1)}a_r^{(1)}
    \quad \textnormal{and} \quad V_{\psi_l}(\mu):=\big[ \mu, \mu^{(1)},\cdots, \mu^{(r-2)} \big]\tau^0. \nonumber 
\end{align}
Then we write the result in the standard basis 
and obtain 
the form of 
${\psi^\vee}^\vee$. The inner biderivation $\delta_{\psi^\vee,C}^{\big( V_{\psi}(\mu)\big)}$ will be used further on to determine ${\Upsilon^\vee}^\vee$ for a triangular \tm module $\Upsilon$.

\begin{thm}[Cartier-Nishi Theorem]\label{thm:dual_dual_case_1}
  Let $\Upsilon=\Upsilon\big(\underline{\psi},\underline{\delta}\big)$ be a triangular \tm module with no nilpotence, such that $r_l:=\rk\psi_l>1$ for 
    $l=1,2,\dots, n$. If $\Upsilon$ allows for LD-biderivations, then
${\Upsilon^\vee}^\vee\cong \Upsilon$ over $K(\Upsilon)$, and 
there exists the following exact sequence of  \tm modules over $K(\Upsilon)$ 
        $$0\lra \Upsilon \lra \Ext^1_\tau(\Upsilon^\vee, C)\lra \G_a^s\lra 0$$
        where $s=\sum\limits_{l=1}^{n}\rk\psi_l -n$
\end{thm}	
\begin{proof}
Since $\Upsilon$ admits LD-biderivations, after changing the field to
$K(\Upsilon)$, we may assume that $\Upsilon$ satisfies
condition~\eqref{eq:warunek_LD}.

The existence of this exact sequence will follow from the isomorphism 
${\Upsilon^\vee}^\vee\cong \Upsilon$ applied to Theorem   \ref{thm:doubledual_dla_trangulowalnego_ogolnie}. So we will focus on establishing this isomorphism. 
According to Remark \ref{rem:postac_macierzowa_double_duala}, it is enough to show that
$${}^D {\delta_l^\vee}^\vee= \delta_l  \quad \textnormal{for}\quad l=1,2,\dots, n-1.$$
Let $\psi_l=\theta+\sum\limits_{s=1}^{r_l}a_{l,s}\tau^s$. 
We will induct on the number of biderivations. 

Let
    $\Upsilon_2=\left[\begin{array}{c|c}
         \psi_2&0  \\ 
         \hline
         \delta_1&\psi_1 
    \end{array}\right]$. From Corollary \ref{cor:dual_dla_biderwaycji} it follows that the correspondence  $\delta_1\longmapsto {\delta_1^\vee}^\vee$ is $\F_q-$linear. Therefore it is enough to find the double dual biderivation for generators of the form $d\tau^h$ for $1\leq h< r_2$ and to show that 
    $$a_{1,r_1}\cdot{\big( d\tau^h\big)^\vee}^\vee\cdot\dfrac 1{a_{2,r_2} } =d\tau^h.$$
    Recall that
     $$\big(d\tau^h\big)^\vee = -\dfrac{d }{a_{1, r_{1} }}E_{h\times r_1-1} \tau,$$

    From Lemma
 \ref{lem:ext_upsilon_dual_carlitz} it follows that  
$${\Upsilon_2^\vee}^\vee\cong
    \Big[\podwzorem{0,0,\dots, 0, K\tau}{r_2-1 - \textnormal{terms}} 
    \ \big|\ \podwzorem{0,0,\dots, 0, K\tau}{r_1-1 - \textnormal{terms}}  \Big].$$
      Recall that $E_{i\mid l}$ denotes the matrix of dimensions $1\times r_l-1,$ which at the $i-$th entry  has value $1$  (counted from left to right) and in addition zeros at all other entries, and 
      $$\block{E_{i\mid 2}}=\Big[\ E_{i\mid 2}  \  \big|\ \boldsymbol{0}\ \Big],\quad 
      \block{E_{i\mid 1}}=\Big[ \ \boldsymbol{0}\ \big| \ E_{i\mid 1}  \  \Big].$$ 
    In order to determine the $F_q[t]-$module structure on ${\Upsilon_2^\vee}^\vee$, 
    one has to find an action of  $t$ on the generator $c\cdot \tau \block{E_{r_2-1\mid 2} }  $ written in the following basis:
\begin{equation*}
    \tau\cdot\block{E_{r_2-1\mid 2} }=\Big[\ 0,\dots, 0, \tau \ \big|\ \boldsymbol{0}\  \Big], \qquad \tau\cdot \block{E_{r_1-1\mid 1} }=\Big[\ \boldsymbol{0} \ \big|\ 0,\dots, 0,\tau \  \Big].
\end{equation*}
Then
\begin{align*}
    t*c\tau \cdot \block{E_{r_2-1\mid 2} } = C_t\cdot  c \tau\cdot \block{E_{r_2-1\mid 2} } = \Big[ \ 0,\dots, 0, \theta c\tau+c^{(1)}\tau^2\  \big| \ \boldsymbol{0} \Big].
\end{align*}
The summand $c^{(1)}\tau^2\cdot \block{E_{r_2-1\mid 2} }$ is reduced by means of the following inner biderivation given by equation \eqref{eq:biderywacje_wewnetrzne_cykl}:
\begin{align*}
    \delta^{\big( \block{V_{\psi_2}(\mu)} \big)}_{\Upsilon^\vee_2,C}= \Bigg[ \delta^{\big( V_{\psi_2}(\mu) \big)}_{\psi_2^\vee,C} \Big|\ 0,\dots,0, -\dfrac{\mu^{(h-1)} d}{a_{1,r_1}} \tau \Bigg]
\end{align*}
where $\block{V_{\psi_2}(\mu)}:=\big[ V_{\psi_2}(\mu) \big|\ 0,\cdots, 0\  \big]$
for $\mu= c^{(1)}\cdot a_{2,r_2}^{(1)}$.
    
We obtain the following equality: 
\begin{align*}
    t*c\tau \cdot \block{E_{r_2-1\mid 2} } 
    &= \Bigg[\textnormal{ terms for } {\psi_2^\vee}^\vee\ \big|\ 0,\dots, 0,  \dfrac{\mu^{(h-1)} d}{a_{1,r_1}} \tau  \Bigg].
\end{align*}
Substituting $\mu=c^{(1)}a_{r_2,2}^{(1)}$  and expressing $ t*c\tau \cdot \block{E_{r_2-1\mid 2} }$ in the basis one gets
\begin{align*}
    t*c\tau \cdot \block{E_{r_2-1\mid 2} } &= \Bigg[\  {\psi_2^\vee}^\vee\ ,\  \dfrac{a_{2,r_2}^{(h)} d}{a_{1,r_1}} \tau^h  \Bigg]_{\Big|_{\tau=c}}.
\end{align*}
Hence ${\big( d\tau^h\big)^\vee}^\vee= \dfrac{a_{2,r_2}^{(h)} d}{a_{1,r_1}} \tau^h=\dfrac{ d}{a_{1,r_1}} \tau^h\cdot a_{2,r_2}$
and
$a_{1,r_1}\cdot{\big( d\tau^h\big)^\vee}^\vee\cdot\dfrac 1{a_{2,r_2} } =
d\tau^h.$

Assume that for a triangular \tm module of dimension $n-1$ we have
$${}^D {\delta_l^\vee}^\vee= \delta_l \quad \textnormal{for}\quad l=1,2,\dots, n-2.$$
 and every double dual biderivation ${\delta_l^\vee}^\vee$ came from determining 
 the action
 $t*c\tau \cdot \block{E_{r_l-1\mid l} }$, where the following inner  biderivation 
$$\delta^{\big( \block{V_{\psi_l}}(\mu) \big)},\quad\textnormal{for}\quad  \mu= c^{(1)}\cdot a_{l,r_l}^{(1)}$$
was used in the reduction process. 

Here, $\block{V_{\psi_l}(\mu) }$ is a bloc matrix where the only nonzero bloc is the $l-$th, counted from right to left, and is equal to $V_{\psi_l}(\mu)$.

Let $\Upsilon$ be a triangular  \tm module of dimension  $n$. By the inductive hypothesis
$${}^D {\delta_{l}^\vee}^\vee= \delta_{l}\quad \textnormal{for}\quad l=1,2,\cdots n-2.$$
Since the assignment $\delta_{n-1}\longmapsto {(\delta_{n-1})^\vee}^\vee$ is $\F_q-$linear,  it is enough to consider the case where  
$\delta_{n-1}=\Big[ d_{n-1}\tau^{h_{n-1}},\cdots,  d_{2}\tau^{h_{2}},d_{1}\tau^{h_{1}}    \Big]^T.$
Recall that by \eqref{eq:postac_delta_vee_latwiejsza} we have
     $$\delta_{n-1}^\vee = -\sum\limits_{j=1}^{n-1} \dfrac{d_i }{a_{j, r_{j} }}E_{h_i\times r_j-1} \tau.$$

From Lemma
 \ref{lem:ext_upsilon_dual_carlitz} it follows that  
$${\Upsilon^\vee}^\vee\cong
    \Big[\podwzorem{0,0,\dots, 0, K\tau}{r_n-1 - \textnormal{terms}}\ \big|\ \cdots
    \ \big|\ \podwzorem{0,0,\dots, 0, K\tau}{r_2-1 - \textnormal{terms}} 
    \ \big|\ \podwzorem{0,0,\dots, 0, K\tau}{r_1-1 - \textnormal{terms}}  \Big].$$
    Recall that $E_{i\mid l}$ denotes the matrix of dimensions $1\times r_{l}-1,$ which at the $i-$th entry  has  value $1$,  zeros at all other entries, and $\block{E_{i\mid l}}$ 
    denotes the bloc matrix where the $l-$th bloc (counted from the right) equals  $E_{i\mid l}$ and all the other blocs are zero.  
    In order to determine the $F_q[t]-$module structure on ${\Upsilon^\vee}^\vee$, 
    one has to find an action of  $t$ on the generators $c\cdot \tau \block{E_{r_{n}\mid n} }  $ written in the following basis:
\begin{equation}\label{eq:baza_przypadek_latwy}
   \tau \cdot\block{E_{r_n\mid n} },\ \tau\cdot\block{E_{r_{n-1}\mid n-1} },\  \cdots,\  
    \tau\cdot\block{E_{r_2\mid 2} },\ \tau\cdot \block{E_{r_1\mid 1} }.
\end{equation}
Then 
\begin{align*}
    t*c\tau \cdot \block{E_{r_{n}-1\mid n} }  = (\theta + \tau) c\tau\cdot \block{E_{r_{n}-1\mid n} } = 
    \Big( \theta c\tau+c^{(1)}\tau^2\Big)\cdot \block{E_{r_{n}-1\mid n} }
\end{align*}
The summand $c^{(1)}\tau^2\cdot \block{E_{r_{n}-1\mid n} }$ is reduced by means of the following inner biderivation given by equation \eqref{eq:biderywacje_wewnetrzne_cykl}:
\begin{align*}
    \delta^{\big(\block{V_{\psi_n}  (\mu) } \tau ^0\big)} &=
\Bigg[  \delta^{\big(V_{\psi_n}(\mu)\big)}_{\psi_n,C}\ \big|  \ 0,\dots, 0,  -\dfrac{\mu^{(h_{n-1}-1)} d_{n-1}}{a_{n-1,r_{n-1}}} \tau\ \big|\cdots 
    \\
    &\cdots \big| \ 0,\dots, 0,  -\dfrac{\mu^{(h_{2}-1)} d_{2}}{a_{2,r_{2}}} \tau\ \big|
    \ 0,\dots, 0,  -\dfrac{\mu^{(h_{1}-1)} d_{1}}{a_{1,r_{1}}} \tau\ \Bigg]
\end{align*}
for $\mu= c^{(1)}\cdot a_{n,r_{n}}^{(1)}$. 
We obtain 
\begin{align*}
    t*c\tau \cdot \block{E_{r_n-1\mid n} } = 
    \Bigg[&\ \textnormal{ terms for } {\psi_2^\vee}^\vee\ \big| \ 0,\dots, 0,  \dfrac{\mu^{(h_{n-1}-1)} d_{n-1}}{a_{n-1,r_{n-1}}} \tau\ \big|\cdots 
    \\
    &\cdots \big| \ 0,\dots, 0,  \dfrac{\mu^{(h_{2}-1)} d_{2}}{a_{2,r_{2}}} \tau\ \big|
    \ 0,\dots, 0,  \dfrac{\mu^{(h_{1}-1)} d_{1}}{a_{1,r_{1}}} \tau\ \Bigg]
\end{align*}
Substituting $\mu=c^{(1)}a_{r_n,n}^{(1)}$ and expressing $ t*c\tau \cdot \block{E_{r_n-1\mid n} }$ in the basis  \eqref{eq:baza_przypadek_latwy} 
one gets
\begin{align*}
    t*c\tau \cdot \block{E_{r_n-1\mid n} } = 
    \Bigg[&\ \textnormal{ terms for } {\psi_2^\vee}^\vee,\ \dfrac{a_{r_n,n}^{(h_{n-1})} d_{n-1}}{a_{n-1,r_{n-1}}} \tau^{h_{n-1}}, \cdots,\\
    &\qquad\qquad\qquad\qquad \dfrac{a_{r_n,n}^{(h_2)} d_{2}}{a_{2,r_{2}}} \tau^{h_2},
    \dfrac{ a_{r_n,n}^{(h_1)} d_{1}}{a_{1,r_{1}}} \tau^{h_1}\ \Bigg]_{\big|_{\tau=c}}.
\end{align*}
Then $${(\delta_{n-1})^\vee}^\vee = \Bigg[ \dfrac{a_{r_n,n}^{(h_{n-1})} d_{n-1}}{a_{n-1,r_{n-1}}} \tau^{h_{n-1}}, \cdots,\dfrac{a_{r_n,n}^{(h_2)} d_{2}}{a_{2,r_{2}}} \tau^{h_2},\ \dfrac{ a_{r_n,n}^{(h_1)} d_{1}}{a_{1,r_{1}}} \tau^{h_1}\Bigg]^T.$$ 

So,
$${}^D{(\delta_{n-1})^\vee}^\vee =\diag\big(a_{n-1,r_{n-1}},\cdots, a_{1,r_{1}}\big) \cdot {(\delta_{n-1})^\vee}^\vee \cdot \dfrac 1{a_{n,r_{n}}}= \delta_{n-1}.$$

\end{proof}

As an immediate consequence of Theorem~\ref{thm:dual_dual_case_1} and 
Corollary~\ref{cor:male_stopnie_delta}, we obtain the following fact.

\begin{cor}\label{cor:Cartier_Nishi_1}
    Let $\Upsilon=\big(\underline{\psi},\underline{\delta}\big)$ be a triangular  \tm module with no nilpotence, such that
    $$1< \rk\psi_1< \rk\psi_2< \cdots<\rk\psi_n.$$
    Then 
${\Upsilon^\vee}^\vee\cong \Upsilon$ over $K(\Upsilon)=K$, and 
there exists the following exact sequence of  \tm modules over $K$
        $$0\lra \Upsilon \lra \Ext^1_\tau(\Upsilon^\vee, C)\lra \G_a^s\lra 0$$
        where $s=\sum\limits_{l=1}^{n}\rk\psi_l -n$
        \qed
\end{cor}

Similarly, from Corollary~\ref{cor:male_stopnie_delta_z_rownosciami} and 
Theorem~\ref{thm:dual_dual_case_1} we obtain the following corollary.
\begin{cor}\label{cor:Cartier_Nishi_2}
    Let $\Upsilon=\big(\underline{\psi},\underline{\delta}\big)$ be a triangular  \tm module, defined over $K$, with no nilpotence, such that
    $$1< \rk\psi_1\leq \rk\psi_2\leq \cdots \leq\rk\psi_n.$$
    Then ${\Upsilon^\vee}^\vee\cong \Upsilon$ over $K(\Upsilon)$ and 
there exists the following exact sequence of  \tm modules  over $K(\Upsilon)$
        $$0\lra \Upsilon \lra \Ext^1_\tau(\Upsilon^\vee, C)\lra \G_a^s\lra 0$$
        where $s=\sum\limits_{l=1}^{n}\rk\psi_l -n$. \qed
\end{cor}

\begin{remark}
 Notice that in order to define the dual $t$-module $\Phi^{\vee}$ using the
methods of Taguchi from \cite{ta}, the assumption that $\Phi$ is strictly pure is essential.
Indeed, in \cite{kk24} it was shown that for a strictly pure $t$-module $\Phi$
without nilpotence, the Cartier--Taguchi dual module coincides with that defined
as $\Ext_{\tau,0}(\Phi,C)$.  However, triangular $t$-modules are in general
not strictly pure (cf.\ Proposition~\ref{prop:strictly_pure_characterisation}).
This indicates that the construction of the dual module via
$\Ext_{\tau,0}^1(\Phi,C)$ is in fact more general.
\end{remark}

The following example shows that in  Theorem \ref{thm:dual_dual_case_1} the nilpotence condition is necessary. 
\begin{ex}\label{ex:dualdual_nie_dziala_z_nilpotentnoscia}
    Consider the  \tm module $$\Upsilon = \left(\begin{array}{c|c}
        \theta+ \tau^3  & 0 \\ \hline
         1              & \theta+  \tau^2
    \end{array}\right).$$
Then
$$\Upsilon^\vee = \left(\begin{array}{ccc|c}
        \theta & \tau^2  & 0 & \tau \\
        \tau   & \theta  & 0 & 0    \\
        0 & 0 & \theta & \big(\theta-\theta^{(1)} \big)\tau \\ \hline
        0 & 0 & \tau & \theta +\tau^2
    \end{array}\right).$$
For determination of ${\Upsilon^\vee}^\vee$ we have to assume that $K$ is a perfect field. Then from the reduction algorithm we obtain that
$${\Upsilon^\vee}^\vee=
\left(\begin{array}{ccc}
        \theta+ \tau^3  & 0 & 0\\ 
         \tau  & \theta &  \big(\theta ^{(1)}-\theta \big)+\tau^2 \\
         0 & \tau & \theta
    \end{array}\right)\ncong \Upsilon.$$
\end{ex}

Our results show that the Cartier--Nishi theorem holds for triangular 
\tm modules allowing for LD-biderivations, and we have provided an explicit 
sufficient condition ensuring that a triangular \tm module allows for LD- biderivations.  
At present, it is natural to ask to what extent this hypothesis is actually 
necessary.  All triangular \tm modules we have tested outside those allowing for LD-biderivations still satisfy the Cartier--Nishi formula, and no counterexample is currently 
known.  This leads to the following question:

\begin{question}
For which triangular \tm modules does the Cartier--Nishi theorem hold?
Does it hold for all triangular \tm modules, or does there exist a triangular
\tm module for which the theorem does not hold?
\end{question}

The situation becomes considerably more complicated beyond the class $LD$, as illustrated by the following example.
\begin{example}
Let 
$$\Upsilon=\left(
\begin{array}{ccc}
 \theta +\tau ^6 & 0 & 0 \\
 a \tau ^9 & \theta +\tau ^4 & 0 \\
 b \tau ^6 & c \tau ^7 & \theta +\tau ^5 \\
\end{array}.
\right)$$
Then, after computing the double dual module, we obtain
 $${\Upsilon^\vee}^\vee=\left(
\begin{array}{ccc}
 \theta +\tau ^6 & 0 & 0 \\
 \delta_{2\times 3} & \theta +\tau ^4 & 0 \\
 \delta_{1\times 3} & 
   \delta_{1\times 2} & \theta +\tau ^5 \\
\end{array}
\right),$$
where
\begin{align*}
     \delta_{2\times 3}&=a^{(4)}\tau ^7 +a\big(\theta -\theta ^{(3)}\big)  \tau ^3,\qquad
    \\
     \delta_{1\times 3}&= -a^{(13)} c^{(10)} \tau ^{10}+ 
     a^{(4)} c^{(5)}  \left( \theta ^{(7)}-\theta^{(4)}\right)  \tau ^7
     +a^{(3)} c \left(\theta  - \theta ^{(3)}\right) \tau ^6 \\
     &+ \left(a^{(8)} c^{(5)}  \left(\theta ^{(11)}-\theta \right)+b^{(5)}\right)\tau ^5\\
      \delta_{1\times 2}&=  
   c^{(10)} \tau ^9+ c^{(5)} \left(\theta -\theta ^{(4)}\right) \tau ^4+c  \left(\theta -\theta ^{(3)}\right)\tau ^3 
\end{align*}
Then $\Upsilon= U\cdot {\Upsilon^\vee}^\vee \cdot V$, where 
$$U=\left(
\begin{array}{ccc}
 1 & 0 & 0 \\
 a \tau ^3 & 1 & 0 \\
 -a^{(8)} c^{(5)} \tau ^5  -ca^{(3)}\big(\theta  +  \theta
   ^{(6)}\big)+b & c^{(5)} \tau ^4 +c \tau ^3 & 1 \\
\end{array}
\right) $$
and
$$V=\left(
\begin{array}{ccc}
 1 & 0 & 0 \\
- a \tau ^3 & 1 & 0 \\
 a^{(4)} c^{(5)} \tau ^7 + a^{(3)} c  \tau ^6 
 +a^{(8)} c^{(5)} \tau ^5  +ca^{(3)}\big(\theta  +  \theta
   ^{(6)}\big)-b & -c^{(5)} \tau ^4 - c \tau ^3 & 1 \\
\end{array}
\right)$$
Thus, for the $t$-module $\Upsilon$, the Cartier--Nishi theorem holds, even though $\Upsilon$ is not a $t$-module that allows for $\mathrm{LD}$-biderivations. It is worth noting, however, that the isomorphism between $\Upsilon$ and ${\Upsilon^\vee}^\vee$ is considerably more intricate than the one obtained in the $\mathrm{LD}$ case.
\end{example}

Moreover, all the examples we have computed indicate that, outside the class $\mathrm{LD}$, the isomorphisms between a triangular $t$-module and its double dual exhibit a considerably more complicated structure. This complexity increases rapidly with the dimension of the $t$-module and with the degrees of the underlying Drinfeld modules and biderivations.

\appendix

{}

\end{document}